\documentclass[preprint,12pt,nopreprintline]{elsarticle}
\usepackage{amssymb}
\usepackage{amsmath}
\usepackage{amsthm}
\usepackage{textcomp}
\usepackage{mathcomp}
\usepackage{adjustbox}
\usepackage{enumitem}
\usepackage{bm}
\usepackage{cases}
\usepackage{gensymb}
\usepackage{breqn}
\usepackage{lineno}
\usepackage{color}
\usepackage{xcolor}
\usepackage[colorlinks=true]{hyperref}
\usepackage{easyReview}
\DeclareUnicodeCharacter{2212}{-}
\usepackage[most]{tcolorbox}
\usepackage{subcaption}
\usepackage{here}
\usepackage{physics2}
\usephysicsmodule{ab}
\usepackage{booktabs}
\newcommand{\pdv}[3][]{\cfrac{\partial^{#1} {#2}}{\partial {#3}^{#1}}} % 偏微分
\newcommand{\mdv}[3]{\cfrac{\partial^2 {#1}}{\partial {#2}\partial {#3}}} %混合微分
\newcommand{\lap}[1]{\nabla^2 {#1}} %ラプラシアン
\newcommand{\sphsum}[2][j]{\sum_{{#1}\in\mathbb{S}_{i}}{#2}}

\newcommand{\wh}[1]{\widehat{#1}}
\newcommand{\wt}[1]{\widetilde{#1}}
\newcommand{\mr}[1]{\mathrm{#1}}

\newcommand{\refF}[1]{\hyperref[#1]{Fig.~\ref{#1}}}
\newcommand{\refFs}[2]{\hyperref[#1]{Figs.~\ref{#1} and \ref{#2}}}
\newcommand{\refT}[1]{\hyperref[#1]{Table~\ref{#1}}}
\newcommand{\refC}[1]{\hyperref[#1]{Section~\ref{#1}}}
\newcommand{\refCon}[1]{\hyperref[#1]{Condition~\ref{#1}}}
\newcommand{\refE}[1]{Eq.~\hyperref[#1]{(\ref{#1})}}
\newcommand{\refEs}[2]{Eqs.~\hyperref[#1]{(\ref{#1})} and \hyperref[#2]{(\ref{#2})}}
\newcommand{\refEss}[2]{Eqs.~\hyperref[#1]{(\ref{#1})}-\hyperref[#2]{(\ref{#2})}}
\newcommand{\Cite}[1]{~\cite{#1}}

\newif\ifMARKED
\MARKEDfalse  % 出版用

\ifMARKED
    \newcommand{\EDITf}[1]{#1}
\else
    \newcommand{\EDITf}[1]{#1}
\fi

\ifMARKED
    \newcommand{\EDITs}[1]{#1}
\else
    \newcommand{\EDITs}[1]{#1}
\fi

\ifMARKED
    \newcommand{\EDIT}[1]{#1}
\else
    \newcommand{\EDIT}[1]{#1}
\fi

\ifMARKED
    \newcommand{\EDITT}[1]{\textcolor[HTML]{FF6600}{#1}}
\else
    \newcommand{\EDITT}[1]{#1}
\fi

\newcommand{\D}{\Delta}
\newcommand{\N}{\nabla}
\graphicspath{{figs/}}
\begin{document}

\begin{frontmatter}

%% Title, authors and addresses

\title{A generalized vertical coordinate transformation based on SPH(2) \\for efficient free surface flow simulations}

\author[add1]{Shujiro Fujioka}
\ead{s-fujioka@doc.kyushu-u.ac.jp}

\author[add2]{Kumpei Tsuji}
\ead{kumpei.tsuji.e1@tohoku.ac.jp}

\author[add3]{Naoto Mitsume}
\ead{mitsume@kz.tsukuba.ac.jp}

\author[add1]{Mitsuteru Asai\texorpdfstring{\corref{mycorrespondingauthor}}{}}
\ead{asai@doc.kyushu-u.ac.jp}

\cortext[mycorrespondingauthor]{Corresponding author}

%% Author affiliation
\affiliation[add1]{organization={Department of Civil Engineering, Kyushu University},
            addressline={744, Motooka, Nishi-ku},
            city={Fukuoka-shi},
            state={Fukuoka},
            postcode={819-0395},
            country={Japan}}

\affiliation[add2]{organization={Department of Civil and Environmental Engineering, Tohoku University},
            addressline={6-6-06, Aza-Aoba, Aramaki, Aoba-ku},
            city={Sendai-shi},
            state={Miyagi},
            postcode={980-8579},
            country={Japan}}

\affiliation[add3]{organization={Institute of Systems and Information Engineering, University of Tsukuba},
            addressline={1-1-1, Tennodai},
            city={Tsukuba-shi},
            state={Ibaraki},
            postcode={305-8577},
            country={Japan}}

%% Abstract
\begin{abstract}
We propose three new particle methods that improve computational efficiency by introducing a generalized Vertical Coordinate Transformation (VCT) for free surface flow problems \EDITf{with complex bottom boundaries}.
The first method is a bottom boundary-fitted particle method (BF-SPH). The BF-SPH is simply an arrangement of the body-fitted-coordinate system in the finite difference method to the particle method. The BF-SPH can accurately impose the bottom boundary conditions, while a simple procedure is performed by transforming the complex bottom into a flat one. The second method is the bottom boundary-fitted ellipsoidal particle method (BFE-SPH), which combines the BF-SPH with the ellipsoidal particle model proposed by Shibata et al. The BFE-SPH can speed up the particle simulation by choosing a reasonable aspect ratio of ellipsoidal particles. The last method is the $\sigma$-SPH method, which automatically selects the aspect ratios of ellipsoidal particles concerning water depth using the $\sigma$-coordinate system. The $\sigma$-coordinate is often employed in numerical simulations of oceanographic fields, such as in the Princeton Ocean Model. However, this is the first attempt to apply the $\sigma$-coordinate to a particle method. \EDITs{Vertical resolution is required from offshore to the coastal region in oceanographic problems such as tsunamis, especially when conducting detailed analysis using a 3-D particle method.} Using $\sigma$-coordinate allows for a stepwise transition to a naturally efficient coordinate system by referencing water depth. In this paper, we have shown that the above three methods can be generalized as Vertical Coordinate Transformations (VCTs), and the VCTs are successfully achieved by employing SPH(2) with the second-order accuracy of the second-order derivatives, including cross derivatives.

\end{abstract}

% %%Graphical abstract
% \begin{graphicalabstract}
% \includegraphics{graphical_abstract.pdf}
% \end{graphicalabstract}

% %%Research highlights
% \begin{highlights}
% \item A vertical coordinate transformation for the SPH is formulated in a unified manner.
% \item The SPH(2) can perform the vertical coordinate transformation properly and accurately.
% \item The $\sigma$-coordinate system has been introduced into SPH for the first time.
% \end{highlights}

%% Keywords

\end{frontmatter}

%% Add \usepackage{lineno} before \begin{document} and uncomment
%% following line to enable line numbers
% \linenumbers

%% main text
%%

%%%%%%%%%%%%%%%%%%%%%%%%%%%%%%%%%%%%%%%%%%%%%%%%%%%%%%%%%%%%%%%%%%%%%
%%  1. Introduction
%%%%%%%%%%%%%%%%%%%%%%%%%%%%%%%%%%%%%%%%%%%%%%%%%%%%%%%%%%%%%%%%%%%%%
\section{Introduction}
\label{sec1:Introduction}
The particle methods, such as the Smoothed Particle Hydrodynamics (SPH) method\Cite{lucy1977numerical, gingold1977smoothed} and the Moving Particle Semi-implicit (MPS) method\Cite{koshizuka1996moving}, discretize the domain into moving Lagrangian particles, eliminating the need for computational grids or meshes. In particle methods, physical quantities and their derivatives are calculated using weight functions based on the distance between the particles, and these values are used to solve the governing equations. These methods have been widely applied in various engineering fields, particularly for simulating moving discontinuities and systems undergoing large deformations, such as free surface flows characterized by breaking, splashing, and fragmentation.\par

However, the particle methods have faced challenges regarding reduced accuracy due to particle movements and disturbances. To address the issue, particle shifting techniques (PSTs) like Particle Shifting (PS)\Cite{xu2009accuracy,lind2012incompressible}, Optimized PS (OPS)\Cite{khayyer2017comparative}, and Density-based PS (DPS)\Cite{morikawa2023corrected}, have been proposed to either equalize particle arrangements or conserve the total volume of particles. In addition, the kernel gradient and Laplacian correction model are widely used for high accuracy, even with disordered particle arrangements. Recently, models such as Fatehi~\&~Manzari's model\Cite{fatehi2011error}, for Laplacian calculations, and Least-Squares MPS (LSMPS)\Cite{tamai2014least} and SPH(2)\Cite{asai2023class}, for the first- and second-derivatives calculations, have been proposed, achieving the 2nd-order accuracy in space. Combining PSTs with high-accuracy approximation models enables simulations to achieve 
an accuracy comparable to mesh-based methods while ensuring stability in analyses of free surface problems. The LSMPS and SPH(2) models can calculate the second-derivatives individually, making them applicable to the coordinate transformations discussed below.\par

Although the particle method is expected to be widely applicable due to improvements in accuracy, the computational cost remains a practical issue, as it requires frequent updates of particle positions and weight functions. Recent studies have focused on parallelization techniques (e.g., Message Passing Interface (MPI)\Cite{ferrari2009new} and Graphics Processing Units (GPUs)\Cite{herault2010sph}) and improved iterative methods for nonlinear differential equations\Cite{chow2018incompressible} to increase computational efficiency. However, high computational costs continue to pose significant challenges because the particle size in the conventional particle methods must be uniform across the domain, and a considerable number of particles are needed to achieve sufficient accuracy and to apply the method to a wide range of simulations. Mesh-based methods often involve adjusting the mesh size according to the importance of the analysis to reduce computational costs, which is impossible with conventional particle methods. Shibata et al. proposed an ellipsoidal particle method\Cite{shibata2016ellipse} that reduces computational costs by introducing a coordinate transformation to change the aspect ratio of particles used in the MPS method, demonstrating the applicability of coordinate transformations to particle methods.\par

In this study, we aim to improve the computational efficiency of the particle method by applying coordinate transformations. This work has been inspired by the body-fitted coordinate system\Cite{thompson1974automatic} in finite-difference methods and the $\sigma$-coordinate system, for example, in the Princeton Ocean Model (POM)\Cite{phillips1957sigma}. The body-fitted coordinate system is obtained by transforming complex physical boundaries into a projected space composed only of flat boundaries. On the other hand, the $\sigma$-coordinate system is a coordinate system used in oceanographic fields. Maintaining a constant number of meshes in the vertical direction regardless of water depth improves efficiency by reducing the number of mesh elements while preserving accuracy in coastal areas, which are critical for simulations. In this study, we have generalized the Vertical Coordinate Transformations (VCTs) before applying them to the SPH method. Then, we propose three new coordinate transformations in the SPH method. The VCTs require the second derivatives, including the cross-derivative terms, and an accurate SPH discretization model, such as SPH(2), should be applied in their implementation. Through the verification, the computational efficiency of VCTs and their accuracy were discussed. At the same time, we investigated the impact of SPH discretization error on the accuracy of VCTs.

%%%%%%%%%%%%%%%%%%%%%%%%%%%%%%%%%%%%%%%%%%%%%%%%%%%%%%%%%%%%%%%%%%%%%
%%  2. Vertical coordinate transformations
%%%%%%%%%%%%%%%%%%%%%%%%%%%%%%%%%%%%%%%%%%%%%%%%%%%%%%%%%%%%%%%%%%%%%
\section{Generalized Vertical Coordinate Transformations (VCTs)}
\label{sec2:VCTs}
\begin{figure}[H]
  \centering
  \begin{subfigure}[b]{\linewidth}
    \centering
    \includegraphics[width=0.8\linewidth]{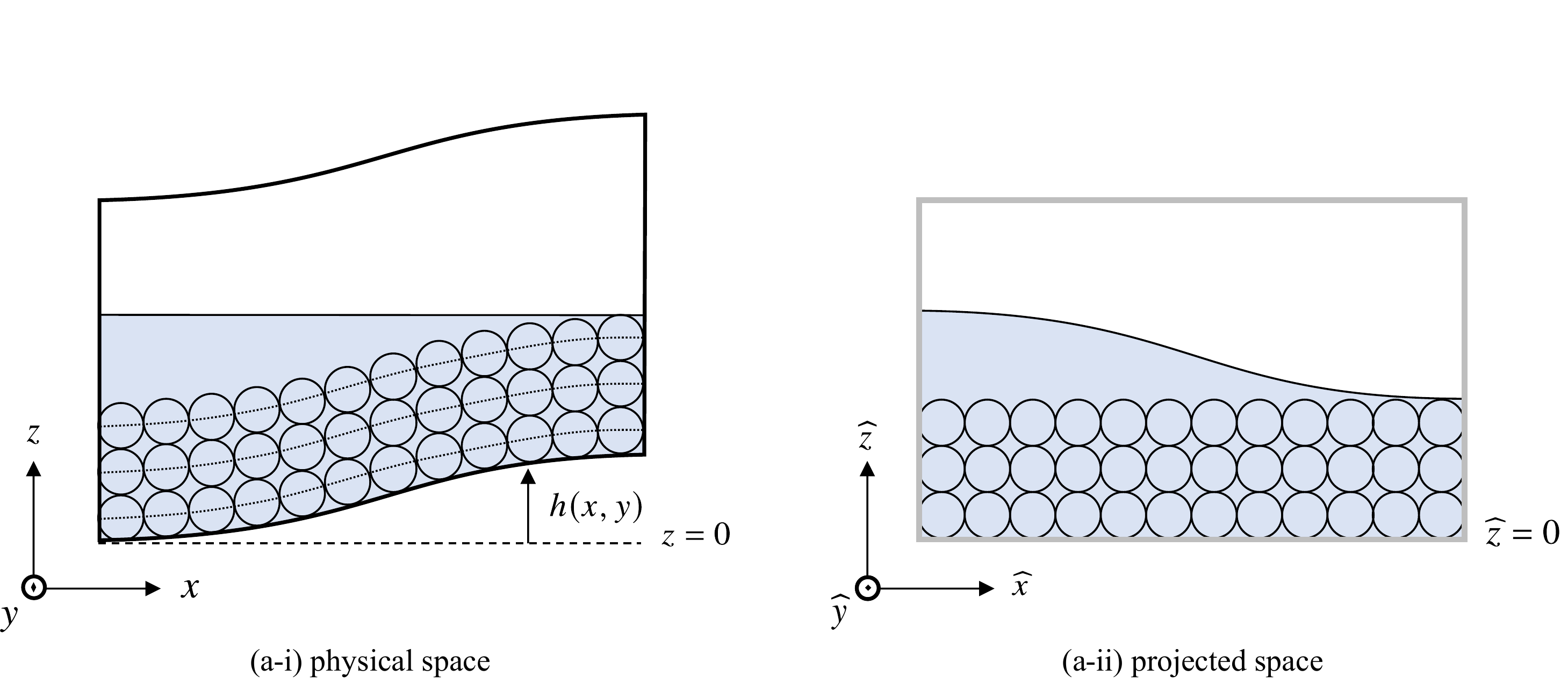}
    \caption{Bottom boundary-fitted particle method (BF-SPH)}\label{fig:s2_b}
  \end{subfigure}
  \begin{subfigure}[b]{\linewidth}
    \centering
    \includegraphics[width=0.8\linewidth]{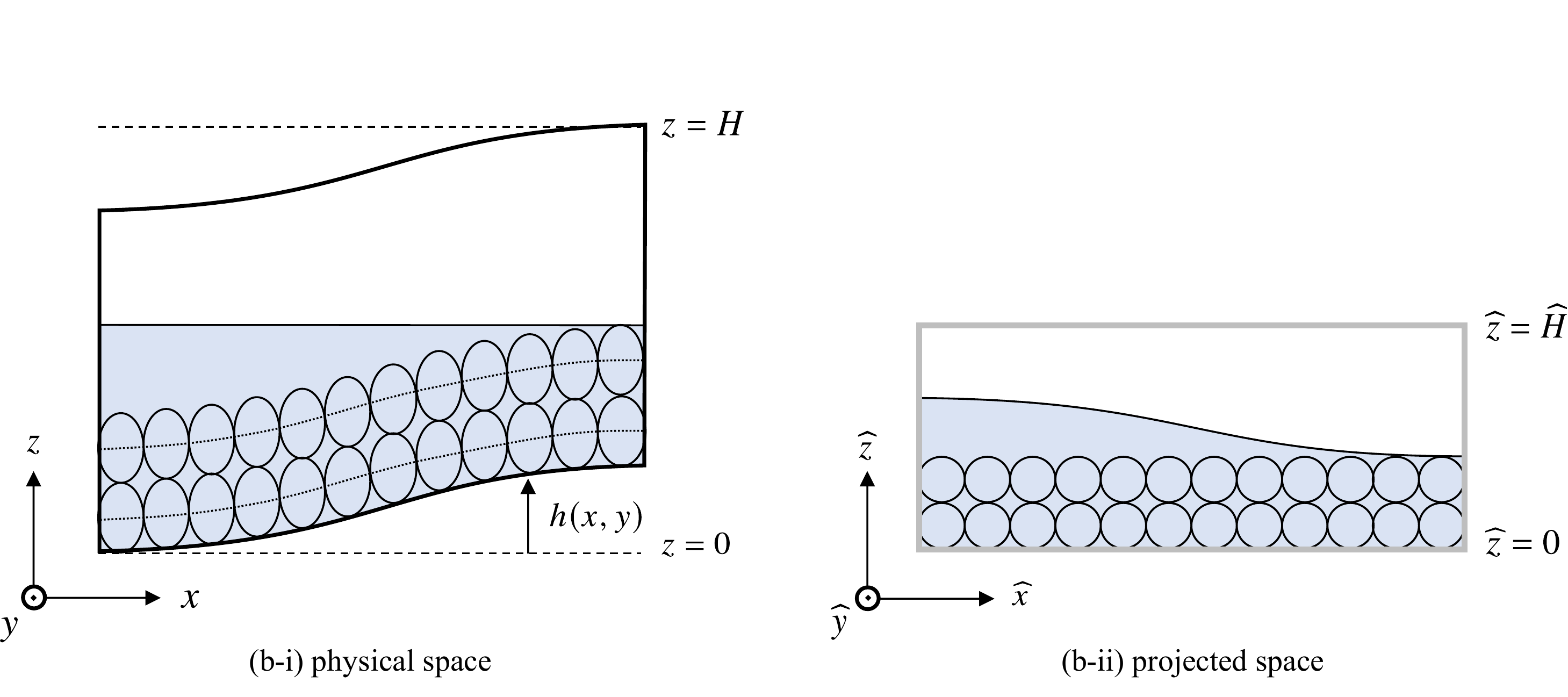}
    \caption{Bottom boundary-fitted ellipsoidal particle method (BFE-SPH)}\label{fig:s2_be}
  \end{subfigure}
  \begin{subfigure}[b]{\linewidth}
    \centering
    \includegraphics[width=0.8\linewidth]{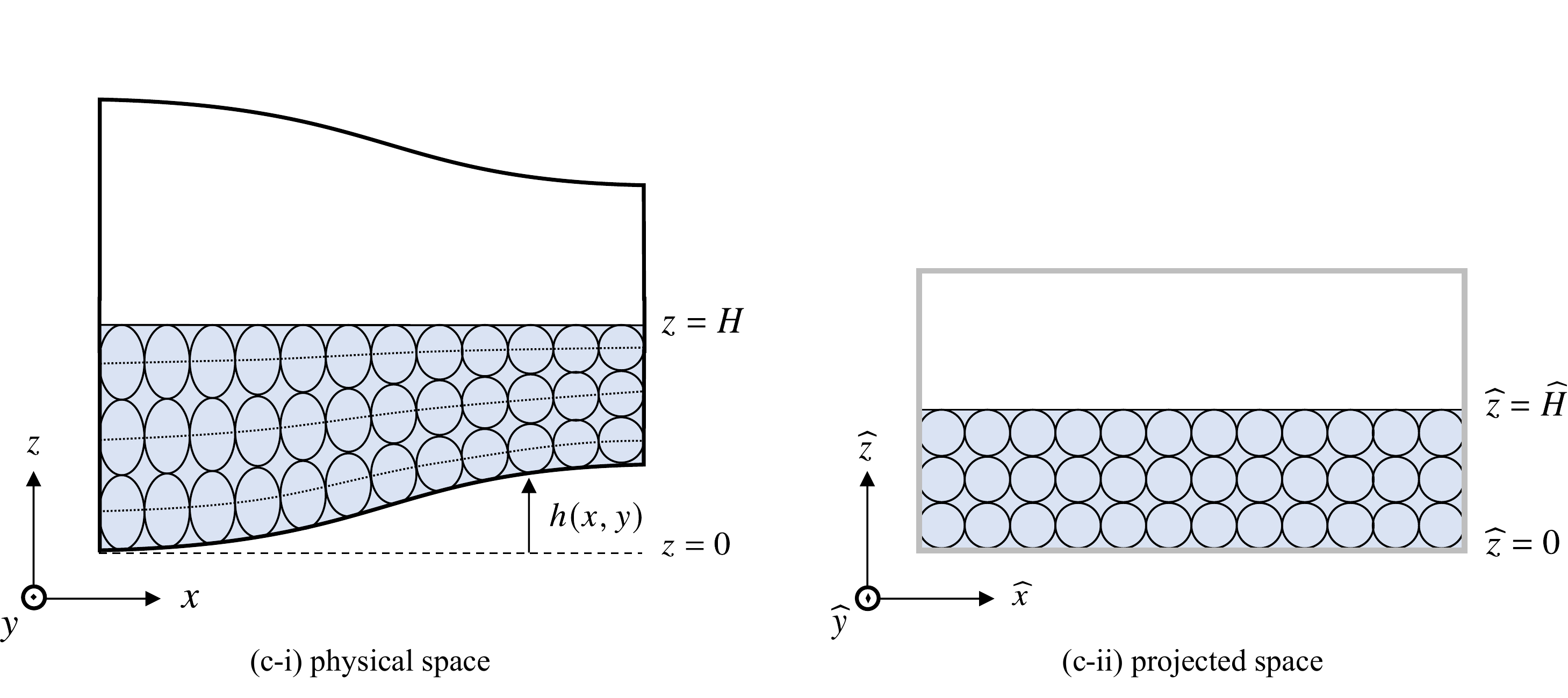}
    \caption{\EDITs{$\sigma$-SPH method using a $\sigma$-coordinate system}}\label{fig:s2_s}
  \end{subfigure}
  \caption{Conceptual diagram of the novel three methods (BF-SPH, BFE-SPH, and $\sigma$-SPH)}\label{fig:s2_VCTs}
\end{figure}

This section proposes three new methods with Vertical Coordinate Transformations (VCTs) as shown in \refF{fig:s2_VCTs}. VCTs involve projecting the physical space onto a normalized coordinate system (projected space), where analyses are performed using the projected space. The details of the VCT procedure will be discussed later. The three new methods are listed as follows;
\begin{enumerate}
    \item Bottom boundary-fitted particle method (BF-SPH) as shown in \hyperref[fig:s2_b]{Fig.~\ref*{fig:s2_VCTs}(a)}
    \item Bottom boundary-fitted ellipsoidal (elliptical in 2-D)\footnote{For consistency in terminology, this paper uses the term ``ellipsoidal'' even for 2-D cases.} particle method (BFE-SPH) as shown in \hyperref[fig:s2_be]{Fig.~\ref*{fig:s2_VCTs}(b)}
    \item $\sigma$-SPH method using a $\sigma$-coordinate system as shown in \hyperref[fig:s2_s]{Fig.~\ref*{fig:s2_VCTs}(c)}
\end{enumerate}
These methods can be formulated in the same manner as VCTs.

%%-------------------------------------------------------------------
%%  2. 1. Generalization of VCTs
%%-------------------------------------------------------------------
\subsection{Coordinate transformation in the vertical direction}
\label{sec2-1:Generalization_of_VCTs}
In this study, the coordinate transformations are considered in the vertical direction, so the coordinate transformations are defined as follows:
\begin{equation}
  \label{eq:VCT}
   \bm{r}=
   \begin{bmatrix}
    
    x \vspace{5pt}\\
    y \vspace{5pt}\\
    z 
   \end{bmatrix}
   \longleftrightarrow
   \wh{\bm{r}}=
   \begin{bmatrix}
    \,\wh x\, \vspace{5pt}\\
    \,\wh y\, \vspace{5pt}\\
    \,\wh z\, 
    \end{bmatrix}
    =
    \begin{bmatrix}
      x \vspace{5pt}\\
      y \vspace{5pt}\\
    \alpha(x,y)\bab{z+\beta(x,y)}
    \end{bmatrix},
\end{equation}
where $\bm{r}$ and $\wh{\bm{r}}$ are the position vectors before and after the coordinate transformation. $\alpha(x,y)$ and $\beta(x,y)$ are arbitrary variables for defining the coordinate transformations. $\alpha$ defines the vertical scale factor ($\alpha>0$) and $\beta$ defines the vertical height correction. The accent $\wh \blacksquare$ indicates that the values are in space after the coordinate transformations. In this paper, the space before the VCTs is referred to as the physical space, and the space after the VCTs is referred to as the projected space.\par
The Jacobian matrix $\bm J$ of \refE{eq:VCT} is expressed as
\begin{equation}
  \label{eq:jacobi}
  \bm J=
  \begin{bmatrix}
    \pdv{\,\wh x}{x}&\pdv{\,\wh y}{x}&\pdv{\,\wh z}{x}\, \vspace{7pt}\\
    \pdv{\,\wh x}{y}&\pdv{\,\wh y}{y}&\pdv{\,\wh z}{y}\, \vspace{7pt}\\
    \pdv{\,\wh x}{z}&\pdv{\,\wh y}{z}&\pdv{\,\wh z}{z}\,
  \end{bmatrix}
  =
  \begin{bmatrix}
    1&0&\alpha_x(z+\beta)+\alpha\beta_x \vspace{7pt}\\
    0&1&\alpha_y(z+\beta)+\alpha\beta_y \vspace{7pt}\\
    0&0&\alpha
  \end{bmatrix},
\end{equation}
\begin{equation}
  \alpha_{r^{\mathrm{I}}}:=\pdv{\alpha}{r^{\mathrm{I}}},\,\beta_{r^{\mathrm{I}}}:=\pdv{\beta}{r^{\mathrm{I}}}\,\,\,,
\end{equation}
where $r^{\mathrm{I}}$ denotes the coordinates along $r^{\mathrm{I}}$-axis ($r^{\mathrm{1}}=x$, $r^{\mathrm{2}}=y$, and $r^{\mathrm{3}}=z$). 
As shown in \refE{eq:jacobi}, the Jacobian determinant $J$ is expressed as
\begin{align}
    \label{eq:det_jacobi}
    J:=\det \bm J=\alpha
\end{align}
Therefore, volume conservation is satisfied when $\alpha$ is constant. However, when $\alpha$ varies, volume changes must be considered. The process for this is described in \refC{sec3-6:volume_conservation}.\par
The gradient and Laplacian in physical space are described using the Jacobian matrix $\bm J$ as follows:
\begin{equation}
  \label{eq:grad_VCT}
  \N\phi=\bm J\wh\N\phi;\quad
  \wh\N=
  \begin{bmatrix}
    \pdv{}{\,\wh x}&\pdv{}{\,\wh y}&\pdv{}{\,\wh z}
  \end{bmatrix}^T \,\,\,,
\end{equation}
\begin{equation}
  \label{eq:lap_VCT}
  \lap\phi=\N\cdot\N\phi=\bm J\wh\N\cdot\bm J\wh\N\phi=:\bm c_{\mathrm{trans}}\cdot\wh{\bm D}\phi \,\,\,,
\end{equation}
\begin{equation}
  \EDITT{
  \bm c_{\mathrm{trans}}:=
  \begin{bmatrix}
    \pab{\alpha_{xx} + \alpha_{yy} }(z+\beta) +
    2\pab{\alpha_x\beta_x + \alpha_y\beta_y} + \alpha\pab{\beta_{xx}      + \beta_{yy}} \vspace{5pt}\\
    1 \vspace{5pt}\\
    1 \vspace{5pt}\\
    \alpha^2+\bab{\alpha_x(z+\beta)+\alpha\beta_x}^2+\bab{\alpha_y(z+\beta)+\alpha\beta_y}^2 \vspace{5pt}\\
    2\bab{\alpha_x(z+\beta)+\alpha\beta_x} \vspace{5pt}\\
    2\bab{\alpha_y(z+\beta)+\alpha\beta_y}
  \end{bmatrix} \,\,\,},  
\end{equation}
\begin{equation}
  \wh{\bm D}:=
  \begin{bmatrix}
    \pdv{}{\,\wh z\,}&\pdv[2]{}{\,\wh x\,}&\pdv[2]{}{\,\wh y\,}&\pdv[2]{}{\,\wh z\,}&\mdv{}{\,\wh x\,}{\,\wh z\,}&\mdv{}{\,\wh y\,}{\,\wh z\,}
  \end{bmatrix}^T\,\,\,,
\end{equation}
\begin{equation}  
  \alpha_{r^{\mathrm{I}}r^{\mathrm{J}}}:=\mdv{\alpha}{r^{\mathrm{I}}}{r^{\mathrm{J}}},\,
  \beta_{r^{\mathrm{I}}r^{\mathrm{J}}}:=\mdv{\beta}{r^{\mathrm{I}}}{r^{\mathrm{J}}}\,\,\,.
\end{equation}\par
In the VCTs, the values of $\alpha$ and $\beta$ should be selected for our purposes, and the governing equations are solved using \refEs{eq:grad_VCT}{eq:lap_VCT}. The next section presents several examples and proposes three novel methods.
%%-------------------------------------------------------------------
%%  2. 2. Applications of VCTs
%%-------------------------------------------------------------------
\subsection{Applications of VCTs}
\label{sec2-2:Applications_of_VCTs}
\subsubsection{Ellipsoidal particle method \EDITs{based on SPH method} (E-SPH)}
The ellipsoidal particle method \EDITs{based on MPS method}\Cite{shibata2016ellipse} reduces the total number of particles by changing their aspect ratio, thus improving efficiency. In a previous study, the aspect ratio was allowed to vary in each coordinate direction. In this paper, as only vertical changes are considered, $\alpha$ and $\beta$ are expressed as
\begin{equation}
  \alpha(x,y)=\cfrac{\wh H}{H};\quad
  \beta(x,y)=0\,\,\,,
\end{equation}
where $\wh H$ and $H$ are arbitrary heights in projected and physical space, respectively. In physical space, particles become horizontally elongated ellipses in case $\alpha>1$ and vertically elongated ellipses in case $0<\alpha<1$.

\subsubsection{Bottom boundary-fitted particle method (BF-SPH)}
The first proposed method is the bottom boundary-fitted particle method inspired by the body-fitted coordinate system\Cite{thompson1974automatic}. This method utilizes a VCT that depends on the elevation of the bottom to transform a complex bottom surface into a flat one, as shown in \hyperref[fig:s2_b]{Fig.~\ref*{fig:s2_VCTs}(a)}. This method makes it easy to represent boundaries that contain curved surfaces, such as seabeds and riverbeds. Therefore, this method is expected to provide a highly accurate analysis of problems \EDITf{with complex bottom boundaries}. In this method, $\alpha$ and $\beta$ are expressed as
\begin{equation}
  \alpha(x,y)=1;\quad
  \beta(x,y)=-h(x,y),
\end{equation}
where $h$  is the bottom elevation in physical space.

\subsubsection{Bottom boundary-fitted ellipsoidal particle method  (BFE-SPH)}
The second method we propose is the bottom boundary-fitted ellipsoidal particle method. This method combines the ellipsoidal particle method with the bottom boundary-fitted particle method, as shown in \hyperref[fig:s2_be]{Fig.~\ref*{fig:s2_VCTs}(b)}. Thus, this method is expected to improve accuracy and computational efficiency. In this method, $\alpha$ and $\beta$ are expressed as
\begin{equation}
  \alpha(x,y)=\cfrac{\wh H}{H};\quad
  \beta(x,y)=-h(x,y).
\end{equation}

\subsubsection{\texorpdfstring{$\sigma$}{sigma}--SPH method using a \texorpdfstring{$\sigma$}{sigma}--coordinate system}
The third method we propose is the $\sigma$-SPH method using a $\sigma$-coordinate system\Cite{phillips1957sigma}. This coordinate system is used in computational models for oceanography, meteorology, and other fields of fluid dynamics. This method allows the resolution to vary according to the importance of the analysis by applying a VCT depending on the depth of the water, as shown in \hyperref[fig:s2_s]{Fig.~\ref*{fig:s2_VCTs}(c)}. As a result, the total number of particles can be reduced, which is expected to improve the computational cost further. In this method, $\alpha$ and $\beta$ are expressed as
\begin{equation}
  \alpha(x,y)=\cfrac{\wh H}{H-h(x,y)};\quad
  \beta(x,y)=-h(x,y).
\end{equation}
If $\alpha$ takes an extremely small or large value, the distortion of particle spacing in each direction may lead to computational instability. Therefore, to prevent such instability, it is recommended to use a constant value for $\alpha$ when it exceeds the applicable range. The representative volume of each particle should be updated as we will discuss the details in \refC{sec3-6:volume_conservation} since $\alpha$ equals to the Jacobian determinant $J$ as shown in \refE{eq:det_jacobi} and is changed in space.
% As shown in \refE{eq:det_jacobi}, we need to consider the change in volume in this method, since $\alpha$ is variable.

%%%%%%%%%%%%%%%%%%%%%%%%%%%%%%%%%%%%%%%%%%%%%%%%%%%%%%%%%%%%%%%%%%%%%
%%  3. Incompressible SPH method
%%%%%%%%%%%%%%%%%%%%%%%%%%%%%%%%%%%%%%%%%%%%%%%%%%%%%%%%%%%%%%%%%%%%%
\section{Incompressible SPH method}
\label{sec3:ISPH}
%%-------------------------------------------------------------------
%%  3. 1. SPH approximations
%%-------------------------------------------------------------------
\subsection{SPH approximations}
\label{sec3-1:SPH_approximations}
\subsubsection{Conventional models}
\EDITs{In the SPH method\Cite{lucy1977numerical,gingold1977smoothed}, the physical quantity $\phi_i:=\phi(\bm{r}_i,t)$ at the location of the target particle $i$ in the domain $\Omega$ can be expressed as a volume integration using a smoothing kernel function $w$, an infinitesimal volume element $d\Omega\in\Omega$, and a position vector $\bm \xi\in d\Omega$ as
\begin{equation}
  \label{eq:kernel_int}
  \phi_i \approx \int_{\Omega}\phi(\bm \xi,t)\, w(\vab{\bm \xi-\bm{r}_i},l)\,d\Omega ,
\end{equation}
% where $\bm{r}_{ij}(:=\bm{r}_j-\bm{r}_i)$ is a relative position vector of particle $i$ and $j$, $t$ is time, and $l$ is the smoothing length. Here, we employ the following cubic spline function as the kernel function: 
where $t$ is time, $\bm{r}_i$ is a position vector of particle $i$, and $l$ is the smoothing length.} Here, we employ the following cubic spline function as the kernel function: 
\begin{equation}
  w(r,l)=\alpha_d^c
  \begin{cases}
    1-\cfrac{3}{2}\,\pab{\cfrac{r}{l}\,}^2+\cfrac{3}{4}\,\pab{\cfrac{r}{l}\,}^3 &\pab{0\le \cfrac{r}{l} \le 1}\\
    \cfrac{1}{4}\,\pab{2-\cfrac{r}{l}\,}^3         &\pab{1\le \cfrac{r}{l} \le 2}\\
    0                  &\pab{2\le \cfrac{r}{l}\, } 
  \end{cases}.
\end{equation}
Here, $\alpha_d^c$ is a constant chosen to satisfy the unity condition, and its value is $10/(7\pi l^2)$ in two-dimensional simulations or $1/(\pi l^3)$ in three-dimensional simulations. Note that \EDITs{the smoothing length $l$} and the effective radius $r_e$ are set to $r_e = 2.0 h$, $l = 1.2 d_0$ with the initial particle distance $d_0$. Let $\mathbb{S}_i$ be defined as the set of neighbor particles $j$ of the target particle $i$ as follows:
\begin{equation}
  \label{eq:SPH_neighbor_set}
  \mathbb{S}_i := \Bab{ \; j = 1, 2, \cdots , N_{\mathrm{SPH}} \; | \;  r_e >  |\bm{r}_{ij}| \; \wedge \; \bm{r}_j \in \Omega \; },
\end{equation}
where $N_{\mathrm{SPH}}$ is the number of SPH particles, and \EDITs{$\bm{r}_{ij}(:=\bm{r}_j-\bm{r}_i)$ is a relative position vector of particle $i$ and $j$}. For SPH, the volume integration in \refE{eq:kernel_int}, $\phi$ and its derivatives (e.g., $\nabla \phi$, $\nabla\cdot \bm \phi$, and $\lap \phi$) can be approximated as

\begin{align}
  \aab{\phi}_i^{(0)}&:=\sphsum{V_j\phi_jw_{ij}}\, ,\\
  \aab{\N\phi}_i^{(0)}&:=\sphsum{V_j\phi_{ij}\N w_{ij}}\, ,\\
  \aab{\N\cdot\bm\phi}_i^{(0)}&:=\sphsum{V_j\bm\phi_{ij}\cdot\N w_{ij}}\, ,\\
  \aab{\lap\phi}_i^{(0)}&:=2\sphsum{V_j\cfrac{\bm{r}_{ij}\cdot\N w_{ij}}{|\bm{r}_{ij}|^2} \phi_{ij}},
\end{align}
where $\phi_{ij}:=\phi_j-\phi_i$, $\bm\phi_{ij}:=\bm\phi_j-\bm\phi_i$. $V_j$ is the representative volume of each neighboring particle $j$, and $\aab{\blacksquare}^{(\mathsf{N})}$ means a model with spatial $\mathsf{N}$th-order accuracy. For simulations using VCTs, the values for each of the second-derivatives, including the crossed derivative, are required, as shown in \refE{eq:lap_VCT}. A model that can compute the respective values of the second-order derivatives is the model by Espan\~ol~\&~Revenga\Cite{espanol2003smoothed} presented as follows:
\begin{equation}
  \label{eq:sec_0th}
  \aab{\mdv{\phi}{r^{\mathrm{I}}}{r^{\mathrm{J}}}}_i^{(0)}:=
  \begin{cases}\vspace{7pt}
    \displaystyle\sphsum{V_j}\pab{\cfrac{4r^{\mathrm{I}}_{ij}r^{\mathrm{J}}_{ij}}{|\bm{r}_{ij}|^2}-\delta^{\mathrm{IJ}}}\cfrac{\bm{r}_{ij}\cdot\N w_{ij}}{|\bm{r}_{ij}|^2}\phi_{ij} & \mr{(2D)}\\
    \displaystyle\sphsum{V_j}\pab{\cfrac{5r^{\mathrm{I}}_{ij}r^{\mathrm{J}}_{ij}}{|\bm{r}_{ij}|^2}-\delta^{\mathrm{IJ}}}\cfrac{\bm{r}_{ij}\cdot\N w_{ij}}{|\bm{r}_{ij}|^2}\phi_{ij} & \mr{(3D)}
  \end{cases},
\end{equation}
where $\delta^{\mathrm{IJ}}$ is Kronecker delta. The effects of these 0th-order models on inaccuracy and numerical instability were demonstrated by Asai et al.\Cite{asai2023class}. Therefore, in this study, corrected models, shown in the next section, are used to improve accuracy.
\subsubsection{Corrected models}
Applying corrections to the SPH approximation models is crucial, especially in regions where neighboring particles are disturbed, to satisfy unity conditions and maintain mathematical consistency.\par
The 1st-order accurate corrected gradient and divergence models were proposed by Randles~\&~Libersky\Cite{randles1996smoothed} and Bonet~\&~Lok\Cite{bonet1999variational} and have been widely used in the SPH method, as follows:
\begin{align}
  \label{eq:gra_1st}
  \aab{\N\phi}_i^{(1)}&:=\sphsum{V_j\phi_{ij}\wt\N w_{ij}},\\
  \label{eq:div_1st}
  \aab{\N\cdot\bm\phi}_i^{(1)}&:=\sphsum{V_j\bm\phi_{ij}\cdot\wt\N w_{ij}},\\
  \wt\N w_{ij}:=\bm L_i\N w_{ij};\quad
  &\bm L_i:=\bab{\sphsum{V_j\pab{\N w_{ij}\otimes\bm{r}_{ij}}}}^{-1}.
\end{align}\par
Furthermore, Asai et al.\Cite{asai2023class} proposed an even more accurate model called SPH(2). This model can evaluate the second-derivatives, including crossed derivatives, with the 2nd-order accuracy in space. Although gradient and divergence with the 2nd-order accuracy can be obtained using SPH(2), the 1st-order accuracy model is used in this study because there were no significant differences in accuracy, according to\Cite{asai2023class}. The formulation of SPH(2) in two dimensions is shown as
\begin{equation}
  \label{eq:sec_SPH2}
  \aab{\bm D\phi}_i^{(2)}
  :=2\bm M_i^{-1}
  \sphsum{V_jF_{ij}\bm q_{ij}\pab{\phi_{ij}-\bm{r}_{ij}\cdot\aab{\N\phi}_i^{(1)}}},
\end{equation}
\begin{equation}
  \bm D :=
  \begin{bmatrix}
    \pdv[2]{}{x} &
    \pdv[2]{}{y} &
    2\mdv{}{x}{y}
  \end{bmatrix}^T,
\end{equation}
\begin{equation}
\bm M_i:=\sphsum{V_jF_{ij}\bm q_{ij}\bm p_{ij}^T};\quad
  F_{ij}:=\cfrac{\bm{r}_{ij}\cdot\wt\N w_{ij}}{|\bm{r}_{ij}|^4},
\end{equation}
\begin{equation}
  \bm p_{ij}:=
  \begin{bmatrix}
    x_{ij}^2-    \displaystyle \sphsum[k]{V_kx_{ik}^2\wt\N w_{ik}}\\
    y_{ij}^2-    \displaystyle \sphsum[k]{V_ky_{ik}^2\wt\N w_{ik}}\\
    x_{ij}y_{ij}-\displaystyle \sphsum[k]{V_kx_{ik}y_{ik}\wt\N w_{ik}}
  \end{bmatrix};\quad
  \bm q_{ij}:=
  \begin{bmatrix}
    x_{ij}^2 \vspace{7pt}\\
    y_{ij}^2 \vspace{7pt}\\
    x_{ij}y_{ij}
  \end{bmatrix}.
\end{equation}
%%-------------------------------------------------------------------
%%  3. 2. Governing equations
%%-------------------------------------------------------------------
\subsection{Governing equations for incompressible fluid}
\label{sec3-2:Governing_equations}
This study focuses on incompressible flows of Newtonian fluids with constant density and variable viscosity. Therefore, the governing equations are the following continuity and Navier–Stokes equations:
\begin{equation}
  \N\cdot\bm{v}=0,
\end{equation}
\begin{equation}
  \label{eq:NS}
  \cfrac{D\bm{v}}{Dt}=-\cfrac{1}{\rho}\N p+\nu\lap{\bm{v}}+\bab{\N\otimes\bm{v}+\pab{\N\otimes\bm{v}}^T}\N \nu+\bm g;\quad\nu=\nu_{\mathrm{K}}+\nu_{\mathrm{E}},
\end{equation}
where $\bm{v}$ is fluid velocity, $\rho$ is fluid reference density, $p$ is fluid pressure, $\nu_{\mathrm{K}}$ is the kinematic viscosity, $\nu_{\mathrm{E}}$ is the eddy viscosity, and $\bm g$ is gravitational acceleration. As demonstrated later, the Smagorinsky model is used for eddy viscosity $\nu_{\mathrm{E}}$.

%%-------------------------------------------------------------------
%%  3. 3. ISPH projection scheme
%%-------------------------------------------------------------------
\subsection{ISPH projection scheme}
\label{sec3-3:projection_scheme}
Following the projection method\Cite{chorin1968numerical}, \refE{eq:NS} is split into two steps as
\begin{align}
  \text{Predictor step:} & \quad \bm{v}^* = \bm{v}^N + \D t\,\Bab{ \nu \N^2 \bm{v}^N + \bab{\N\otimes\bm{v}^N+\pab{\N\otimes\bm{v}^N}^T}\N \nu^N +\bm{g} }, \\
  \text{Corrector step:} & \quad \bm{v}^{N+1} = \bm{v}^* - \D t\,\pab{ \frac{1}{\rho} \N p^{N+1} },
\end{align}
where subscripts $N$, $*$, and $N+1$ indicate current, predictor, and next-time steps, respectively. The following model\EDITs{\Cite{monaghan1992smoothed}}, called the ``summation'' model, with high numerical stability, is widely used for calculating pressure gradients.
\begin{equation}
  \label{eq:gra_sum}
  \aab{\N p}_i^{+}:=\rho_i\sphsum{m_j\pab{\cfrac{p_j}{\rho_j^2}+\cfrac{p_i}{\rho_i^2}}\N w_{ij}}.
\end{equation}\par
The pressure in the $N+1$ step is calculated by solving the pressure Poisson equation as follows:
\begin{equation}
  \label{eq:ppe}
  \lap p^{N+1}=\cfrac{\rho}{\D t}\N\cdot\bm{v}^*.
\end{equation}
In this study, we use the pressure Poisson equation in the stabilized ISPH method proposed by Asai et al.\Cite{asai2012stabilized} as 
\begin{equation}
  \aab{\lap p^{N+1}}_i=\cfrac{\rho}{\D t}\aab{\N\cdot\bm{v}^*}_i+\gamma\cfrac{\rho-\aab{\rho^N}_i}{\D t^2},
\end{equation}
where the second stabilization term becomes \EDITs{a non-zero positive (or negative) value when the numerical density is smaller (or larger)} than the reference density. The positive value $\gamma$~$(\ll1)$ is a coefficient to maintain the total fluid volume as discussed in\Cite{asai2012stabilized}. The coefficient $\gamma$ is fixed at $\Delta t \times 10^2$~[-] in this paper. Finally, we update the position based on the velocity $\bm{v}^{N+1}$ as follows:
\begin{equation}
  \bm{r}^{N+1}=\bm{r}^N+\bm{v}^{N+1}\D t.
\end{equation}
Using Vertical Coordinate Transformations (VCTs), solve the following equation:
\begin{align}
  \label{eq:VCT_predictor}
  \text{Predictor step:} &\quad \bm{v}^*=\bm{v}^N+\D t\, \Bab{\nu\bm c_{\mathrm{trans}}\cdot\wh{\bm D}{\bm{v}^N}+ \bab{\bm J\, \wh\N\otimes\bm{v}^N+\pab{ \, \bm J \, \wh\N\otimes\bm{v}^N}^T} \, \bm J\, \wh\N \nu^N +\bm g},\\
  \text{Corrector step:} &\quad \bm{v}^{N+1}= \bm{v} ^*-\D t\pab{\cfrac{1}{\rho} \, \bm J\wh\N p^{N+1}},
\end{align}
\begin{equation}
  \label{eq:VCT_ppe}
  \bm c_{\mathrm{trans}}\cdot\aab{\wh{\bm D} p^{N+1}}_i=\cfrac{\rho}{\D t}\, \bm J \, \aab{\wh\N\cdot\bm{v}^*}_i+\gamma\cfrac{\rho-\aab{\wh\rho^N}_i}{\D t^2},
\end{equation}
\begin{equation}
  \label{eq:VCT_update}
  \bm{r}^{N+1}=\bm{r}^N+\bm{v}^{N+1}\D t,
\end{equation}
where $\aab{\wh\rho^N}_i$ denotes the numerical density computed in the projected space. The particle position in the projected space is updated based on $\bm{r}^{N+1}_{ij}$ and \refE{eq:VCT}. \refEss{eq:VCT_predictor}{eq:VCT_update} consist only of physical quantities in physical space and derivatives that can be calculated in projected space. Therefore, the computational cost increase from the application of VCTs is minimal.
%%-------------------------------------------------------------------
%%  3. 4. Smagorinsky model
%%-------------------------------------------------------------------
\subsection{Smagorinsky model}
\label{sec3-4:Smagorinsky}
For eddy viscosity, we use the Smagorinsky model\Cite{smagorinsky1963general} as follows:
\begin{align}
  \nu_{\mathrm{E}}&=\pab{C_sf_s\D}^2\sqrt{2\bm S:\bm S},\\
  \bm S&:=\cfrac{1}{2}\bab{\N\otimes\bm{v}+\pab{\N\otimes\bm{v}}^T},
\end{align}
where $C_s$ is the Smagorinsky constant, set to the commonly used value of $C_s = 0.2$ based on a previous study\Cite{asai2012stabilized}; $\D$ represents the filter width (taken as $\D = r_e$), and $f_s$ is the damping function near the wall. In the SPH method, one approach to correct the eddy viscosity near the wall is to switch to the Reynolds-Averaged Navier-Stokes (RANS) model\Cite{launder2002closure} in the vicinity of the wall\Cite{nakayama2022wall}. \EDITs{The van Driest wall damping function\Cite{vandriest1956onturbulent} is commonly used in mesh-based methods to model turbulence near solid boundaries, where the distance from the wall is well defined throughout the computational domain. In SPH methods, however, particle interactions are determined based on local neighbor searches using background cells, and as a result, the distance to the wall is defined only for particles located near the wall surface. This characteristic makes it challenging to apply the van Driest damping function in a consistent manner within SPH frameworks.} This study uses the following damping functions, analogous to the van Driest wall damping function.
\begin{align}
\EDITs{
    f_s=1-\exp\pab{-\cfrac{10\,d_0}{d_s}\cfrac{d_{\mathrm{wall}}}{r_e}},
}
\end{align}
where $d_0$ is the initial particle spacing, $d_s$ is the reference particle spacing, and $d_{\mathrm{wall}}$ is distance from wall surface. In this study, $d_s$ is set to 1.0~cm. This damping function is applied to particles within a distance of effective radius $r_e$ from the wall surface. The eddy viscosity shows a significant gradient near the wall and gradually diminishes with increasing distance. Therefore, the damping function is adjusted to provide greater damping with smaller initial particle spacing $d_0$.

%%-------------------------------------------------------------------
%%  3. 5. Boundary conditions
%%-------------------------------------------------------------------
\subsection{Boundary conditions}
\label{sec3-5:Boundary_conditions}
For the judgment of the free surface, which provides the Dirichlet boundary condition of zero pressure, we follow the same methods as Matsunaga and Koshizuka\Cite{matsunaga2022stabilized} and Marrone et al.\Cite{marrone2010fast}.\par
\begin{figure}[H]
  \centering
  \begin{subfigure}[b]{0.49\linewidth}
    \centering
    \includegraphics[width=0.8\linewidth]{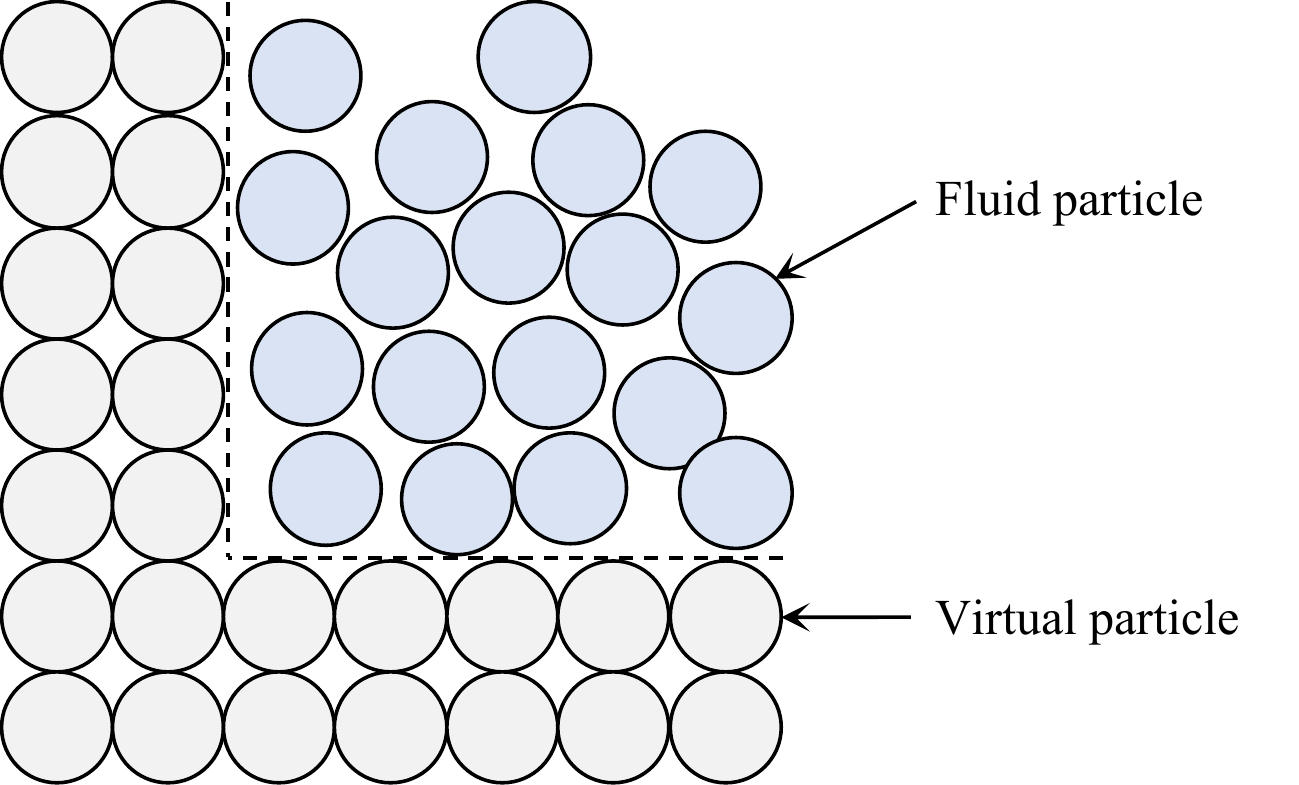}
    \caption{Virtual particle}
  \end{subfigure}
  \begin{subfigure}[b]{0.49\linewidth}
    \centering
    \includegraphics[width=0.731\linewidth]{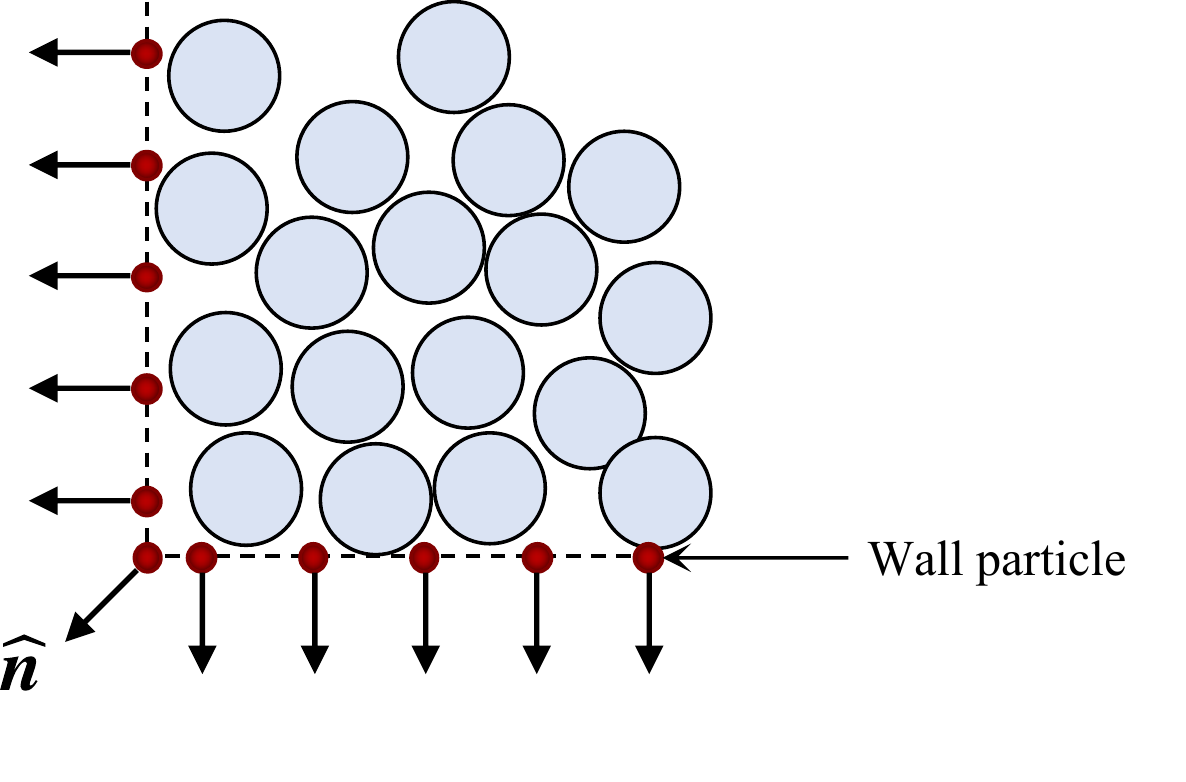}
    \caption{Wall particle}
  \end{subfigure}
  \caption{Conceptual diagram of two types of particles: virtual and wall particles}
  \label{fig:s3_bc}
\end{figure}

Regarding walls, two types of particles are prepared: virtual and wall particles, as shown in \refF{fig:s3_bc}. Virtual particles outside the wall are used to calculate the particle concentration of the PSTs and the numerical density of the stabilization terms in the Stabilized ISPH method. Wall particles on the wall surface are used for calculations relating to velocity and pressure.\par
For the velocity at the wall surface, the Dirichlet conditions are given as follows:
\begin{align}
  \text{No-slip:} &\quad \bm{v}=\bm{v}_{\mathrm{wall}},\\
  \text{Free-slip:} &\quad \bm{v}=\pab{\bm I-2\bm n\otimes\bm n}\,\bm{v}_{\mathrm{wall}},
\end{align}
where $\bm I$ is the identity matrix, and $\bm n$ is the outer normal vector to the wall surface in physical space. For the free-slip condition, the value of 
$\bm{v}_{\mathrm{wall}}$ is obtained using the SPH approximation with particles placed on the boundary. When using VCTs, $\bm n$ can be calculated using the normal vector $\wh{\bm n}$ in projected space as
\begin{equation}
  \bm n=\bm R_x\bm R_y\wh{\bm n},
\end{equation}
\begin{align}
  \bm R_x&:=
  \begin{bmatrix}
  1&                          0&                           0\\
  0&\cos\pab{\arctan\pdv{h}{y}}&-\sin\pab{\arctan\pdv{h}{y}}\\
  0&\sin\pab{\arctan\pdv{h}{y}}& \cos\pab{\arctan\pdv{h}{y}}\\
  \end{bmatrix},\\%;\quad
  \bm R_y&:=
  \begin{bmatrix}
  \cos\pab{\arctan\pdv{h}{x}}&0& \sin\pab{\arctan\pdv{h}{x}}\\
  0&                          1&                           0\\
  -\sin\pab{\arctan\pdv{h}{x}}&0& \cos\pab{\arctan\pdv{h}{x}}\\
  \end{bmatrix}.
\end{align} 
\par
\EDITs{For the pressure at a wall particle $j$, assuming that the viscous term $\nu \nabla^2 \bm{v}$ and the term involving velocity gradient tensor $[\nabla \otimes \bm{v} + (\nabla \otimes \bm{v})^T]\nabla \nu$ are the same as those of fluid particle $i$, the pressure Neumann boundary condition is given as follows:}
\begin{align}
  \pdv{p^{N+1}}{n}&\approx
  \cfrac{p_j^{N+1}-p^{N+1}_i}{\bm{r}_{ij}}
  =\rho\bab{\nu\aab{\lap{\bm{v}^N}}_i+\aab{\N\otimes\bm{v}^N+\pab{\N\otimes\bm{v}^N}^T}_i\aab{\N \nu}^N_i+\bm g},\\
  \label{eq:PNC}
  p_j^{N+1}&=p_i^{N+1}+\rho\bm{r}_{ij}\cdot\bab{\nu\aab{\lap{\bm{v}^N}}_i+\aab{\N\otimes\bm{v}^N+\pab{\N\otimes\bm{v}^N}^T}_i\aab{\N \nu}^N_i+\bm g}.
\end{align}
\EDITT{The pressure Neumann condition serves to prevent fluid particles from penetrating the wall. Therefore, when calculating the properties of a target fluid particle $i$, the pressures of the neighboring wall particles $j$ are extrapolated using \refE{eq:PNC}.}

%%-------------------------------------------------------------------
%%  3. 6. The representative volume change in the \texorpdfstring{$\sigma$}{sigma}--SPH
%%-------------------------------------------------------------------
\subsection{The representative volume change in the \texorpdfstring{$\sigma$}{sigma}--SPH}
\label{sec3-6:volume_conservation}

When the Jacobian determinant $J$ is not constant, the representative volume for one SPH particle should be updated to preserve the total volume. Of the three proposed methods, Only the $\sigma$-SPH meets this condition. In this section, we explain one of the ways to correct the representative volume using a stabilization term in the Stabilized ISPH method. \par
Since each particle conserves its initial representative volume, a numerical density in the projected space must be $\rho J^0_i/J^N_i$, where $J^N$ is $J$ in the $N$ step. The pressure Poisson equation~(\ref{eq:VCT_ppe}) with volume conservation is therefore as follows:
\begin{equation}
  \label{eq:VCT_ppe_w0}
  \bm c_{\mathrm{trans}}\cdot\aab{\wh{\bm D} p^{N+1}}_i=\cfrac{\rho}{\D t}\,\,\bm J\aab{\wh\N\cdot\bm{v}^*}_i+\gamma\,\cfrac{\rho\,J^0_i/J_i^N
-\aab{\wh\rho^N}_i}{\D t^2}.
\end{equation}
\EDITT{Here, $\aab{\wh\rho^N}$ is the smoothed value computed using the SPH approximation, whereas $J_0$ depends on the initial configuration and may become discontinuous during long-time simulations due to particle motion. Hence, the stabilization term may take excessively large or excessively small values. A previous study\Cite{asai2012stabilized} has shown that when the stabilization parameter $\gamma$ is too large, numerical instability occurs, whereas when it is too small, volume reduction occurs. Therefore, directly using \refE{eq:VCT_ppe_w0} may lead to numerical instability or volume reduction. }\par
\EDITT{In this study, the following equation with $\bar J^N$, a smoothed version of $J^0$, is employed to prevent these issues.}
\begin{equation}
  \label{eq:VCT_ppe_w}
  \bm c_{\mathrm{trans}}\cdot\aab{\wh{\bm D} p^{N+1}}_i=\cfrac{\rho}{\D t}\,\,\bm J\aab{\wh\N\cdot\bm{v}^*}_i+\gamma\,\cfrac{\rho\,\bar J^N_i/J_i^N
-\aab{\wh\rho^N}_i}{\D t^2}.
\end{equation}
The $\bar J^N$ is determined to smooth the volume of the surrounding particles. The volume change $\D V_{ij}$ between the target particle $i$ and the neighboring particle $j$ defined as follows:
\begin{align}
  \quad\D V_{ij}:=C_V \, \cfrac{r_e-|\wh{\bm{r}}_{ij}|}{r_e}\,\,\pab{V_j-V_i
  }=C_V\,\cfrac{r_e-|\wh{\bm{r}}_{ij}|}{r_e}\,\,\pab{\cfrac{J_j^N}{\bar J_j^{N-1}}-\cfrac{J_i^N}{\bar J_i^{N-1}}}V=:\D J_{ij} V,
\end{align}
where $C_V$ is a parameter that regulates the volume smoothing, and $V_i=J^N_i/\bar J^{N-1}_iV$. In this case, the volumes of the $i$ and $j$ particles change as
\begin{align}
  V'_i=V_i+\D V_{ij};\quad V'_j=V_j-\D V_{ij},
\end{align}
where $V'_i=J^N_i/\bar J^N_iV$. Using this relationship, $\bar J^N_i$ and $\bar J^N_j$ can be obtained as
\begin{align}
  \label{eq:VCT_bar_a}
  \bar J^N_i=\cfrac{ J_i^N\bar J^{N-1}_i}{ J^N_i+\bar J^{N-1}_i\D J_{ij}};\quad
  \bar J^N_j=\cfrac{ J_j^N\bar J^{N-1}_j}{ J^N_j-\bar J^{N-1}_j\D J_{ij}}.
\end{align}

\EDITs{\refE{eq:VCT_bar_a} is evaluated sequentially for each neighbor $j$ of target particle $i$, with $\bar{J}^{N-1}$ being updated to $\bar{J}^N$ after each evaluation, prior to proceeding to the next neighbor.}\par
Although the density-based particle shifting method\Cite{morikawa2023corrected,tsuji2024reliable}, \EDITf{Volume Conservation Shifting (VCS)\Cite{khayyer2023enhanced}} and others can be used to control the representative volume, one of the most simple techniques based on the stabilization term of the stabilized ISPH method is introduced here for reference. In addition, when the representative volume change is expressed using the above method, the particle arrangement in the projected space becomes either sparse or dense depending on the volume change. Therefore, using SPH(2), which Tsuji et al.\Cite{tsuji2024reliable} demonstrated to be less affected by the sparseness or denseness of particle arrangements, is expected to enable highly accurate simulations.

%%%%%%%%%%%%%%%%%%%%%%%%%%%%%%%%%%%%%%%%%%%%%%%%%%%%%%%%%%%%%%%%%%%%%
%%  4. V&V for free surface simulation with VCTs
%%%%%%%%%%%%%%%%%%%%%%%%%%%%%%%%%%%%%%%%%%%%%%%%%%%%%%%%%%%%%%%%%%%%%
\section{V\&V for free surface simulation with VCTs}
\label{sec4:VandV}

This section performs validation and verification through 2-D free surface flow simulations using the Stabilized ISPH method with the VCTs presented in \refC{sec2-2:Applications_of_VCTs}. Simulations are conducted with the standard ISPH method and the ISPH(2), as shown in \refT{tab:case_names}. In both simulations, OPS\Cite{khayyer2017comparative}, Dynamic Stabilization\Cite{tsuruta2013short}, and XSPH\Cite{monaghan1994simulating} are applied to stabilize, following the approach of Asai et al.\Cite{asai2023class}. \EDITf{In this study, the shift parameter for the OPS is set to 0.2, and the allowable overlap in the DS is set to 40\% of the particle diameter. The XSPH method is applied only to free-surface particles, with the smoothing parameter for velocity set to $5.0\times10^{-3}$.} The no-slip boundary condition is applied to the solid walls. The gravitational acceleration is $|\bm g|=9.8\times10^2$~cm/s$^2$. The physical properties of water are $\rho=1.0$~g/cm$^3$ and $\nu_{\mathrm{K}}=9.8\times10^{-3}$~cm$^2$/s.

\begin{table}[H]
  \centering
  \caption{Case names in \refC{sec4:VandV}}
  \label{tab:case_names}
  \begin{tabular}{cccc}
  \toprule
  % Case names  & \begin{tabular}[c]{@{}c@{}}Second derivatives\\(viscosity term and PPE)\end{tabular} 
  \textbf{Case names}  
  & Second derivatives 
  & Velocity divergence 
  & Pressure gradient \\ \midrule
  \textbf{ISPH}        
  & 0th-order (\refE{eq:sec_0th}) 
  & 1st-order (\refE{eq:div_1st}) 
  & Summation (\refE{eq:gra_sum}) \\ 
  \textbf{ISPH(2)}     
  & SPH(2) (\refE{eq:sec_SPH2})   
  & 1st-order                     
  & 1st-order  (\refE{eq:gra_1st}) \\ \bottomrule
  \end{tabular}
\end{table}
%%-------------------------------------------------------------------
%%  4. 1. Ellipsoidal particle method (Hydrostatic pressure problem)
%%-------------------------------------------------------------------
\subsection{Hydrostatic pressure problem with the ellipsoidal particle method \EDITs{based on SPH} (E-SPH)}
A hydrostatic pressure problem with a rectangular tank is calculated using \EDITT{the ellipsoidal particle method \EDITs{based on SPH} (E-SPH)}. The water height and width are 50~cm each in the projected space. The computational conditions are time step width~$\D t=1.0\times10^{-3}$~s, and initial particle spacing~$d_0=0.5$~cm. In this section, the \EDITs{verification} is carried out for the vertical scale factor~$\alpha$ in the range of 0.5 to 2.5. That means the rectangular tank changes from 25~cm ($\alpha$=0.5) to 125~cm ($\alpha$=2.5).\par
\begin{figure}[H]
  \centering
  \begin{subfigure}[b]{0.49\linewidth}
      \centering
      \includegraphics[width=0.94\linewidth]{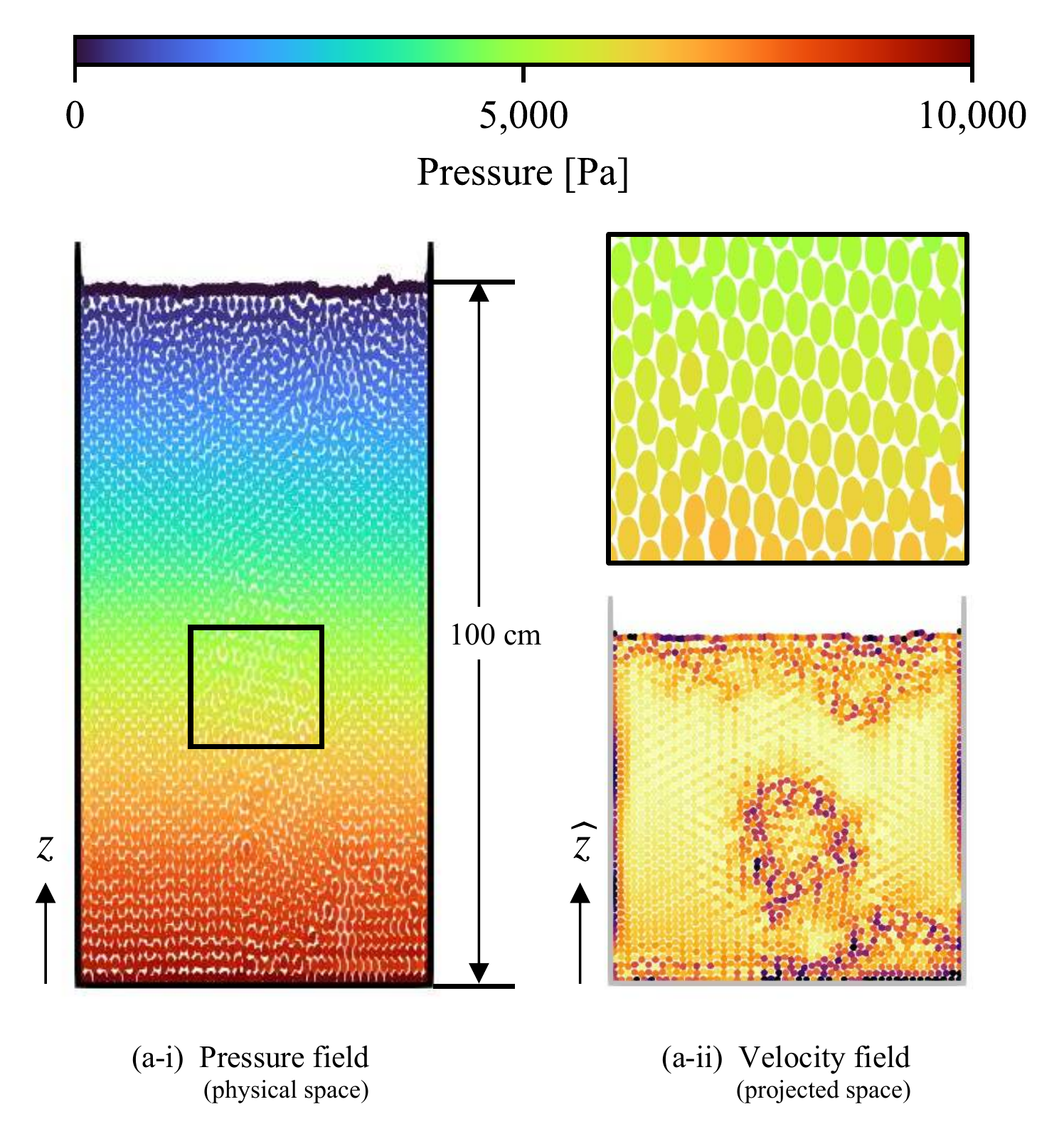}
      \caption{ISPH}
  \end{subfigure}
  \begin{subfigure}[b]{0.49\linewidth}
      \centering
      \includegraphics[width=0.875\linewidth]{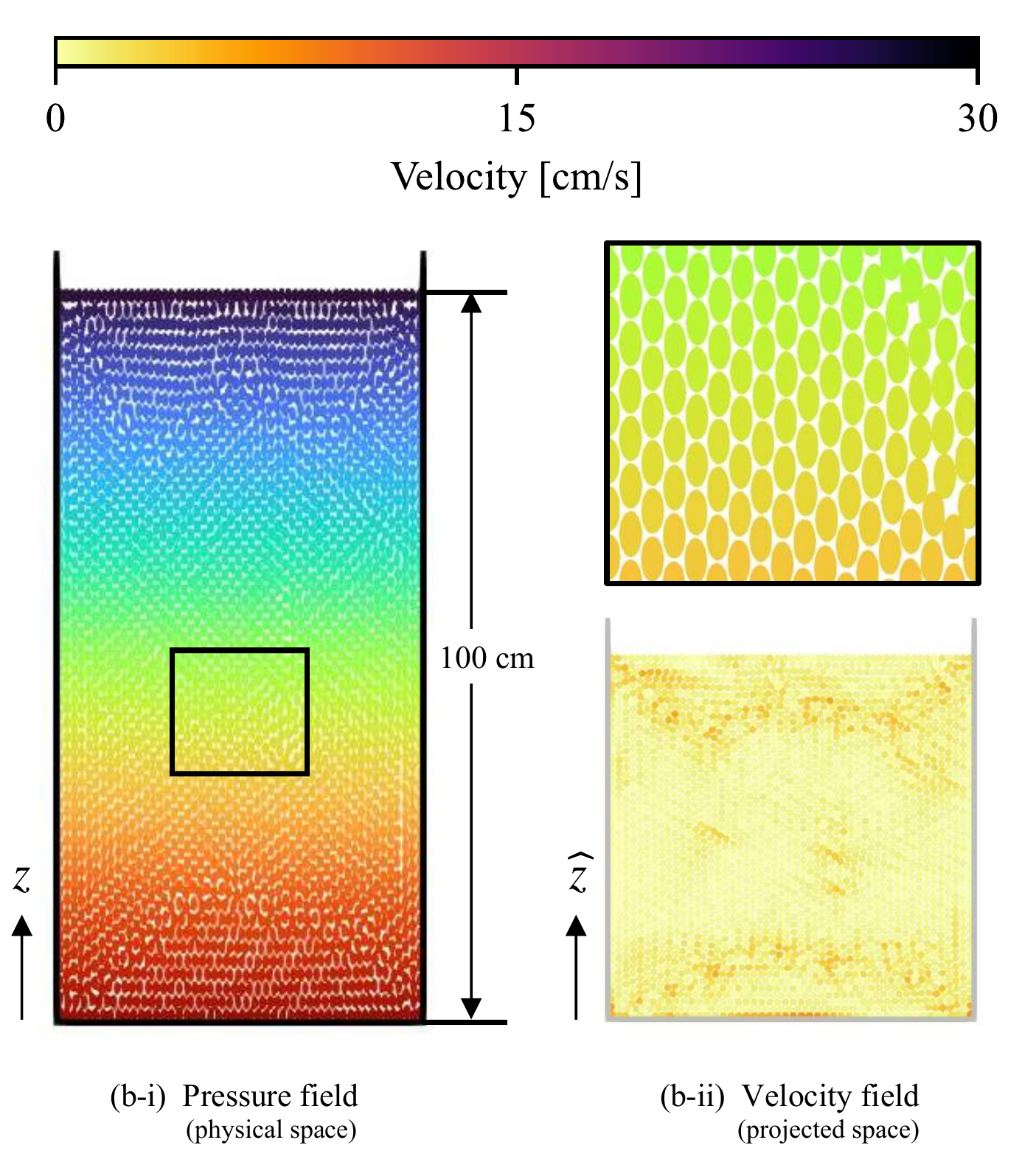}
      \caption{ISPH(2)}
  \end{subfigure}
  \caption{Pressure and velocity fields at 10~s calculated by the E-SPH for the hydrostatic pressure problem at $\alpha=0.5$}\label{fig:s4_e_0.5}
\end{figure}
\begin{figure}[H]
  \centering
  \begin{subfigure}[b]{0.49\linewidth}
      \centering
      \includegraphics[width=0.9247\linewidth]{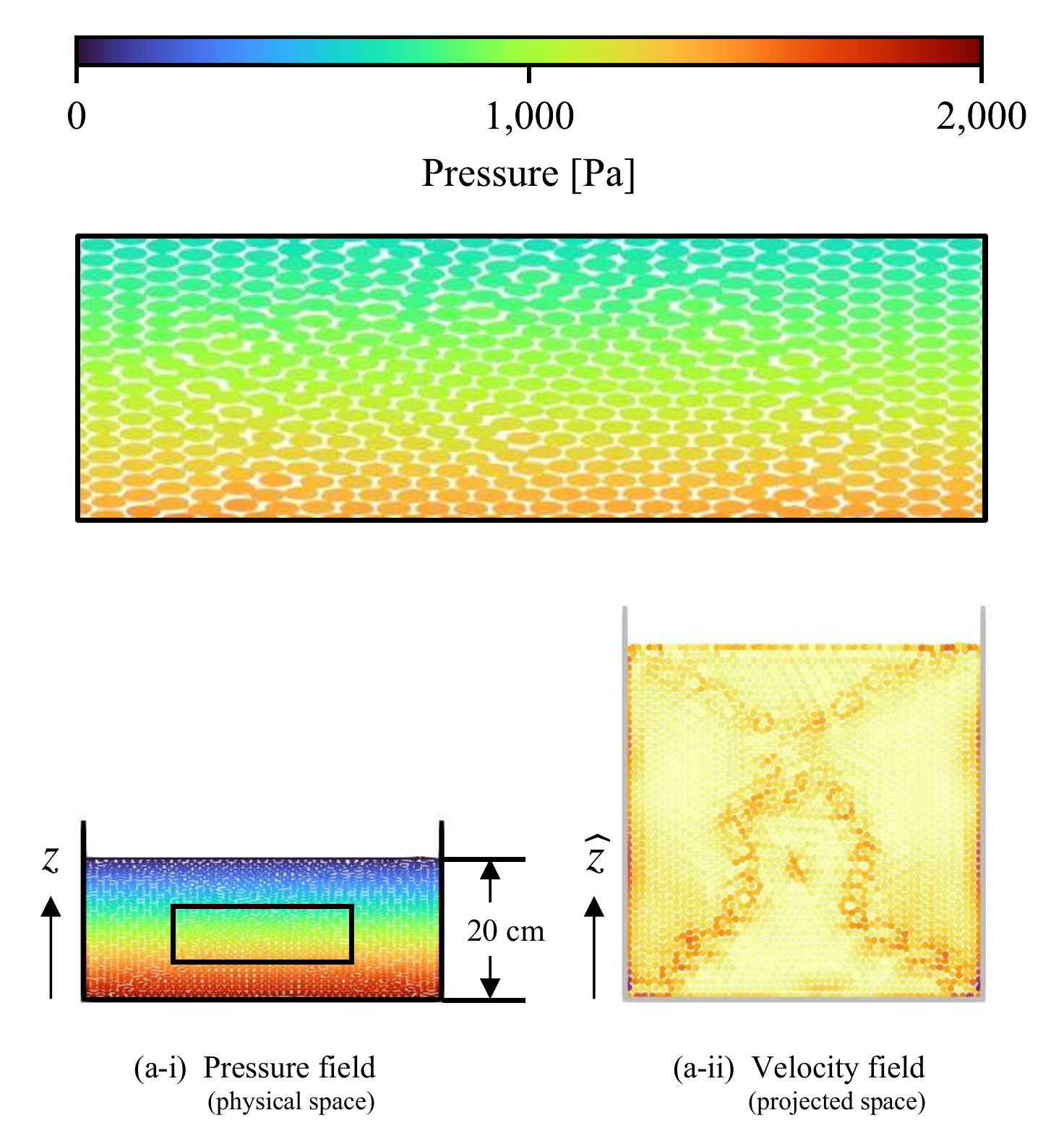}
      \caption{ISPH}
  \end{subfigure}
  \begin{subfigure}[b]{0.49\linewidth}
      \centering
      \includegraphics[width=0.875\linewidth]{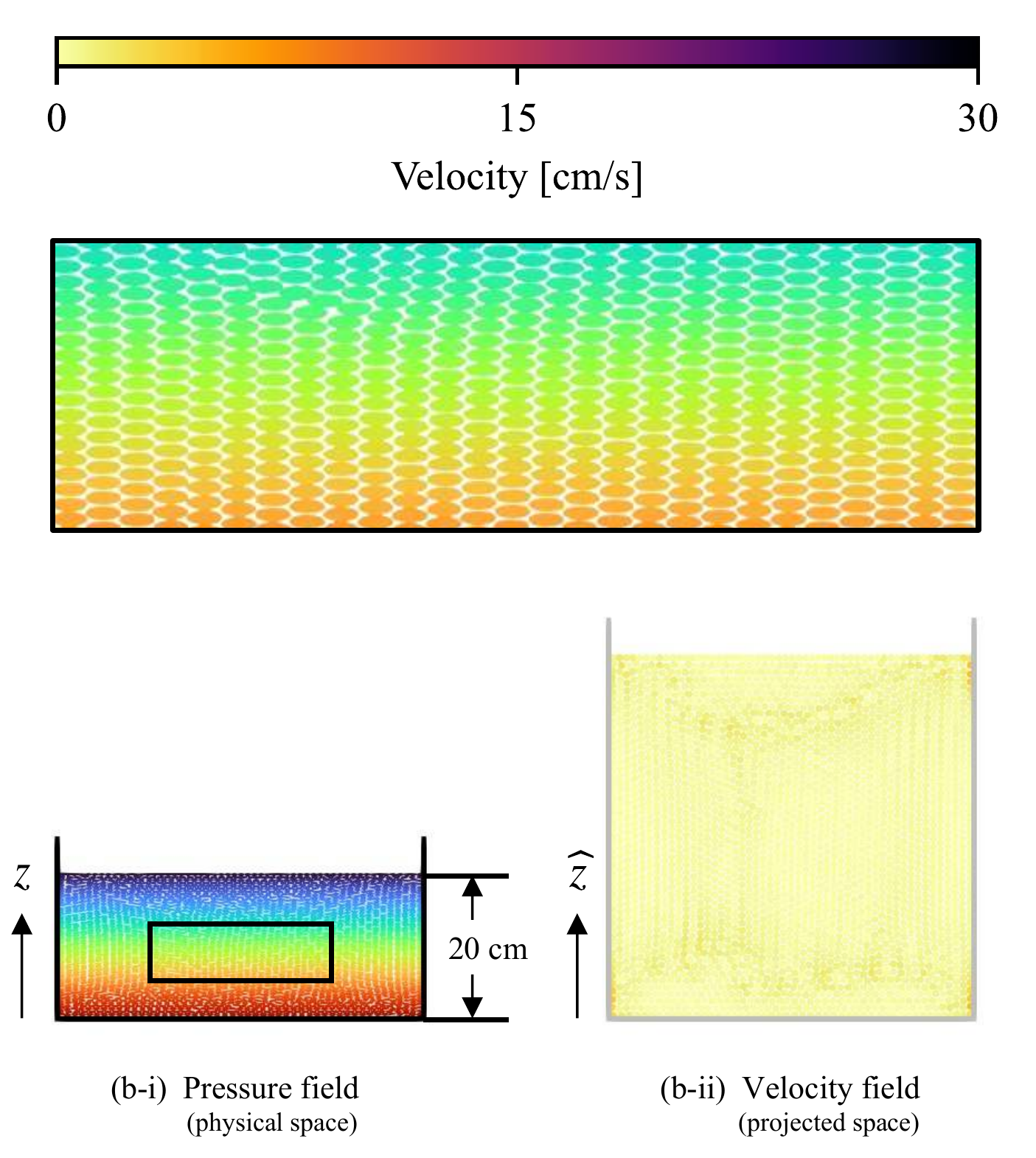}
      \caption{ISPH(2)}
  \end{subfigure}
  \caption{Pressure and velocity fields at 10~s calculated by the E-SPH for the hydrostatic pressure problem at $\alpha=2.5$}\label{fig:s4_e_2.5}
\end{figure}
\begin{figure}[H]  
    \centering
    \begin{subfigure}[b]{0.49\linewidth}
      \centering
      \includegraphics[width=\linewidth]{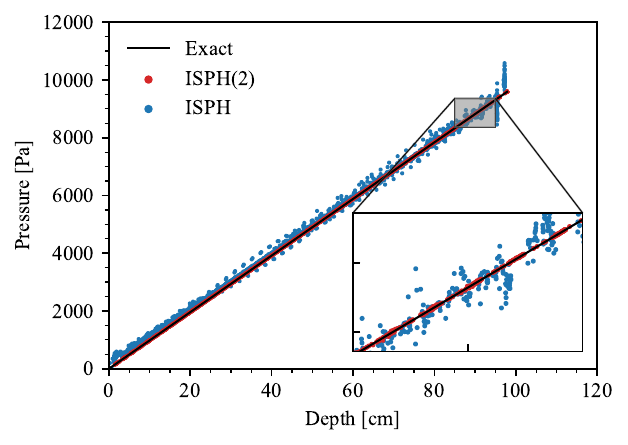}
      \subcaption{$\alpha=0.5$}          
    \end{subfigure} 
    \begin{subfigure}[b]{0.49\linewidth}
      \centering
      \includegraphics[width=0.947\linewidth]{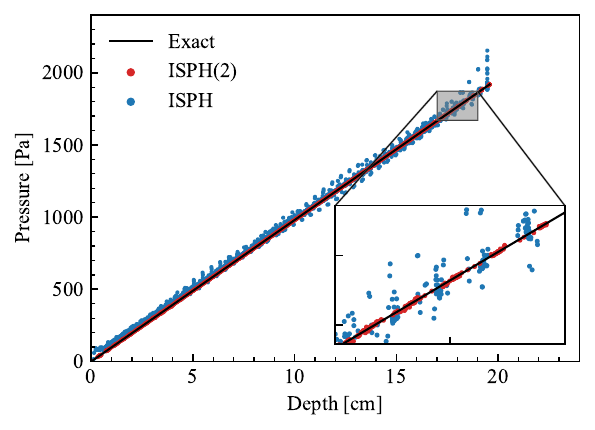}
      \subcaption{$\alpha=2.5$}
    \end{subfigure}
  \caption{Pressure profiles at 10~s calculated using the E-SPH at $\alpha=0.5,\,2.5$}\label{fig:s4_e_pres}
\end{figure}
\begin{figure}[H]
  \centering
    \centering
    \includegraphics[width=0.482\linewidth]{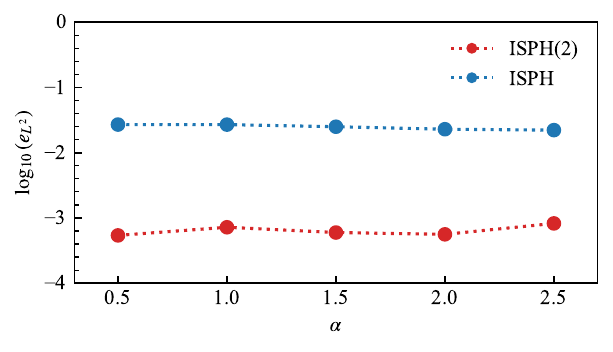}
  \caption{Logarithmic relative $L^2$ errors in pressure at 10~s obtained using the E-SPH for each $\alpha$}\label{fig:s4_e_l2}
\end{figure}

\refF{fig:s4_e_0.5} and \ref{fig:s4_e_2.5} show the pressure and velocity fields at 10~s calculated by E-SPH based on the conventional ISPH and ISPH(2) with $\alpha=0.5$ and $\alpha=2.5$, respectively. In both cases, a coordinate transformation is applied to expand or contract the projected space vertically. The pressure distribution in the ISPH shows significant spatial fluctuations, and the maximum value of non-physical velocity reaches 30~cm/s. In contrast, ISPH(2) exhibits a smoother pressure distribution, and the non-physical velocity is suppressed ideally.
\refF{fig:s4_e_pres} presents the pressure distribution along the vertical axis, where the depth is defined as $z_{\mathrm{FS}} - z_i$, with $z_{\mathrm{FS}}$ representing the average $z$-coordinate of the free surface particles. It is seen that, compared to ISPH, ISPH(2) agrees well with the exact solution of hydrostatic pressure.
\refF{fig:s4_e_l2} shows the logarithmic relative $L^2$ errors $e_{L^2}$ in pressure at 10~s for all $\alpha$ values. It is observed that the relative $L^2$ errors in ISPH(2) are two orders smaller than in ISPH for all $\alpha$ values.


%%-------------------------------------------------------------------
%%  4. 2. Bottom boundary-fitted particle method (Hydrostatic pressure problem)
%%-------------------------------------------------------------------
\subsection{Hydrostatic pressure problem with the bottom boundary-fitted particle method (BF-SPH)}
A hydrostatic pressure problem in a tank with a cosine-shaped bottom is calculated using the bottom boundary-fitted particle method (BF-SPH). The bottom is defined by $h(x)=R_{\cos}\cos(2\pi  x/50)$~cm. The water height and width are each 50~cm in physical space. The computational conditions are time step width~$\D t=1.0\times10^{-3}$~s, and initial particle spacing~$d_0=0.5$~cm. In this section, the \EDITs{verification} is carried out for the coefficient $R_{\cos}$ that defines the bottom shape in the range of 1.0 to 5.0~cm.\par
\begin{figure}[H]
  \centering
      \begin{subfigure}[b]{0.49\linewidth}
          \centering
          \includegraphics[width=0.92\linewidth]{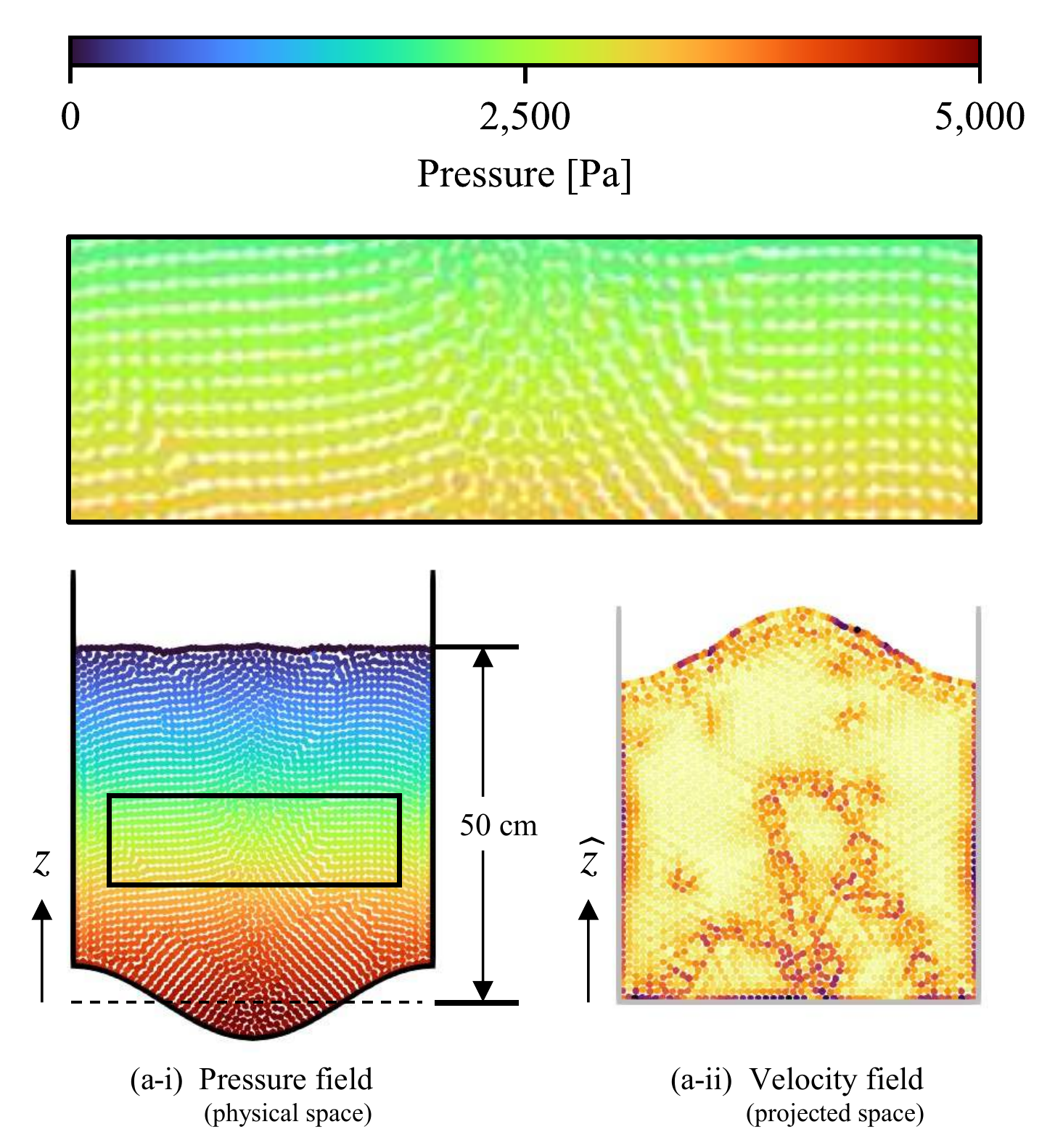}
          \caption{ISPH}
      \end{subfigure}
      \begin{subfigure}[b]{0.49\linewidth}
          \centering
          \includegraphics[width=0.875\linewidth]{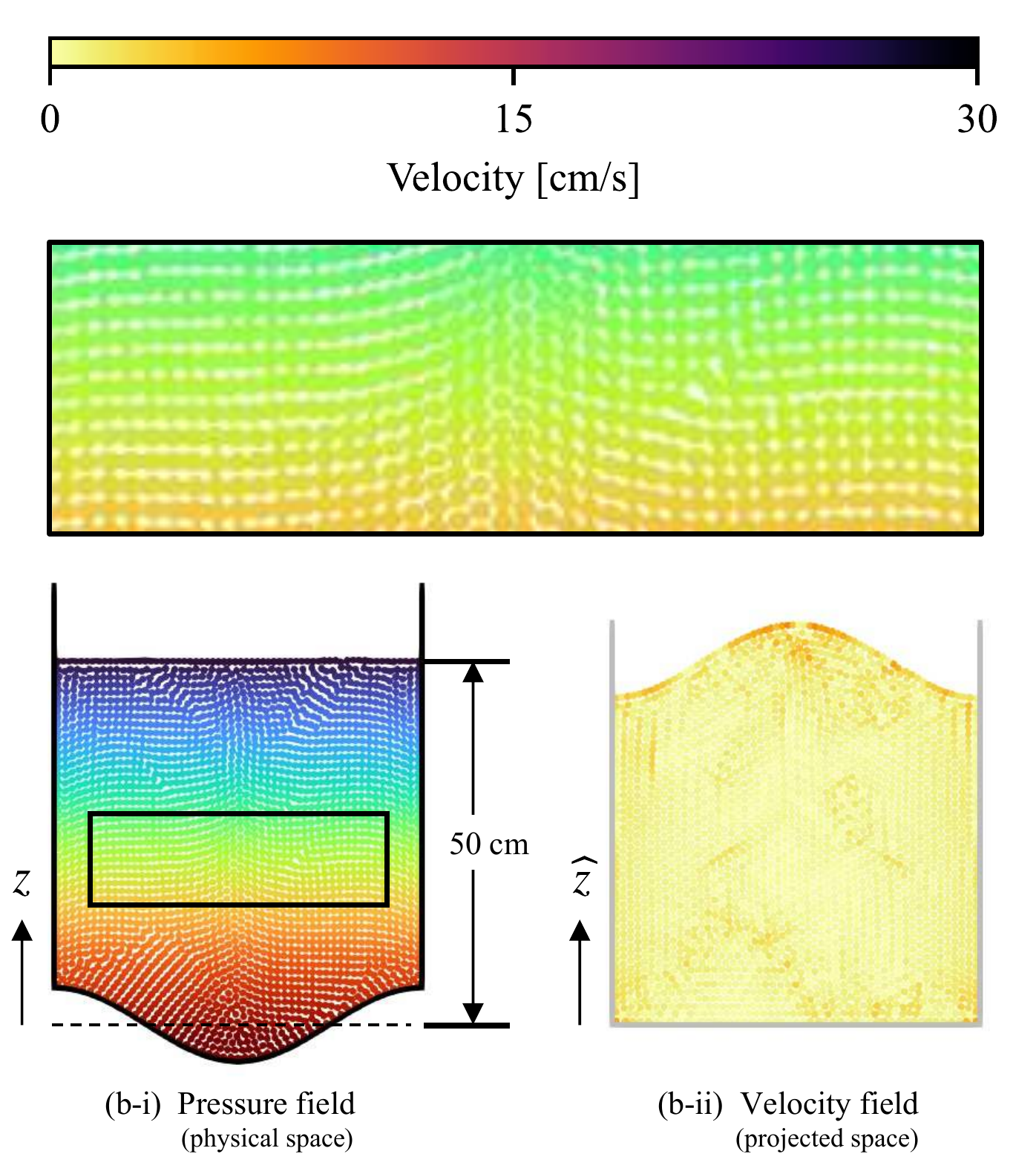}
          \caption{ISPH(2)}
      \end{subfigure}
  \caption{Pressure and velocity fields at 10~s calculated with the BF-SPH for the hydrostatic pressure problem at $R_{\cos}$}\label{fig:s4_b_5.0}
\end{figure}
\begin{figure}[H]  
  \centering
  \includegraphics[width=0.4998\linewidth]{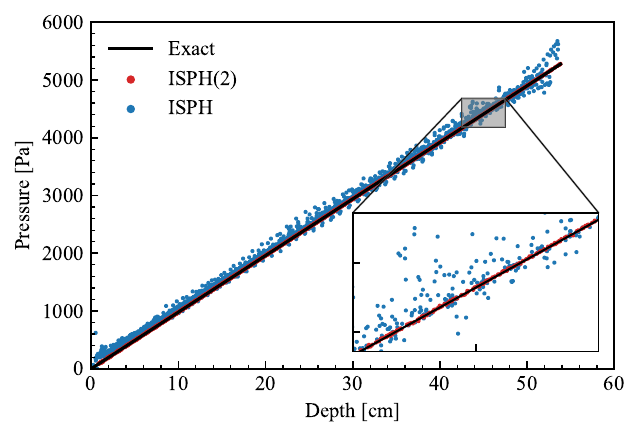}
  \caption{Pressure profiles at 10~s calculated using the BF-SPH at $R_{\cos}=5.0$}\label{fig:s4_b_pres}
\end{figure}
\begin{figure}[H]
  \centering
    \centering
    \includegraphics[width=0.482\linewidth]{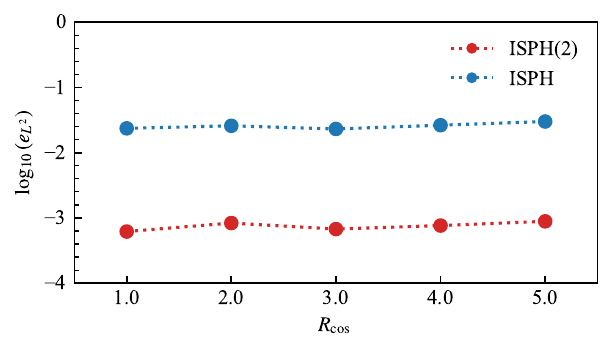}
  \caption{Logarithmic relative $L^2$ errors in pressure at 10~s obtained using the BF-SPH at each ${R_\cos}$}\label{fig:s4_b_l2}
\end{figure}

\refF{fig:s4_b_5.0} shows the pressure and velocity fields at 10~s calculated by the method, compared with conventional ISPH and ISPH(2) at $R_{\cos}=5.0$~cm. In both cases, a coordinate transformation is applied to make the bottom flat in the projected space. The pressure distribution in the ISPH shows significant spatial fluctuations, and the maximum value of non-physical velocity reaches 30~cm/s. In contrast, ISPH(2) exhibits smoother pressure distribution and suppresses the non-physical velocity.
\refF{fig:s4_b_pres} presents the pressure distribution along the vertical axis, where the depth is defined as $z_{\mathrm{FS}} - z_i$, with $z_{\mathrm{FS}}$ representing the average $z$-coordinate of the free surface particles. It is seen that, compared to ISPH, ISPH(2) agrees well with the exact solution of hydrostatic pressure.
\refF{fig:s4_b_l2} shows the logarithmic relative $L^2$ errors $e_{L^2}$ in pressure at 10~s for all $R_{\cos}$ cases. It is observed that the relative $L^2$ errors in ISPH(2) are two orders smaller than in ISPH for all $R_{\cos}$ cases. 


%%-------------------------------------------------------------------
%%  4. 3. Bottom boundary-fitted ellipsoidal particle method (Dam break problem)
%%-------------------------------------------------------------------
\subsection{Dam break problem with the bottom boundary-fitted ellipsoidal particle method (BFE-SPH)}
\begin{figure}[H]
  \centering
  \includegraphics[width=0.9856\linewidth]{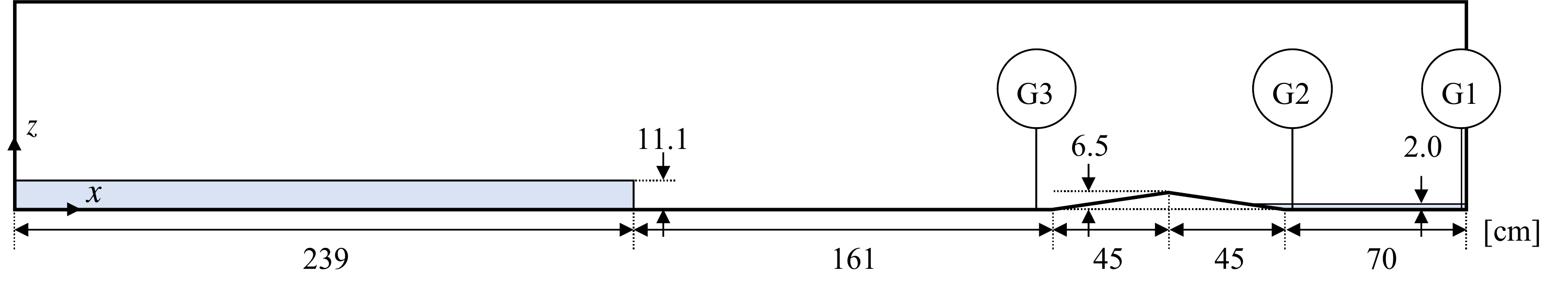}
  \caption{Dam break problem setup for comparing computational cost}\label{fig:s4_be_model}
\end{figure}

A dam break problem resulting in flow over a triangular bump is calculated using the bottom boundary-fitted ellipsoidal particle method (BFE-SPH). The setup for the dam-break problem is shown in \refF{fig:s4_be_model}. This validation compares the experiment results from\Cite{soares2007experiments} with the numerical simulation results. In the experiment, water levels were measured sequentially for 45~s in real-time on three water level gauges: G1, G2, and G3. The $x$-coordinate values of Gauge 1, 2, and 3 are 557.5~cm, 492.5~cm, and 393.5~cm, respectively. The time step width~$\D t$ is set as $d_0/0.4\times10^{-3}$~s. In ISPH, calculations are performed with a vertical scale factor of $\alpha = 2.0$ and initial particle spacing of $d_0 = 0.2$~cm, while in ISPH(2), $\alpha$ is set to 1.0 (BF-SPH) and 2.0 (BFE-SPH), and $d_0$ to 0.4 and 0.2~cm.\par
\begin{figure}[H]
  \centering
  \begin{subfigure}[b]{\linewidth}
      \centering
      \includegraphics[width=\linewidth]{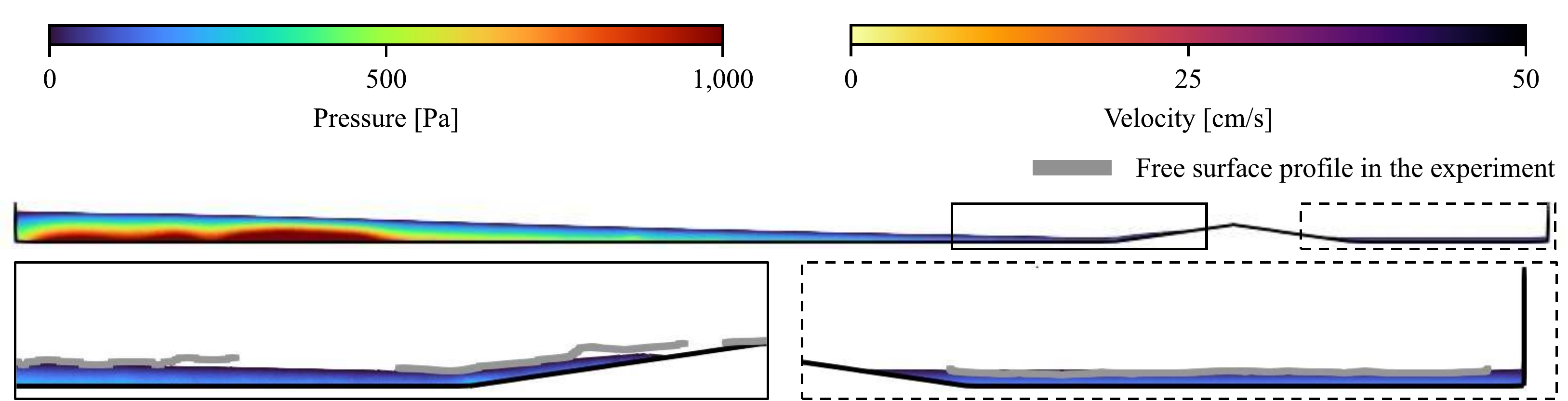}
      \caption{1.8~s}\label{fig:s4_be_w_a}
  \end{subfigure}
  \begin{subfigure}[b]{\linewidth}
      \centering
      \includegraphics[width=\linewidth]{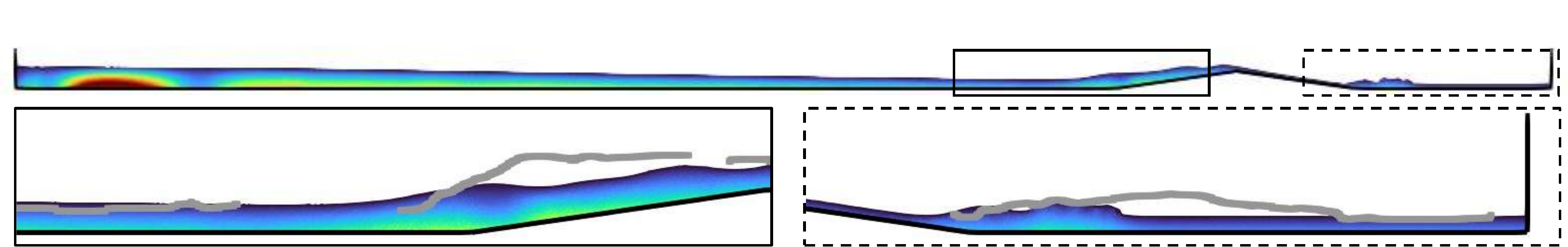}
      \caption{3.0~s}\label{fig:s4_be_w_b}
  \end{subfigure}
  \begin{subfigure}[b]{\linewidth}
    \centering
    \includegraphics[width=\linewidth]{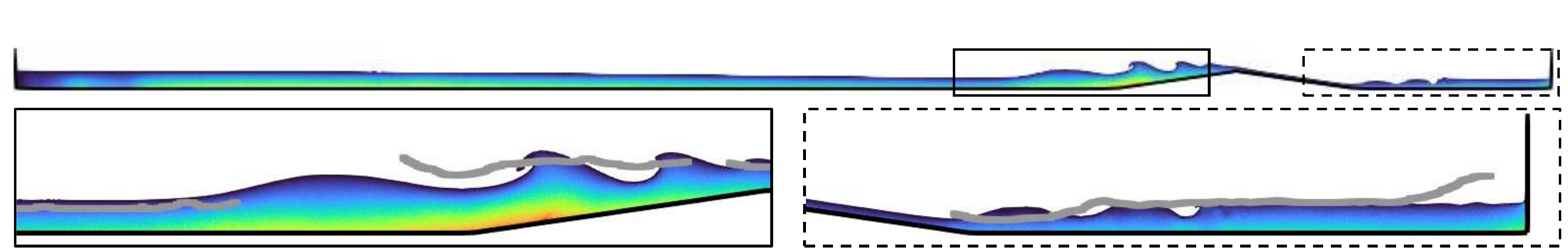}
    \caption{3.7~s}\label{fig:s4_be_w_c}
  \end{subfigure}
  \begin{subfigure}[b]{\linewidth}
    \centering
    \includegraphics[width=\linewidth]{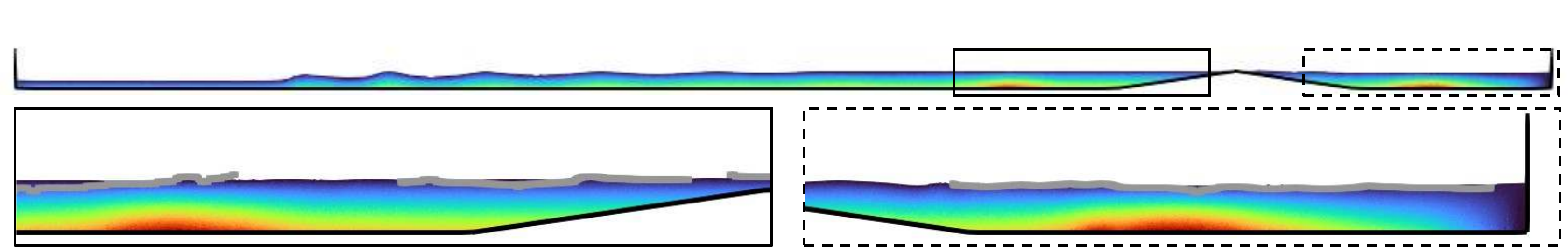}
    \caption{8.4~s}\label{fig:s4_be_w_d}
  \end{subfigure}
  \begin{subfigure}[b]{\linewidth}
    \centering
    \includegraphics[width=\linewidth]{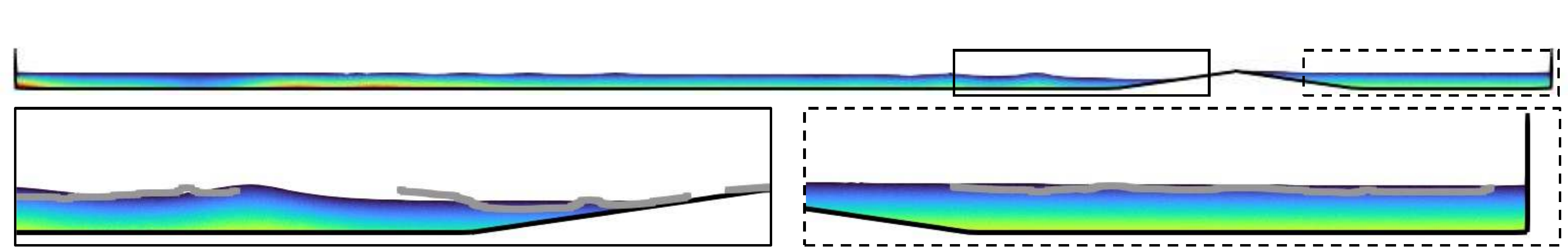}
    \caption{15.5~s}\label{fig:s4_be_w_e}
  \end{subfigure}
  \begin{subfigure}[b]{\linewidth}
    \centering
    \includegraphics[width=\linewidth]{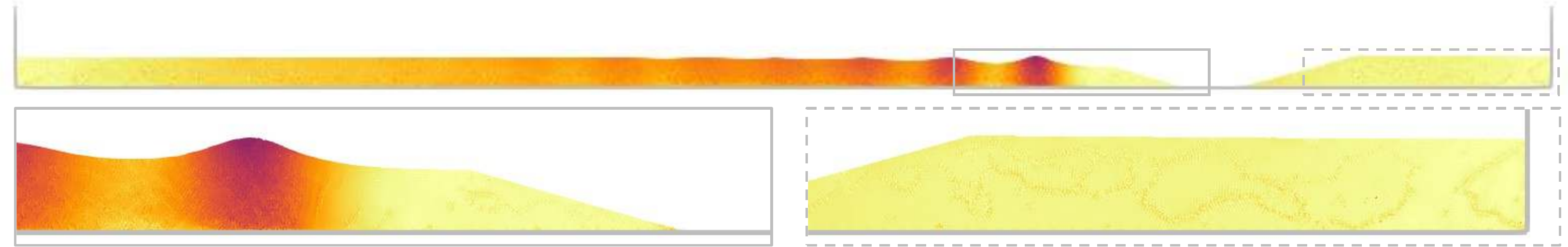}
    \caption{15.5~s (velocity field in projected space)}\label{fig:s4_be_w_f}
  \end{subfigure}
  \caption{Pressure and velocity fields calculated with the BFE-SPH ($\alpha=2.0$, $d_0=0.2$) with ISPH for the dam break problem}\label{fig:s4_be_w}
\end{figure}
\begin{figure}[H]
  \centering
  \begin{subfigure}[b]{\linewidth}
      \centering
      \includegraphics[width=\linewidth]{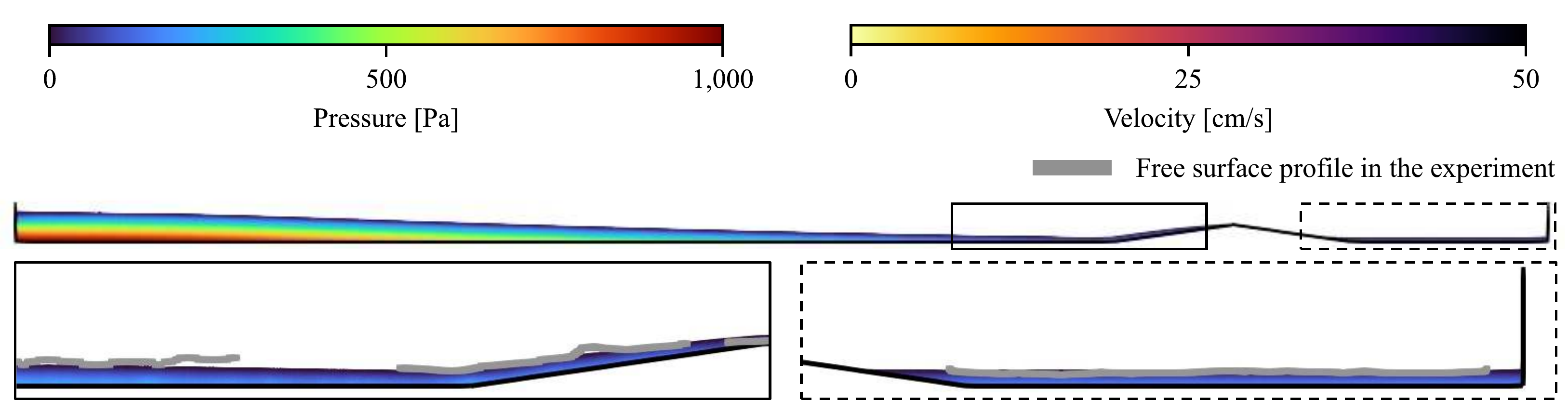}
      \caption{1.8~s}
  \end{subfigure}
  \begin{subfigure}[b]{\linewidth}
      \centering
      \includegraphics[width=\linewidth]{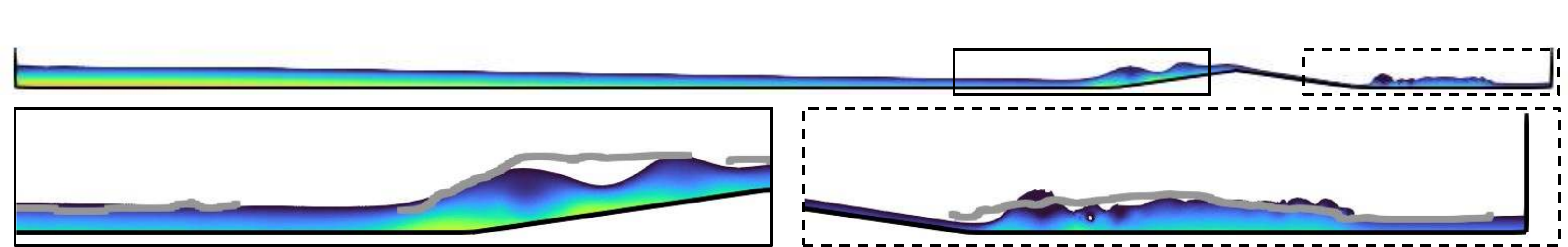}
      \caption{3.0~s}
  \end{subfigure}
  \begin{subfigure}[b]{\linewidth}
    \centering
    \includegraphics[width=\linewidth]{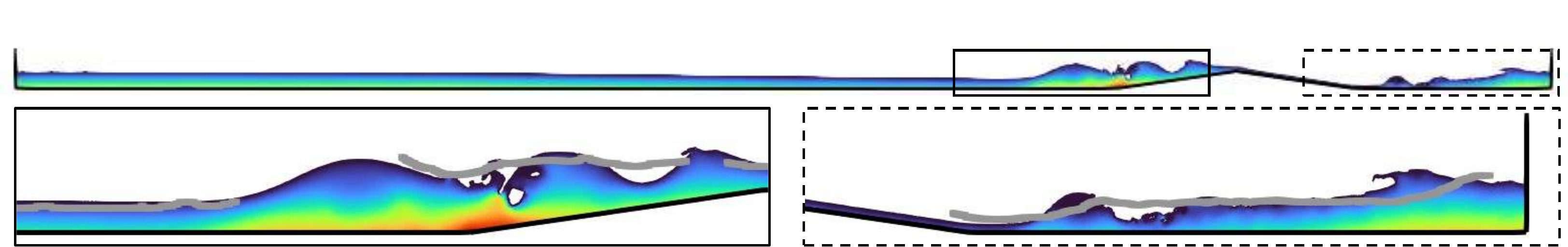}
    \caption{3.7~s}
  \end{subfigure}
  \begin{subfigure}[b]{\linewidth}
    \centering
    \includegraphics[width=\linewidth]{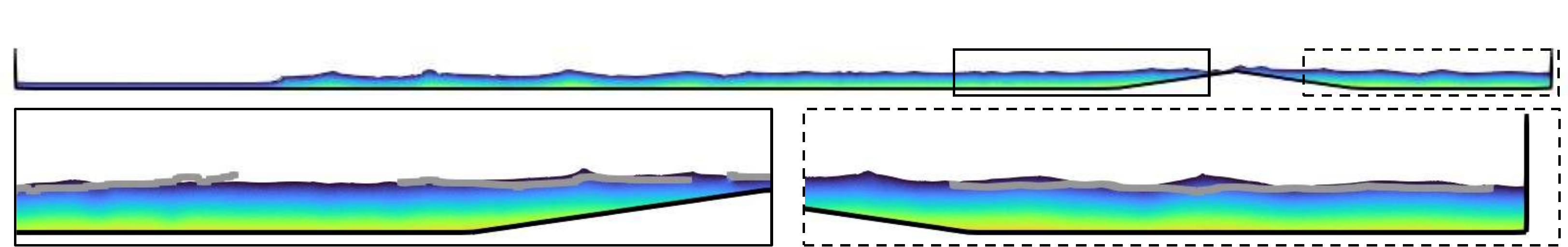}
    \caption{8.4~s}
  \end{subfigure}
  \begin{subfigure}[b]{\linewidth}
    \centering
    \includegraphics[width=\linewidth]{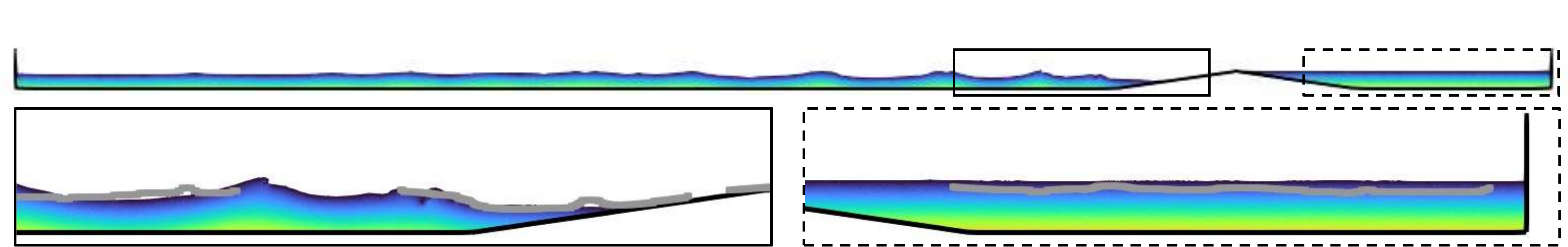}
    \caption{15.5~s}
  \end{subfigure}
  \begin{subfigure}[b]{\linewidth}
    \centering
    \includegraphics[width=\linewidth]{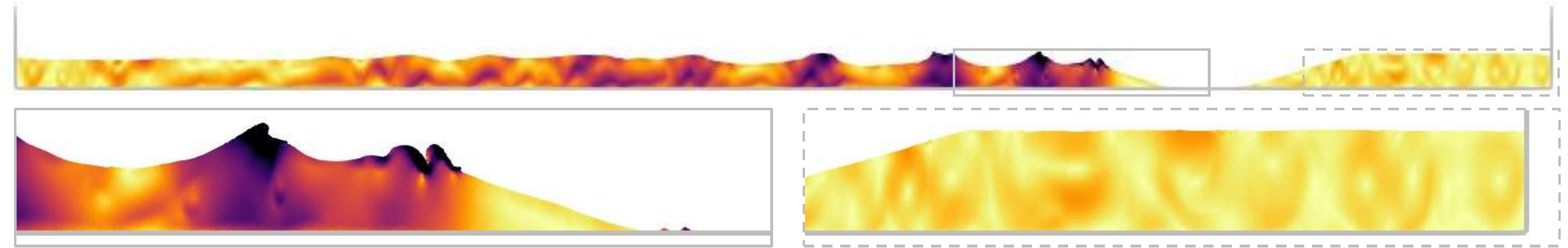}
    \caption{15.5~s (velocity field in projected space)}\label{fig:s4_be_2_f}
  \end{subfigure}
  \caption{Pressure and velocity fields calculated with the BFE-SPH ($\alpha=2.0$, $d_0=0.2$) with ISPH(2) for the dam break problem}\label{fig:s4_be_2}
\end{figure}

\refF{fig:s4_be_w} and \refF{fig:s4_be_2} show the pressure and velocity fields calculated using ISPH and ISPH(2) with $\alpha=2.0$ and $d_0=0.2$~cm, respectively. The gray lines indicate the free surface profile in the experiment. In the experiment, it was observed that the water flowed over the dry channel, and upon reaching the bump, part of the wave was reflected, forming a bore that traveled back upstream. In contrast, the other part moved over the bump, leading to wave propagation on an upward dry slope ($t=1.8$~s). After passing over the bump, the water flowed down the dry slope until it reached a resting pool of water, where the rapid wavefront slowed abruptly, generating a bore that traveled downstream ($t=3.0$~s). The bore reflected off the downstream wall, traveling back toward the bump ($t=3.7$~s), but the water was initially unable to cross the crest. A second reflection against the downstream wall was necessary for the wave to cross the bump and travel upstream again ($t=8.4$~s). Multiple reflections of the flow were observed against both the bump and the channel ends ($t=15.5$~s). 
\EDITs{\refF{fig:s4_be_w} with ISPH shows non-physical pressure fields, especially in the region from $x=0$ to $x=200$ and in the enlarged area presented in \hyperref[fig:s4_be_w_d]{Fig.~\ref*{fig:s4_be_w}(d)}. In contrast, \refF{fig:s4_be_2} with ISPH(2) shows a smooth and physically consistent pressure distribution across the entire domain. Furthermore, when comparing the free surface profiles recorded in the experiment with the numerical results, the ISPH method shows a lower free surface position than the experimental data at 3.0~s (\hyperref[fig:s4_be_w_b]{Fig.~\ref*{fig:s4_be_w}(b)}) and 3.7~s (\hyperref[fig:s4_be_w_c]{Fig.~\ref*{fig:s4_be_w}(c)}). In contrast, the ISPH(2) method provides better agreement with the experimental free surface profile. In addition, the ISPH seems to include too much artificial damping shown in the velocity field of \hyperref[fig:s4_be_w_f]{Fig.~\ref*{fig:s4_be_w}(f)}, although the ISPH(2) can suppress it as shown in \hyperref[fig:s4_be_2_f]{Fig.~\ref*{fig:s4_be_2}(f)}.}\par

\begin{figure}[H]
  \centering
  \begin{subfigure}[b]{\linewidth}
      \centering
      \includegraphics[width=0.911\linewidth]{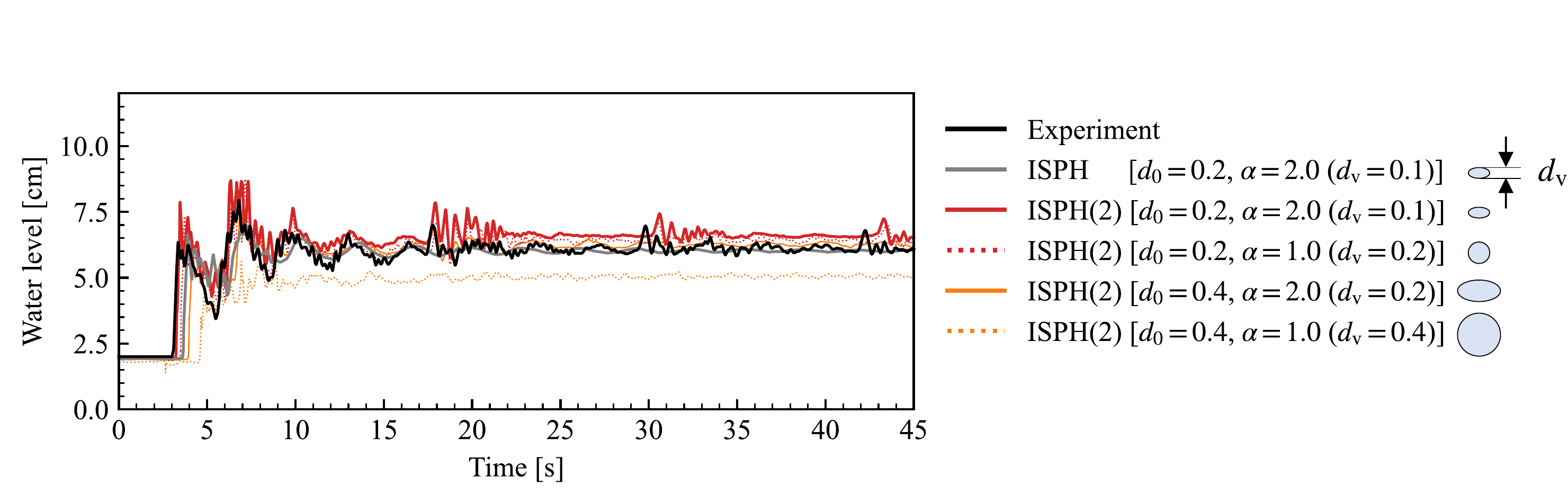}
      \caption{Gauge 1 ($x=557.5$~cm)}\label{fig:s4_be_comp1}
  \end{subfigure}
  \begin{subfigure}[b]{\linewidth}
      \centering
      \includegraphics[width=0.911\linewidth]{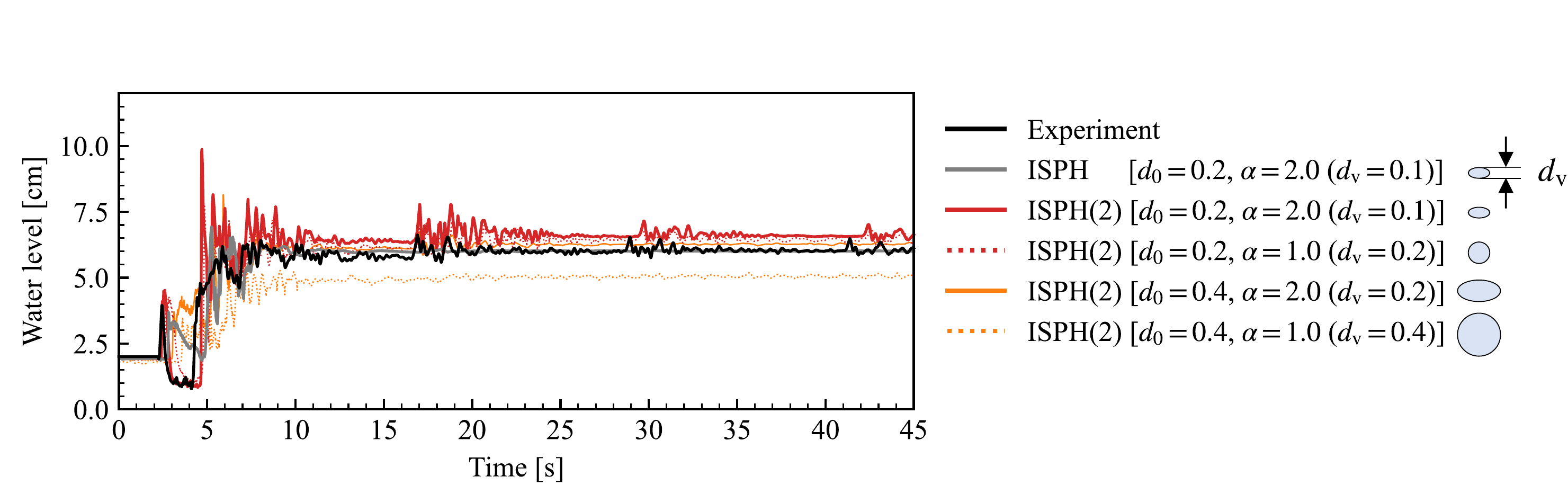}
      \caption{Gauge 2 ($x=492.5$~cm)}\label{fig:s4_be_comp2}
  \end{subfigure}
  \begin{subfigure}[b]{\linewidth}
    \centering
    \includegraphics[width=0.911\linewidth]{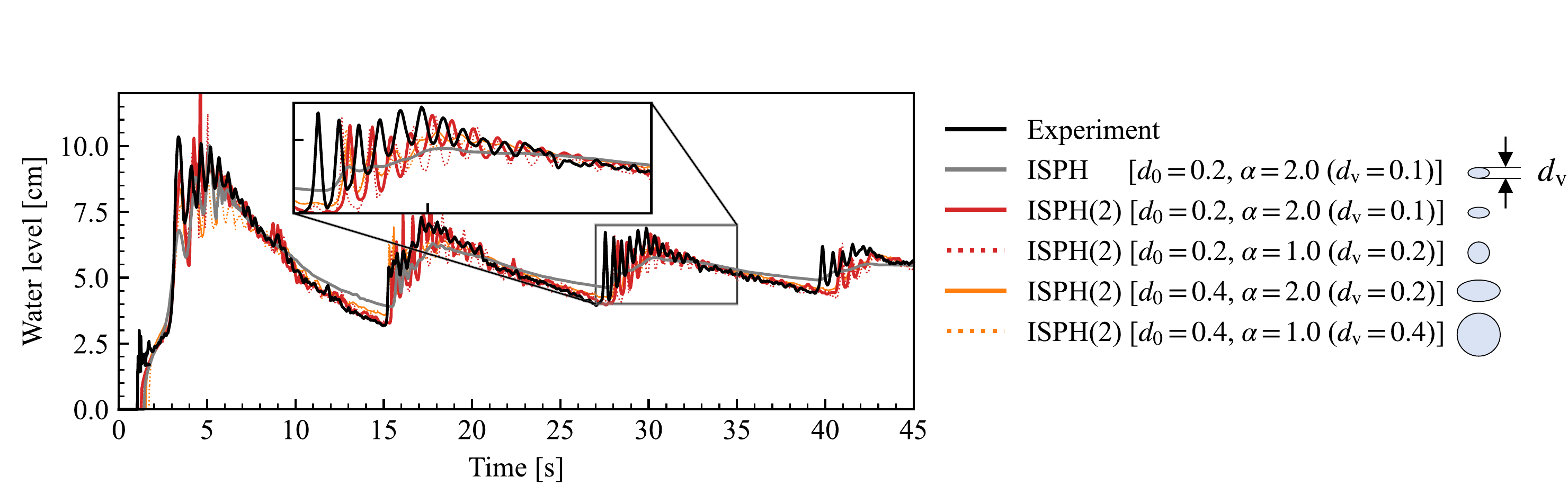}
    \caption{Gauge 3 ($x=393.5$~cm)}\label{fig:s4_be_comp3}
  \end{subfigure}
  \caption{Time histories of the water levels obtained using the BF-SPH and BFE-SPH at each gauge}\label{fig:s4_be_comp}
\end{figure}

\begin{figure}[H]
    \centering
    \includegraphics[width=\linewidth]{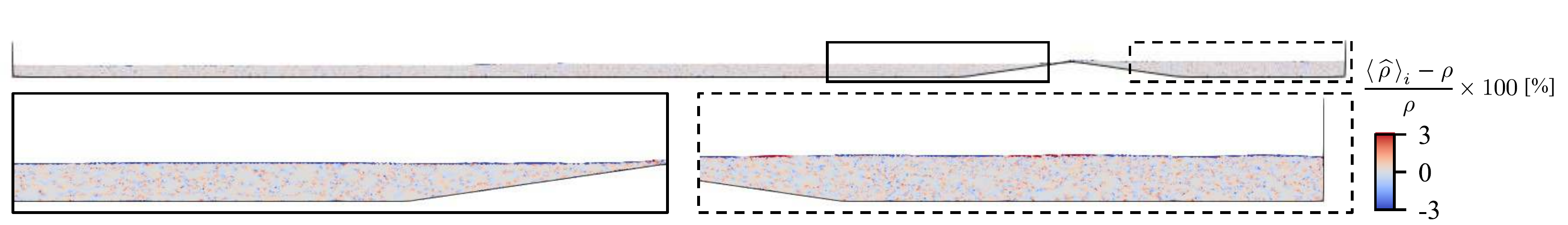}
  \caption{\EDITf{Density error field at 45~s calculated with the BFE-SPH ($\alpha=2.0$, $d_0=0.2$) with ISPH(2) for the dam break problem}}\label{fig:s4_be_dens}
\end{figure}

\refF{fig:s4_be_comp} shows the time histories of the water levels on Gauges 1, 2, and 3 ($x=557.5$, 492.5, and 393.5~cm). 
The high-frequency disturbances observed at each gauge are due to the soliton fission, i.e., disturbances formed by solitary waves that disperse as they move upslope.
In the case of ISPH (gray line), the flow is significantly attenuated compared to the experimental results for each gauge, especially for Gauge 3, where no high-frequency disturbances are calculated after 10~s, as shown in the enlarged figure. In contrast, for ISPH(2), the flow does not attenuate as compared to ISPH and converges to the experimental results with increasing vertical resolution. In this figure, $d_0$ indicates the particle distance in the projected space, and $d_{\mathrm{v}}$ means the effective vertical resolution in the real space. In particular, the finest resolution model with ISPH(2) (red line) shows the sinking and subsequent rise of the water level at Gauge 2 from 3 to 5~s, as shown in \hyperref[fig:s4_be_comp2]{Fig.~\ref*{fig:s4_be_comp}(b)}. However, the water level at 45~s on the left side of the triangular bump (Gauge 3) is generally consistent with the experiment, while the water level on the right side (Gauges 1 and 2) is higher than in the experiment.  \refF{fig:s4_be_dens} shows the density error field at 45~s, indicating that the density error is very small and that volume is conserved. 
Therefore, it is unlikely that volume conservation is the cause of the difference in water level error on the left and right sides of the triangular bump, and it may be that a 2-D simulation has been performed, resulting in a slight tendency for water to accumulate on the right side of the triangular bump. It may be necessary to perform a 3-D analysis of the wave-breaking phenomenon, as this shows that it can be improved by performing a 3-D simulation, see, for example,\Cite{biscarini2010cfd}.
However, it was confirmed that the calculations were performed with a high degree of reproducibility of the water level compared to the results, for example, when the same problem was solved with the shallow water equation\Cite{lavoie2017comparison}.

\begin{figure}[H]
  \centering
  \begin{subfigure}[b]{\linewidth}
      \centering
      \includegraphics[width=0.74\linewidth]{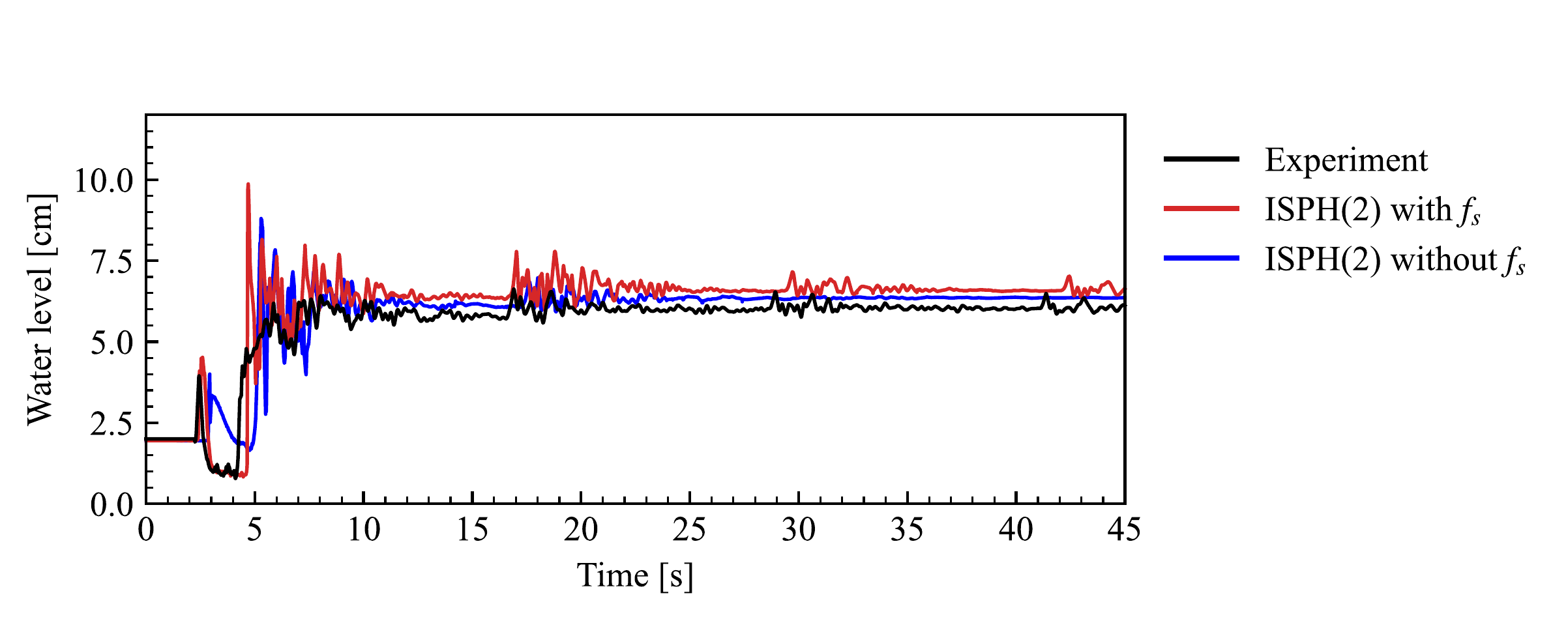}
      \caption{Gauge 1 ($x=557.5$~cm)}
  \end{subfigure}
  \begin{subfigure}[b]{\linewidth}
      \centering
      \includegraphics[width=0.74\linewidth]{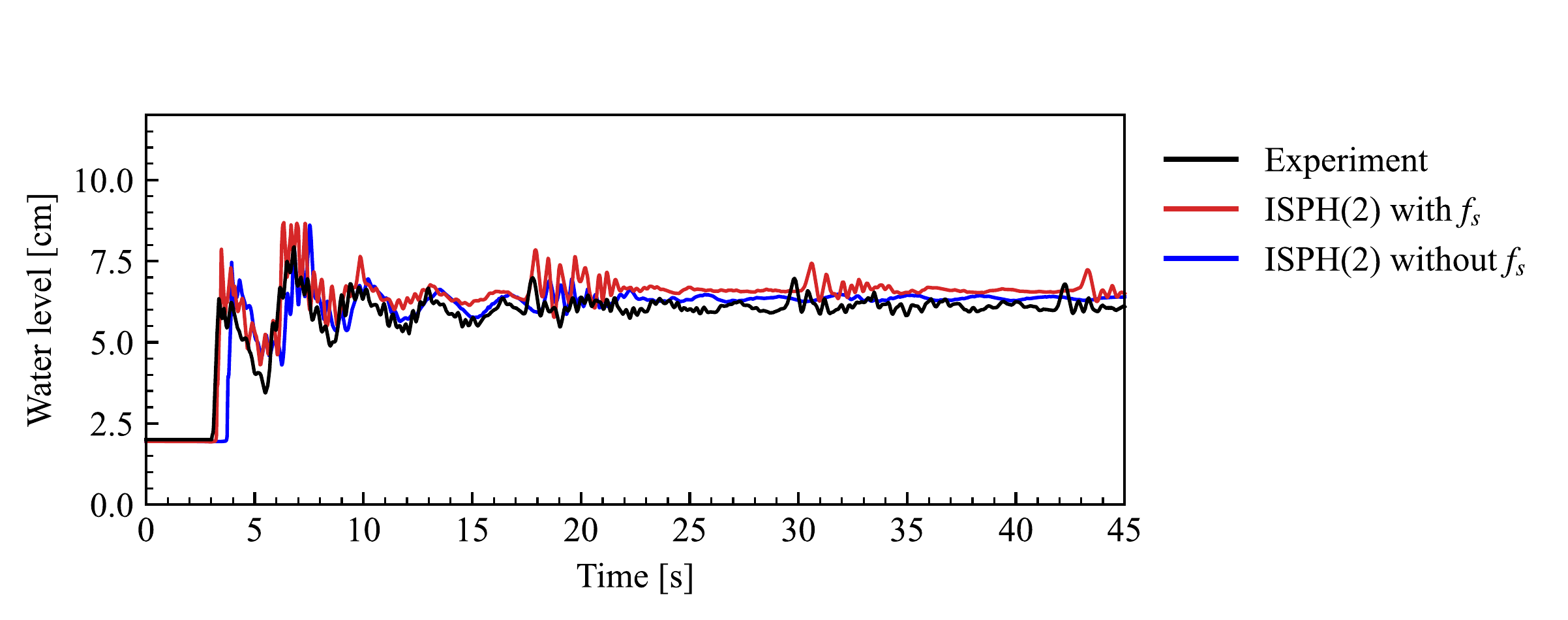}
      \caption{Gauge 2 ($x=492.5$~cm)}
  \end{subfigure}
  \begin{subfigure}[b]{\linewidth}
    \centering
    \includegraphics[width=0.74\linewidth]{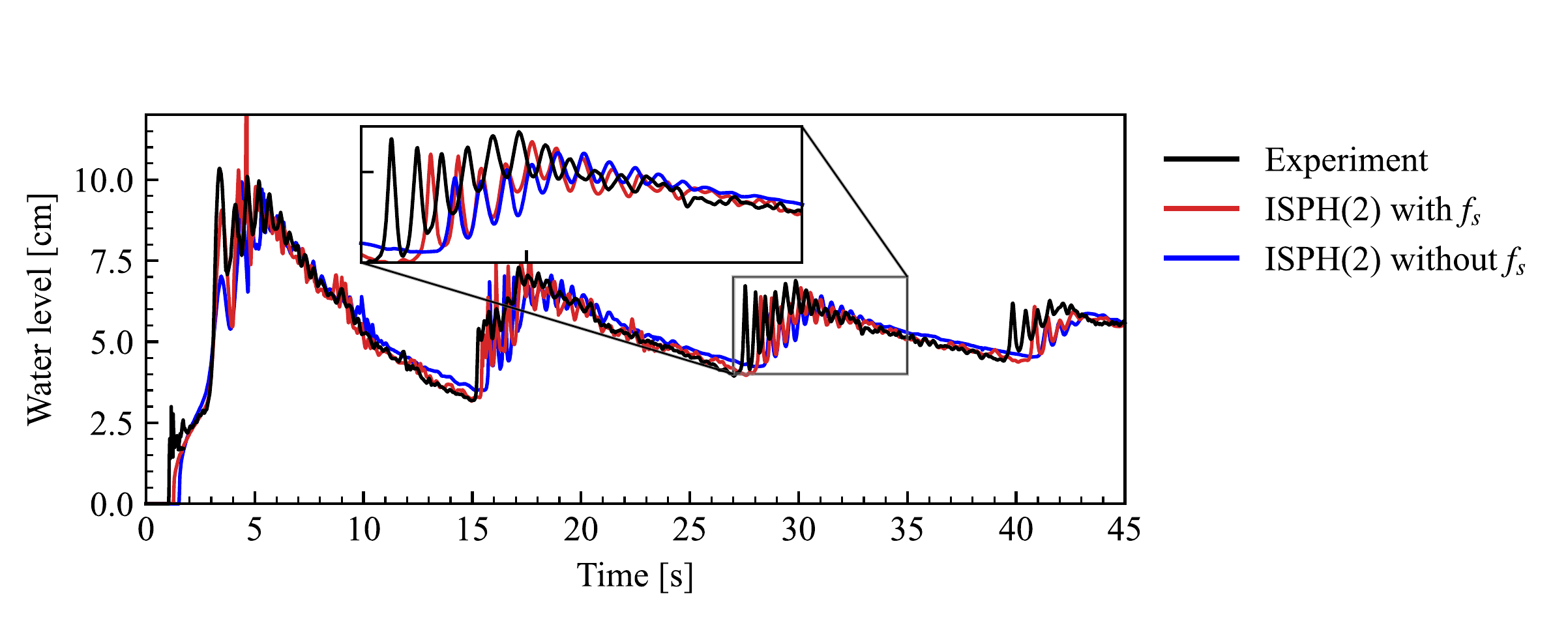}
    \caption{Gauge 3 ($x=393.5$~cm)}
  \end{subfigure}
  \caption{Time histories of the water levels obtained using the BFE-SPH ($\alpha=2.0$, $d_0=0.2$) with ISPH(2), both with and without the damping function $f_s$ at each gauge}\label{fig:s4_be_fs}
\end{figure}

The last discussion in this example is on the necessity of a wall function in the eddy viscosity models, such as the Smagorinsky model. 
\refF{fig:s4_be_fs} shows the time histories of the water levels on Gauges 1, 2, and 3 ($x=557.5$, 492.5, and 393.5~cm). In these results, the finest resolution model with ISPH(2) is utilized for cases with and without the wall function $f_s$ defined in \refC{sec3-4:Smagorinsky}. The no-slip boundary condition without $f_s$ may cause too much friction on the solid boundary and induce slower water movement than the experimental test. As discussed in the former CFD communities\Cite{hughes2001large,weickert2010investigation}, the wall function seems to decrease this tendency.

%%-------------------------------------------------------------------
%%  4. 4. \texorpdfstring{$\sigma$}{sigma}--SPH method (Hydrostatic pressure problem)
%%-------------------------------------------------------------------
\subsection{Hydrostatic pressure and dynamic problem using the \texorpdfstring{$\sigma$}{sigma}--SPH method}
Hydrostatic pressure and dynamic problem simulations are performed using the $\sigma$-SPH method with ISPH(2) in this section. For the hydrostatic pressure problem, the same analytical model is used to perform analyses with both the bottom boundary-fitted particle method (BF-SPH) and the bottom boundary-fitted ellipsoidal particle method (BFE-SPH), and a comparison of computational efficiency is carried out. For the dynamic problem, the influence of the volume conservation techniques described in \refC{sec3-6:volume_conservation} on temporal volume changes is examined.
\subsubsection{Hydrostatic pressure problem with the \texorpdfstring{$\sigma$}{sigma}--SPH method}
\begin{figure}[H]
  \centering
  \includegraphics[width=0.861\linewidth]{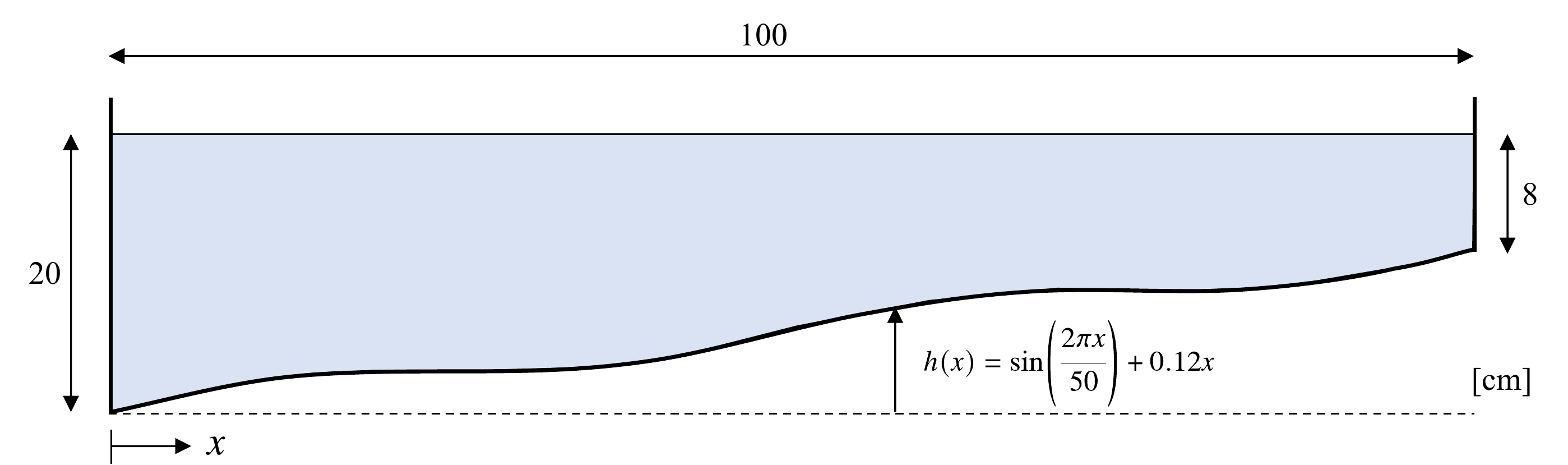}
  \caption{Hydrostatic pressure problem setup for comparing computational cost}\label{fig:s4_s_model}
\end{figure}

A hydrostatic pressure problem in a tank with a bottom-shaped sine and slope as shown in \refF{fig:s4_s_model} is calculated using the BF-SPH, BFE-SPH, and $\sigma$-SPH with ISPH(2). All particle models are standardized with the same vertical resolutions at the shallowest part ($x = 100$~cm). In BFE-SPH, $\alpha$ is set to 2.5, while in $\sigma$-SPH, \EDITs{$\alpha$ is configured to vary from 1.0 to 2.5.} \EDITf{In the projected space for $\sigma$-SPH, the computational model uses a fluid domain with dimensions $100\times20$ in the horizontal~($x$) and vertical~($z$) directions. The parameters $H$ and $\wh H$ for the $\sigma$-SPH method are both set to 20.0, respectively.}  The computational conditions are time step width~$\D t=5.0\times10^{-4}$~s, and initial particle spacing~$d_0=0.2$~cm in BF-SPH and $d_0=0.5$~cm in the other methods.\par
\begin{figure}[H]
  \centering
  \begin{subfigure}[b]{\linewidth}
      \centering
      \includegraphics[width=0.951\linewidth]{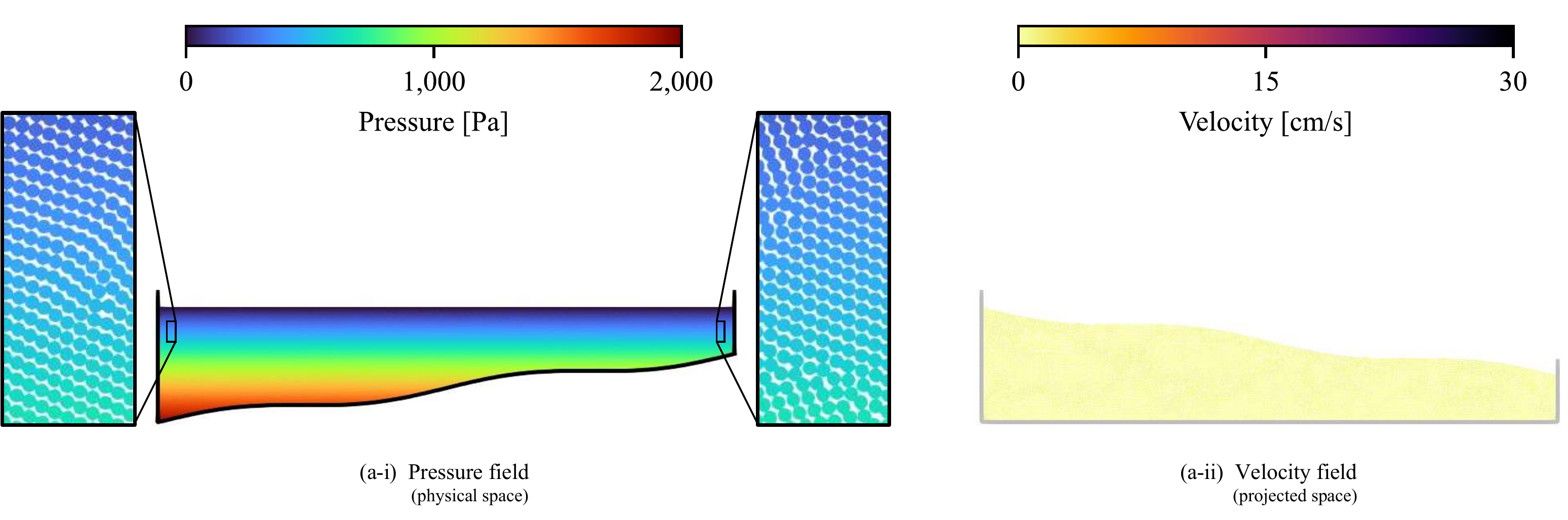}
      \caption{BF-SPH}
  \end{subfigure}
  \begin{subfigure}[b]{\linewidth}
      \centering
      \includegraphics[width=0.951\linewidth]{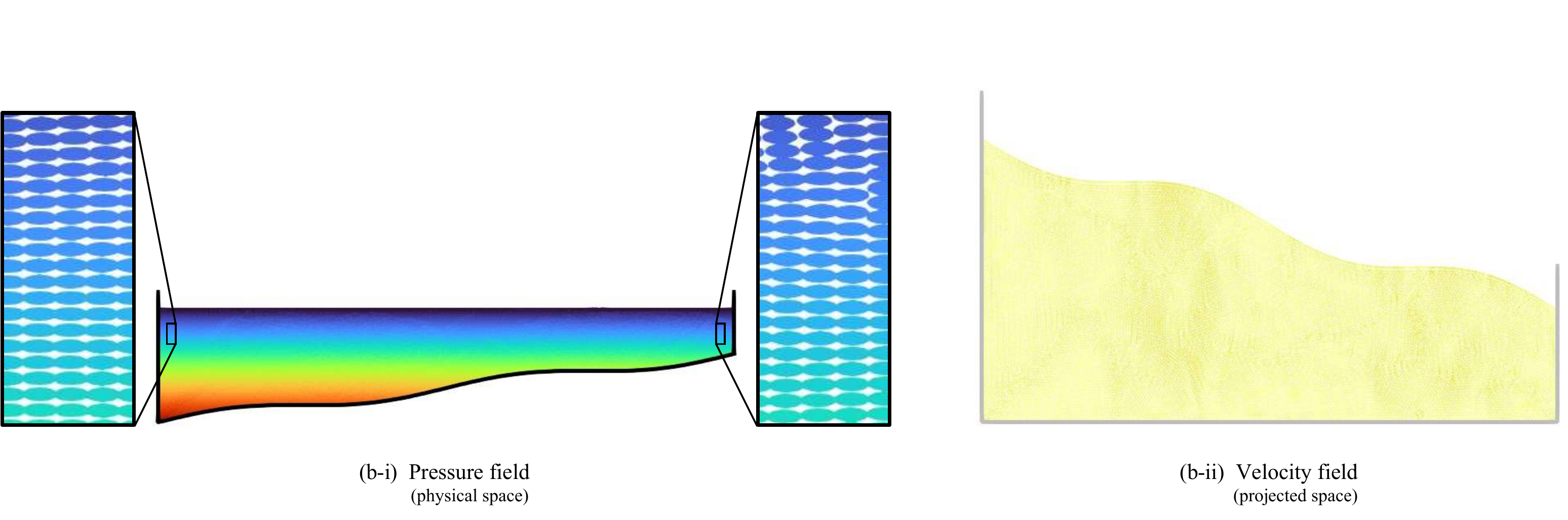}
      \caption{BFE-SPH}
  \end{subfigure}
  \begin{subfigure}[b]{\linewidth}
    \centering
    \includegraphics[width=0.951\linewidth]{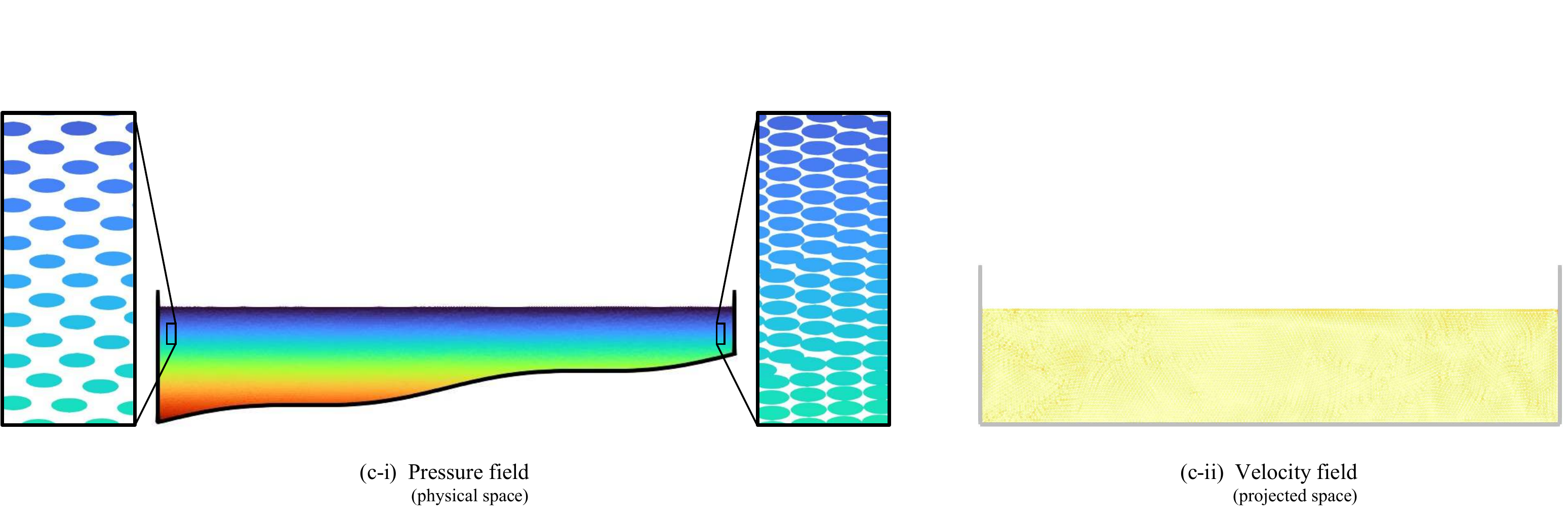}
    \caption{\EDIT{$\sigma$-SPH}}
  \end{subfigure}
  \caption{Pressure and velocity fields at 10~s calculated with the BF-SPH, BFE-SPH, and $\sigma$-SPH for the hydrostatic pressure problem}\label{fig:s4_s_comp}
\end{figure}
\begin{figure}[H]
  \centering
  \begin{subfigure}[b]{0.325\linewidth}
      \centering
      \includegraphics[width=0.99\linewidth]{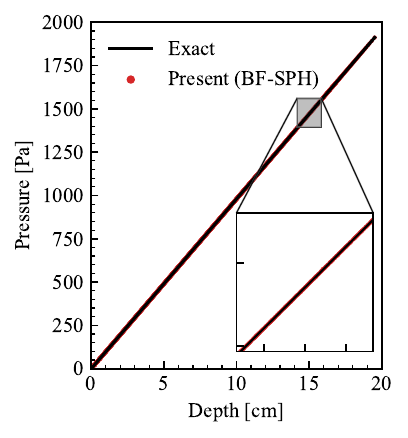}
      \caption{BF-SPH}
  \end{subfigure}
  \begin{subfigure}[b]{0.325\linewidth}
      \centering
      \includegraphics[width=0.99\linewidth]{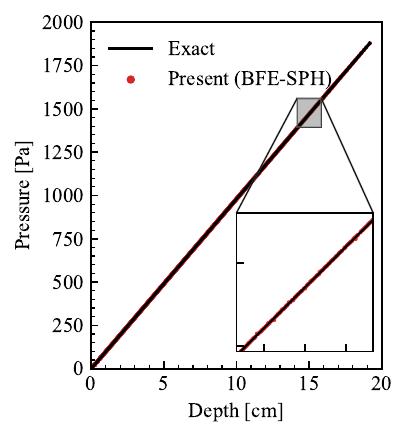}
      \caption{BFE-SPH}
  \end{subfigure}
  \begin{subfigure}[b]{0.325\linewidth}
    \centering
    \includegraphics[width=0.99\linewidth]{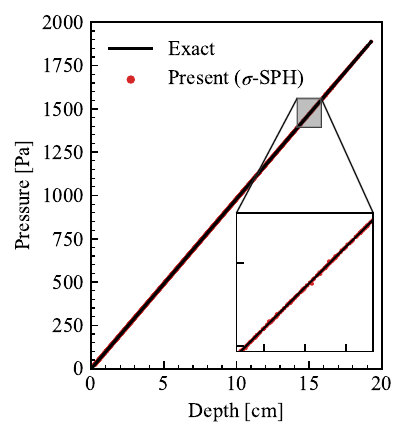}
    \caption{\EDIT{$\sigma$-SPH}}
  \end{subfigure}
  \caption{Pressure profiles at 10~s calculated using the BF-SPH, BFE-SPH, and $\sigma$-SPH for the hydrostatic pressure problem}\label{fig:s4_s_comp_p}
\end{figure}
\begin{figure}[H]
    \centering
    \begin{minipage}[b]{0.49\linewidth}
      \centering
      \includegraphics[height=0.557\linewidth]{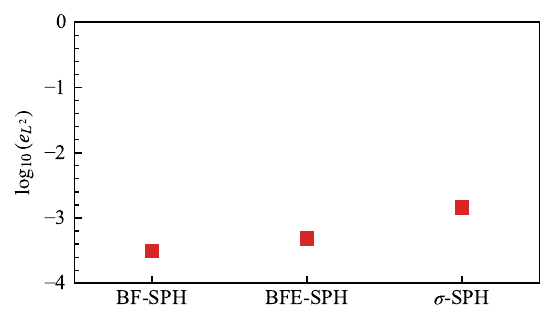}
      % \subcaption{Logarithmic relative $L^2$ errors in pressure at 10~s}\label{fig:s4_s_l2}
    \end{minipage}
    % \begin{minipage}[b]{0.49\linewidth}
    %   \centering
    %   \includegraphics[height=0.557\linewidth]{Figure_19_b.pdf}
    %   \subcaption{Number of particles and computational time}\label{fig:s4_s_time}
    % \end{minipage}
    % \caption{Logarithmic relative $L^2$ errors in pressure, number of particles, and computational time obtained using the BF-SPH, BFE-SPH, and $\sigma$-SPH for hydrostatic pressure problem}\label{fig:s4_s_l2_time}
    \caption{\EDITf{Logarithmic relative $L^2$ errors in pressure obtained using the BF-SPH, BFE-SPH, and $\sigma$-SPH for hydrostatic pressure problem}}\label{fig:s4_s_l2}
\end{figure}

\refF{fig:s4_s_comp} shows the pressure and velocity fields calculated using each method. The vertical resolution remains constant at $x=100$~cm. In contrast, for $\sigma$-SPH, the resolution varies with depth, demonstrating that it becomes coarser near $x=0$~cm compared to other methods. The pressure field is smoothly distributed for all methods. Furthermore, the velocity field calculated with the $\sigma$-SPH exhibits almost no flow.
\refF{fig:s4_s_comp_p} presents the pressure distribution along the vertical axis, where the depth is defined as $z_{\mathrm{FS}} - z_i$, with $z_{\mathrm{FS}}$ representing the average $z$-coordinate of the free surface particles. The results for each method agree well with the exact solution. 
\EDITf{\refF{fig:s4_s_l2} shows the logarithmic relative $L^2$ errors in pressure obtained using each method. For BF-SPH, BFE-SPH, and $\sigma$-SPH, the errors decrease with increasing resolution, and overall, they are comparable across all methods. \refT{fig:s4_s_time} summarizes the number of particles, average memory usage per step, and computational time for each method. The number of particles in BFE-SPH and $\sigma$-SPH is approximately 2/5 and 1/4 of that in BF-SPH, respectively. The memory usage per step also decreases proportionally with the number of particles. Furthermore, the computation times for BFE-SPH and $\sigma$-SPH are approximately 1/3 and 1/6 of that in BF-SPH, respectively.}

\begin{table}[H]
  \centering
  \caption{\EDITf{Number of particles, average memory usage per step, and computational time obtained using the BF-SPH, BFE-SPH, and $\sigma$-SPH for hydrostatic pressure problem}}\label{fig:s4_s_time}
  \begin{tabular}{cccc}
      \toprule
      \textbf{VCT's Types} & \textbf{\# of particles}&\textbf{Avg. Memory [MB/step]}& \textbf{Comput. time [h]}\\ \midrule
      \textbf{BF-SPH}        & 37,285                      & 447.5                        & 28.34\\
      \textbf{BFE-SPH}       & 15,436                      & 179.2                         & 9.72\\
      \textbf{$\sigma$-SPH}  & 9,086                       & 109.6                         & 4.41\\ \bottomrule
  \end{tabular}
\end{table}

\subsubsection{\EDIT{Dynamic problem with the \texorpdfstring{$\sigma$}{sigma}--SPH method}}
\begin{figure}[H]
  \centering
  \includegraphics[width=0.553\linewidth]{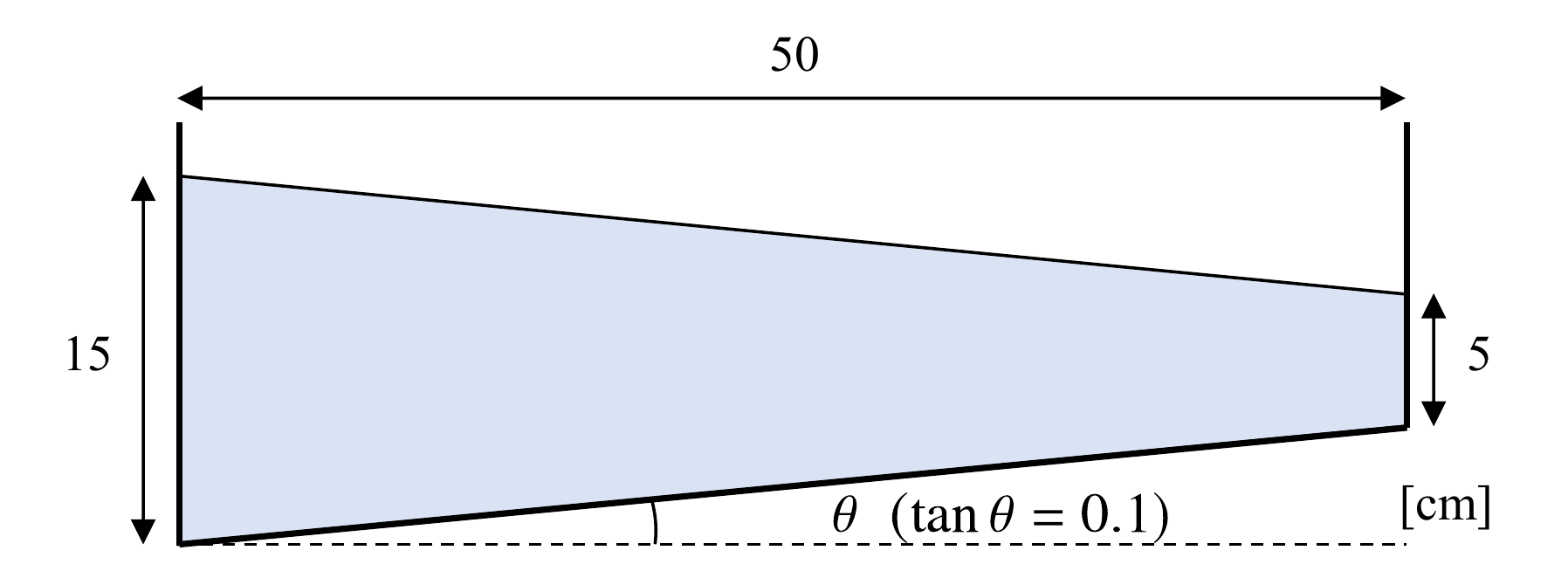}
  \caption{Dynamic problem setup for verifying volume conservation}\label{fig:s4_s_model2}
\end{figure}

A dynamic problem in a tank with a slope, as shown in \refF{fig:s4_s_model2}, is calculated using the $\sigma$-SPH with ISPH(2). 
\EDITs{The calculations are performed under the following three conditions:
% \begin{enumerate}[label=\textsf{(\Alph*)}]
\begin{enumerate}[label=(\Alph*)]
  \item \label{item:a}without volume conservation technique (solving \refE{eq:VCT_ppe});
  \item \label{item:b}with volume conservation technique but without volume smoothing (solving  \refE{eq:VCT_ppe_w0});
  \item \label{item:c}with both volume conservation technique and volume smoothing (solving \refE{eq:VCT_ppe_w}).
\end{enumerate}
These conditions are described in \refC{sec3-6:volume_conservation}.} The computational conditions are time step width $\D t=1.0\times 10^{-4}$~s, and initial particle spacing $d_0 = 0.25$~cm. In this section, the validation is carried out for the vertical scale factor $\alpha$ in the $0.\dot 6$ to 2.0. \EDITf{In the projected space, the computational model uses a fluid domain with dimensions $50\times10$ in the horizontal~($x$) and vertical~($z$) directions. The parameters $H=7.5$ and $\wh H=5.0$ are adopted for the $\sigma$-SPH method.} The coefficient $C_V$ for the volume smoothing varies from $\Delta t \times 10^{-3}$ to $\Delta t \times 10^5$~[-], and its appropriate value is evaluated. The volume conservation is checked using errors $e_{\rho}$ in the density and errors $e_{V}$ between the initial and current total volume.
\begin{equation}
    e_{\rho}=\cfrac{\aab{\rho}^N-\rho}{\rho}\times 100,
\end{equation}
\begin{equation}
    \aab{\rho}^N_i=\cfrac{J_i^N}{\bar{J}^N_i}\aab{\,\wh\rho\,}^N_i;\quad\aab{\,\wh\rho\,}^N_i:=\sphsum{V_j \rho_j w_{ij}},
\end{equation}
\begin{equation}
    e_{V}=\cfrac{V^N-V^0}{V^0}\times 100,
\end{equation}
\begin{equation}
    V^N=\sum_{i}^{N_{\mathrm{SPH}}}\aab{V}_i;\quad \aab{V}_i:=\cfrac{m_i}{J_i\aab{\,\wh\rho\,}^N_i}.
\end{equation}\par
\begin{figure}[H]
  \centering
  \begin{minipage}[b]{\linewidth}
    \centering
    \includegraphics[width=0.938\linewidth]{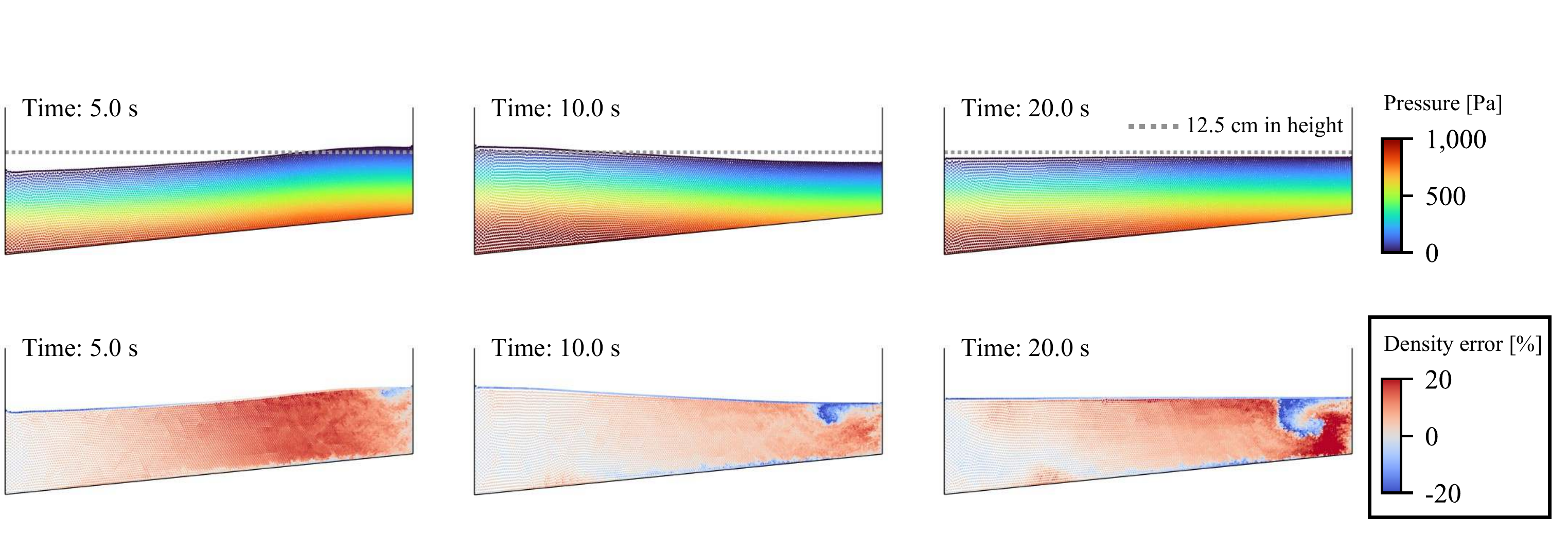}
    \subcaption{Condition A (solving \refE{eq:VCT_ppe})}\label{fig:s4_s_vol_comp_wo}
  \end{minipage}
  \begin{minipage}[b]{\linewidth}
    \centering
    \includegraphics[width=0.938\linewidth]{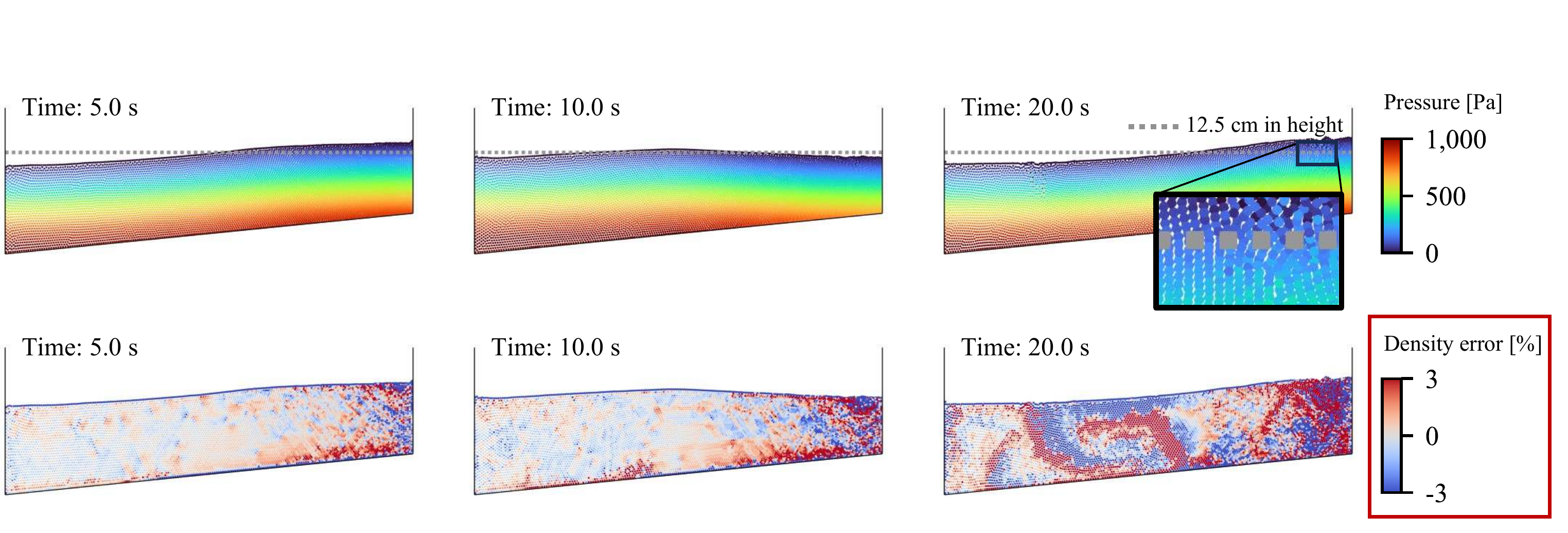}
    \subcaption{Condition B (solving \refE{eq:VCT_ppe_w0})}\label{fig:s4_s_vol_comp_w0}
  \end{minipage}
  \begin{minipage}[b]{\linewidth}
    \centering
    \includegraphics[width=0.938\linewidth]{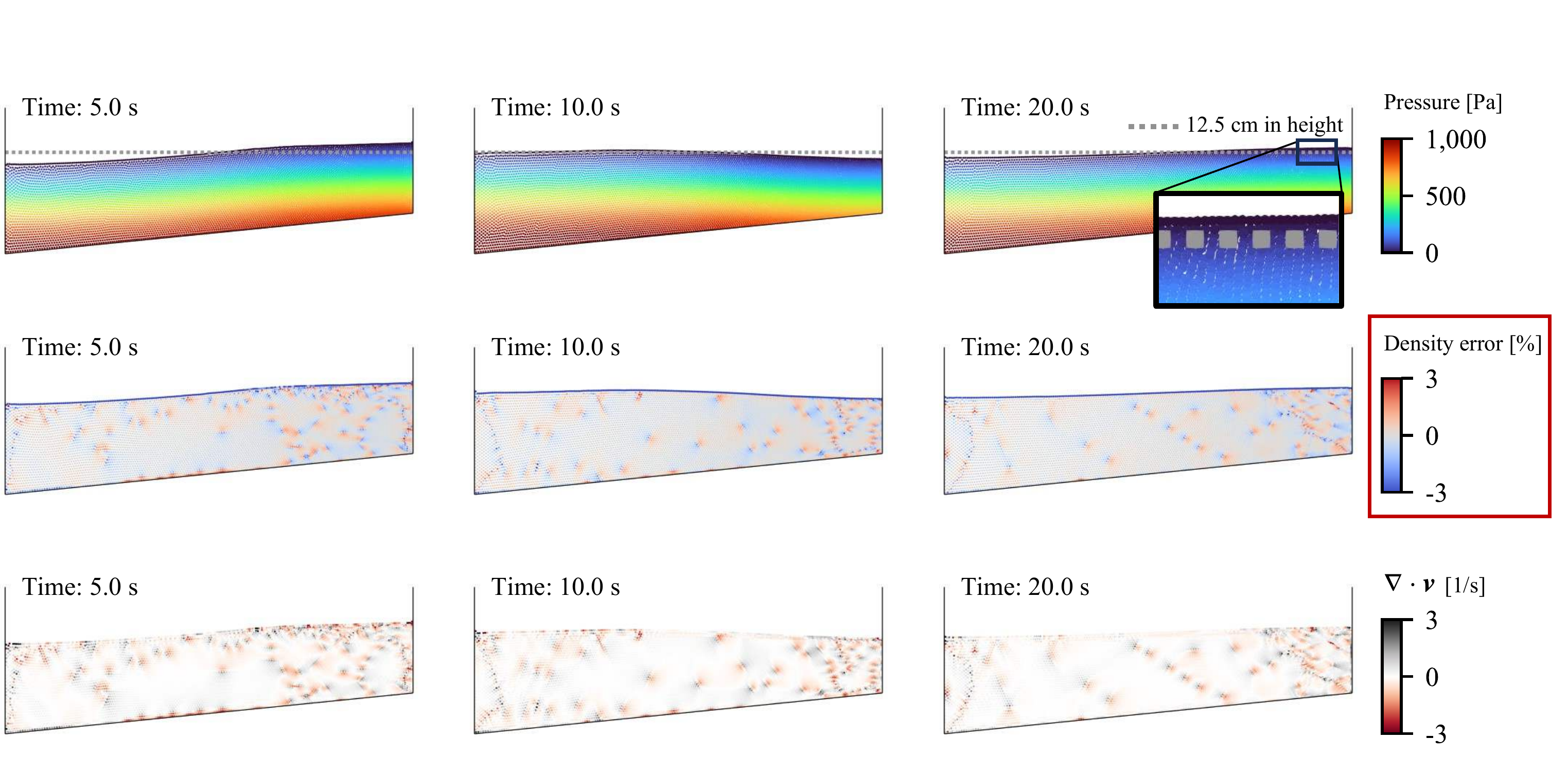}
    \subcaption{Condition C (solving \refE{eq:VCT_ppe_w} with $C_v = \D t \times 10^2$)}\label{fig:s4_s_vol_comp_w}
  \end{minipage}
  \caption{\EDITf{Pressure, density error, and velocity divergence fields obtained using the $\sigma$-SPH under Conditions A, B and C ($C_V=\D t\times 10^{2}$) in a dynamic problem}}\label{fig:s4_s_vol_comp}
\end{figure}

\EDITf{\refF{fig:s4_s_vol_comp} shows the pressure, density error $e_\rho$, and velocity divergence fields obtained using the $\sigma$-SPH under Conditions~A--C, where Condition~C uses $C_v = \D t \times 10^2$. 
Under Condition~A, the pressure distribution remains smooth.  However, the density error increases when particles move from the left to the right side at $t = 5.0$~s, decreases upon their return to the left side at $t = 10.0$~s, and eventually shows a mixture of positive and negative values at $t = 20.0$~s. 
Furthermore, the free surface position at 20.0~s falls below the reference height 12.5~cm (indicated by the gray dashed line), suggesting a reduction in fluid volume.
\EDITT{This result is considered to be since $J^0$ becomes discontinuous with particle motion, causing the stabilization term to deviate from the appropriate range.}
For Condition~B, both the density error and pressure fluctuations remain small up to approximately 5~s. However, the density error increases after 10.0~s, and pressure fluctuations are observed at 20.0~s. 
Under Condition~C, the density error remains within approximately $\pm 3$~\% up to 20~s, and the pressure field is also smooth. The velocity divergence field is generally close to zero, although some regions exhibit relatively large values. These high-divergence regions spatially correspond to areas exhibiting large density errors, which is characteristic of the stabilized ISPH method. This method inherently reduces density deviations by applying pressure correction.}\par
\begin{figure}[H]
  \centering
  \begin{minipage}[b]{\linewidth}
    % \centering
    % \includegraphics[width=14.34cm]{Figure_22_a.pdf}
    \includegraphics[width=0.9975\linewidth]{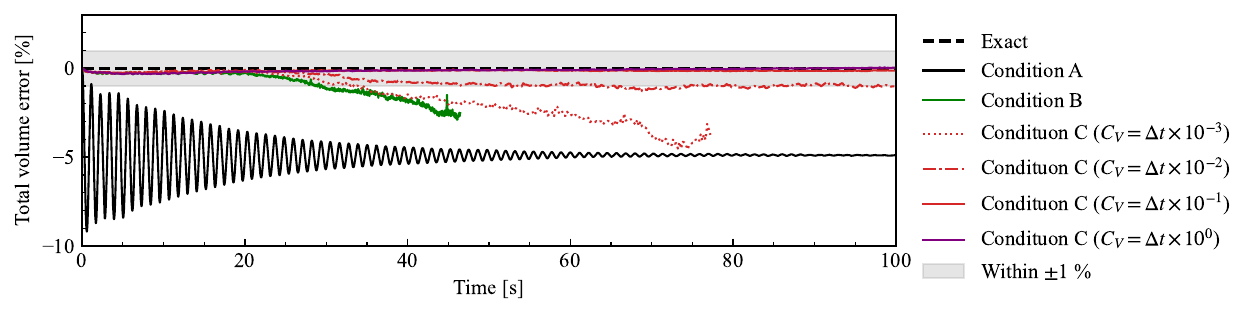}
    \subcaption{Conditions~A--C, where Condition~C uses $C_V$ ranging from $\D t \times 10^{-3}$ to $\D t \times 10^0$
    }\label{fig:s4_s_vol_error1}
  \end{minipage}
  \begin{minipage}[b]{\linewidth}
    % \centering
    % \includegraphics[width=14.34cm]{Figure_22_b.pdf}
    \includegraphics[width=0.9888\linewidth]{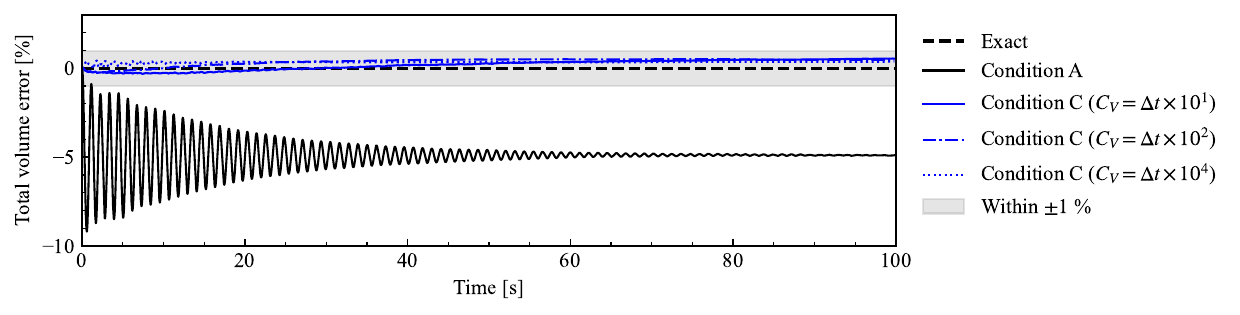}
    \subcaption{Conditions~A and C, where Condition~C uses $C_V$ ranging from $\D t \times 10^{1}$ to $\D t \times 10^4$}\label{fig:s4_s_vol_error2}
  \end{minipage}
  \caption{\EDITs{Time histories of total volume errors $e_V$ obtained using the $\sigma$-SPH under Conditions A--C in a dynamic problem}}\label{fig:s4_s_vol_error}
\end{figure}

% The pressure field is smoothly distributed with or without volume conservation techniques. In addition, the pressure fields at 100~s show that, without the volume conservation techniques, the free surface position is lower than 12.5 cm in height (gray dashed line), indicating a decrease in volume. In contrast, with the volume conservation techniques, the free surface position remains almost identical, demonstrating volume conservation. Without the volume conservation technique, the density error is positive when a particle moves from the left side to the right side ($t = 0.5$ s), smaller errors when it returns to the left side ($t = 1.0$ s), and a mixture of positive and negative errors after a sufficient amount of time ($t = 100$ s) as shown in \hyperref[fig:s4_s_vol_comp_wo]{Fig.~\ref*{fig:s4_s_vol_comp}(a)}. In contrast, the density errors remain within $\pm 20$~\% with the volume preservation techniques as shown in \hyperref[fig:s4_s_vol_comp_w]{Fig.~\ref*{fig:s4_s_vol_comp}(c)}.\par
\EDITs{\refF{fig:s4_s_vol_error} shows the time histories of total volume errors $e_V$ obtained using $\sigma$-SPH under Conditions A--C.
Under Condition~A (black line), the total volume exhibits oscillations until approximately 80~s, after which it stabilizes at around $-5$~\%.
For Condition~B (green line), the error remains within $\pm 1$~\% up to roughly 20~s. However, the volume gradually decreases afterward, and the simulation result diverges at around 50~s. 
\EDITT{This result is thought to occur for the same reason noted above, as the stabilization term deviates from the appropriate range due to the discontinuity of $J^0$.}
Under Condition~C, for the cases with $C_V\le \D t\times 10^0$ (red and purple lines in \hyperref[fig:s4_s_vol_error1]{Fig.~\ref*{fig:s4_s_vol_error}(a)}), the total volume initially decreases. Afterward, when $C_V=\D t\times 10^0,\, \D t\times 10^{-1}$, the error subsequently increases, but stays within $\pm 1$~\% up to 100~s. Meanwhile, when $C_V=\D t\times 10^{-2},\, \D t\times 10^{-3}$, the error continues to decrease over time. 
\EDITT{This phenomenon is thought to result from the inadequacy of the volume smoothing relative to the particle displacement, leading to a condition similar to Condition~B.} 
Additionally, for the cases with $C_V > \D t \times 10^0$ (blue lines in \hyperref[fig:s4_s_vol_error2]{Fig.~\ref*{fig:s4_s_vol_error}(b)}), the error remains within $\pm 1$~\% throughout the simulation. However, in the case of $C_V = \D t\times 10^5$, the simulation diverges immediately after it begins.}\par
\begin{figure}[H]
  \centering
  \includegraphics[width=0.8921\linewidth]{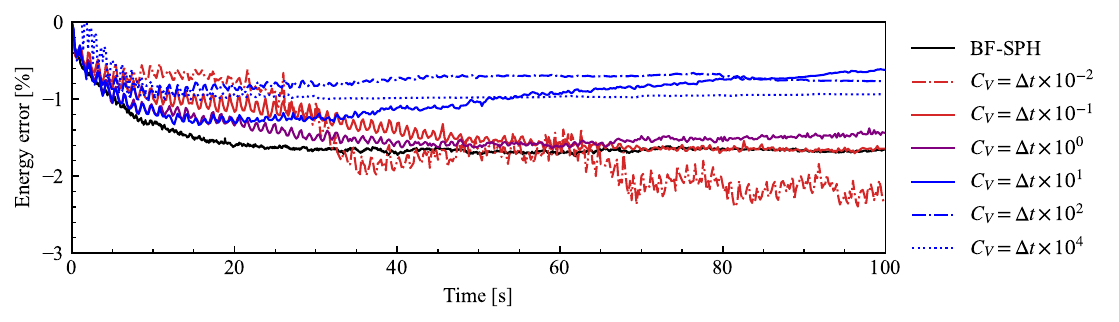}
  \caption{\EDITf{Time histories of energy errors obtained using the $\sigma$-SPH under Condition~C in a dynamic problem, where the energy errors are measured relative to the initial energy}}\label{fig:s4_s_ene_error}
\end{figure}

\EDITf{\refF{fig:s4_s_ene_error} shows the time histories of the energy errors obtained using the $\sigma$-SPH under Condition~C in which the total volume error remains within $\pm 1$~\%. The energy error is evaluated relative to the initial energy. For comparison, results obtained using the BF-SPH method are also included as reference. At $t = 100$~s, the energy error remains within $\pm 1$~\% relative to the BF-SPH results for all tested values of $C_V$. Among these, the case with $C_V = \Delta t \times 10^{-1}$ exhibits the best agreement with the BF-SPH reference results.}\par
\begin{figure}[H]
  \centering
  \begin{minipage}[b]{\linewidth}
    % \centering
    % \includegraphics[width=14.34cm]{Figure_22_a.pdf}
    \includegraphics[width=0.9504\linewidth]{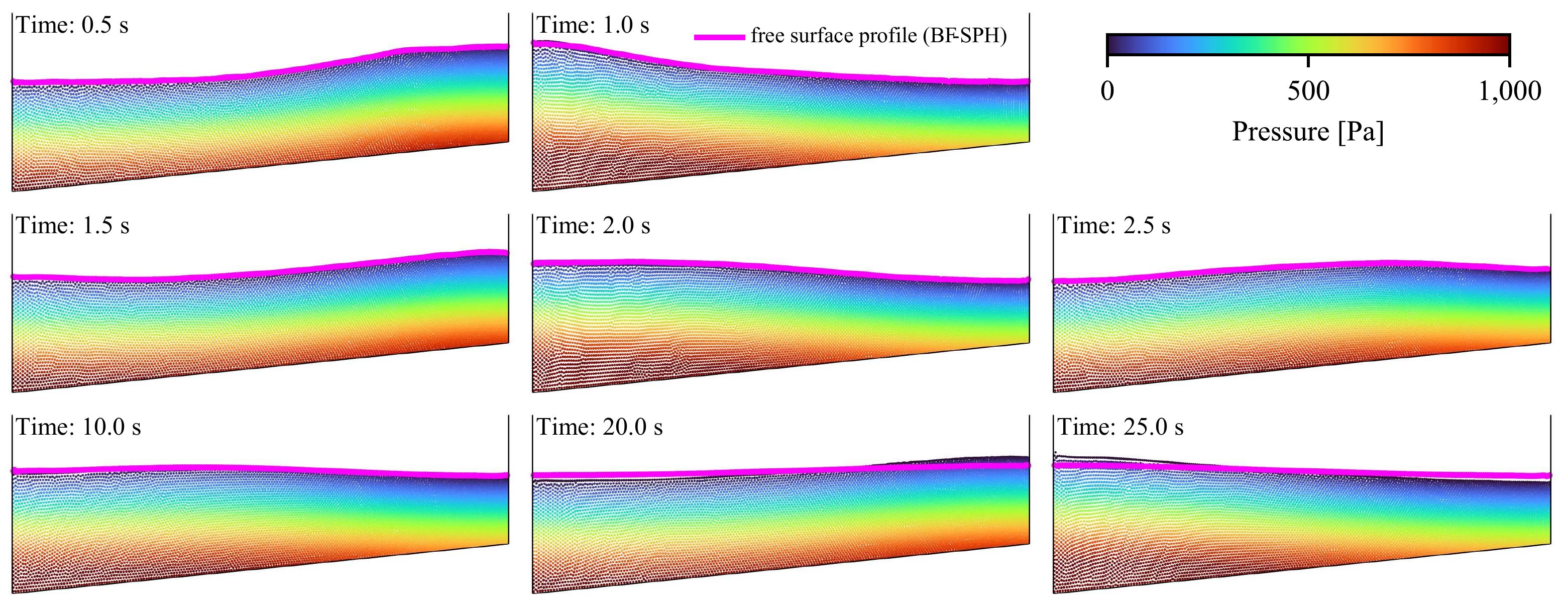}
    \subcaption{$C_V=\D t \times 10^{-1}$}\label{fig:s4_s_comp_bf_1}
  \end{minipage}
  \begin{minipage}[b]{\linewidth}
    % \centering
    % \includegraphics[width=14.34cm]{Figure_22_b.pdf}
    \includegraphics[width=0.9504\linewidth]{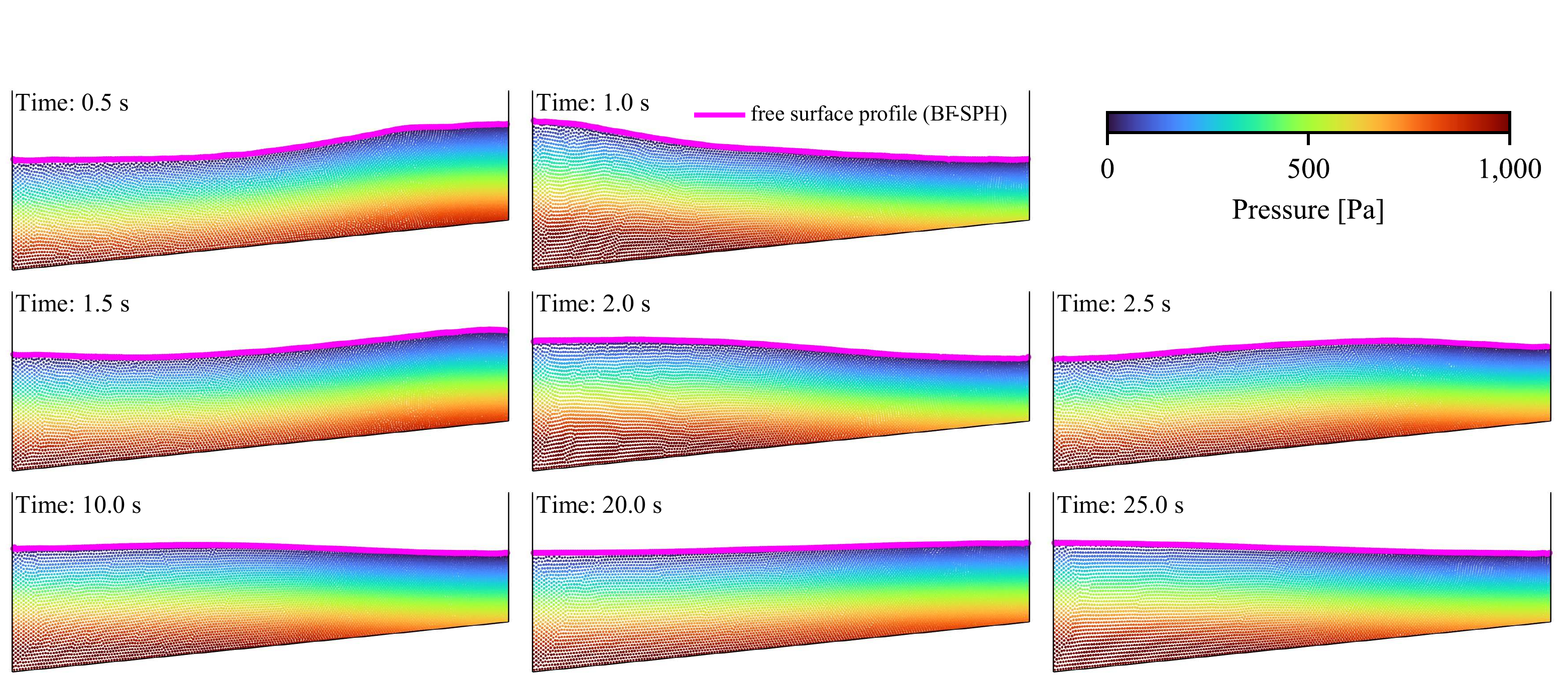}
    \subcaption{$C_V=\D t \times 10^{2}$}\label{fig:s4_s_comp_bf_2}
  \end{minipage}
  \begin{minipage}[b]{\linewidth}
    % \centering
    % \includegraphics[width=14.34cm]{Figure_22_b.pdf}
    \includegraphics[width=0.9504\linewidth]{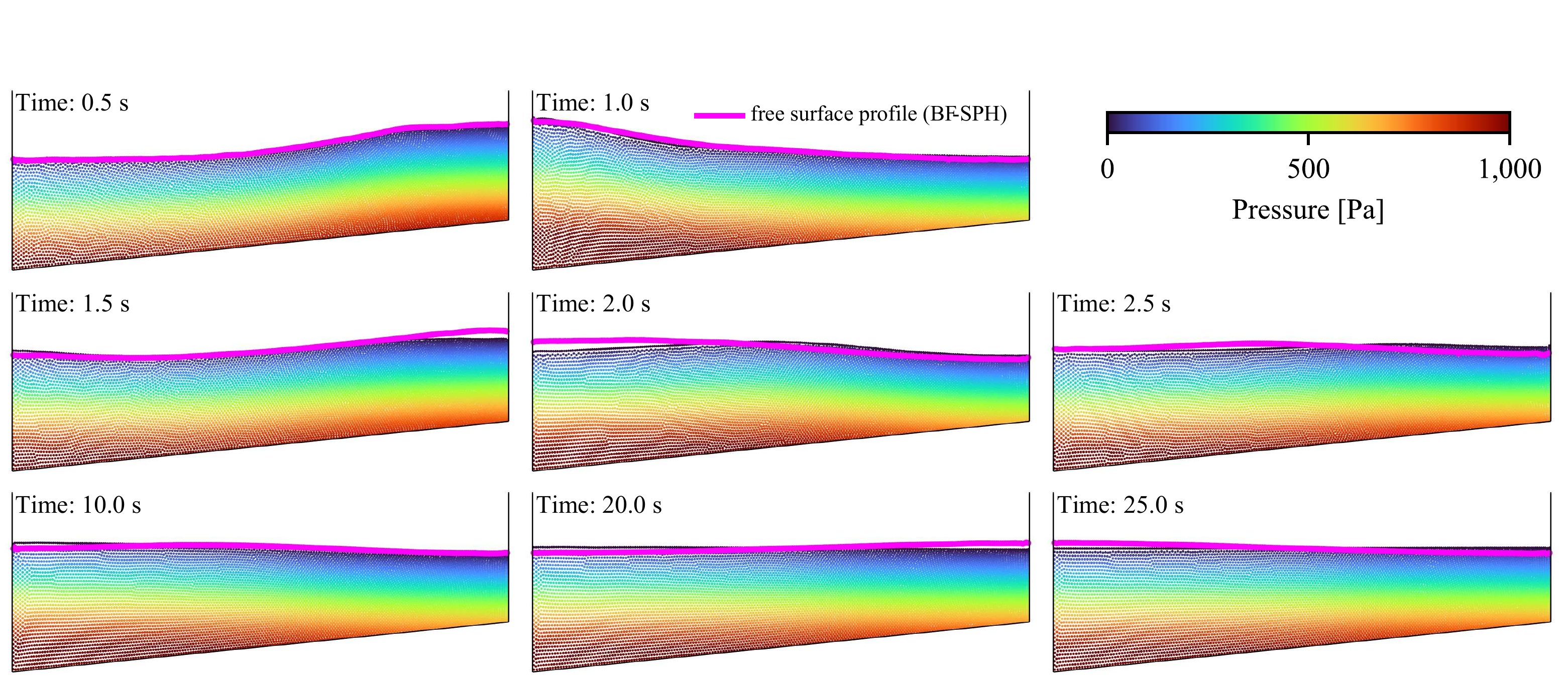}
    \subcaption{$C_V=\D t \times 10^{4}$}\label{fig:s4_s_comp_bf_3}
  \end{minipage}
  \caption{\EDITf{Pressure fields obtained using the $\sigma$-SPH under Condition~C and comparison of the free surface profile with that from BF-SPH in a dynamic problem}}\label{fig:s4_s_comp_bf}
\end{figure}

\EDITf{\refF{fig:s4_s_comp_bf} shows the pressure fields obtained using the $\sigma$-SPH under Condition~C, along with the comparisons of the free surface profile against the BF-SPH results (shown in pink). In the case of $C_V = \Delta t \times 10^{-1}$, the free surface motion is consistent with the BF-SPH result up to 10.0~s. However, deviations appear at 20.0~s and become more pronounced at 25.0~s. For $C_V = \Delta t \times 10^{2}$, the free surface position remains in good agreement with the BF-SPH result throughout the entire simulation. In contrast, when $C_V = \D t \times 10^{4}$, the wave attenuation is observed, and the free surface begins to deviate from the BF-SPH profile as early as 1.5~s. After that, the motion no longer matches the BF-SPH reference solution. 
While the case with $C_V = \D t \times 10^{-1}$ demonstrates good performance in terms of volume conservation and minimal energy loss, 
the accurate reproduction of the wave motion is considered more critical in this study.  Therefore, the value $C_V = \D t \times 10^2$ is adopted for the subsequent simulations.}
%In terms of volume conservation and energy loss, $C_V = \Delta t \times 10^{-1}$ performs well. However, since the reproducibility of the motion is considered more important, this study applies using $C_V = \Delta t \times 10^2$.}

\subsubsection{\EDIT{3-D dam break problem with the \texorpdfstring{$\sigma$}{sigma}--SPH method}}
\begin{figure}[H]
  \centering
  \includegraphics[width=0.9609\linewidth]{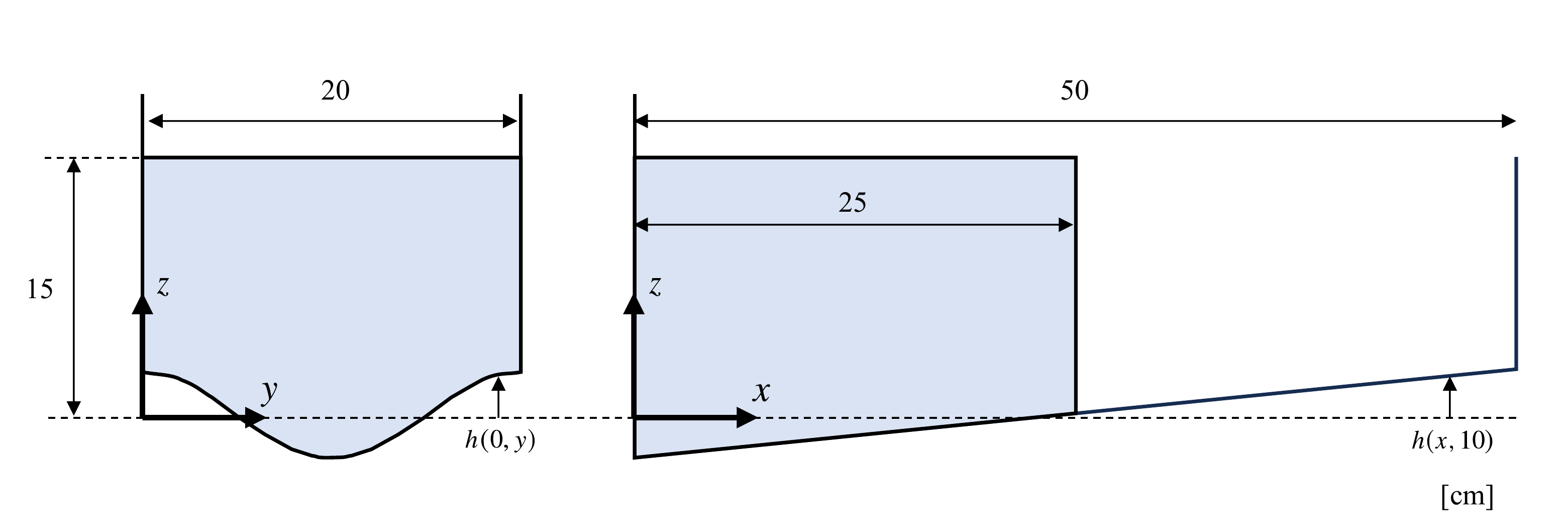}
  \caption{3-D dam break problem setup ($yz$-plane at $x=0$, $xz$-plane at $y=10$)}\label{fig:s4_s_model_dam}
\end{figure}

A three-dimensional dam break problem is simulated based on the computational model setup in \refF{fig:s4_s_model_dam}. The bottom elevation $h$ is defined as
\begin{align}
  h(x,y)=0.1x+2\cos\pab{\cfrac{2\pi x}{20}}.
\end{align}
The simulation conditions include a time step size of $\D t=1.0\times 10^{-4}$~s and initial particle spacing $d_0 = 0.25$~cm. 
The simulation is conducted under Condition~C using $C_V=\D t \times 10^2$. 
In the projected space, the computational domain measures $25\times20\times10$ in the horizontal~($x$), depth~($y$), and vertical~($z$) directions, respectively. The parameters $H=15.0$ and $\wh H=10.0$ are adopted for the $\sigma$-SPH method. \par
\begin{figure}[H]
  \centering
  \begin{minipage}[b]{\linewidth}
    \centering
    \includegraphics[width=\linewidth]{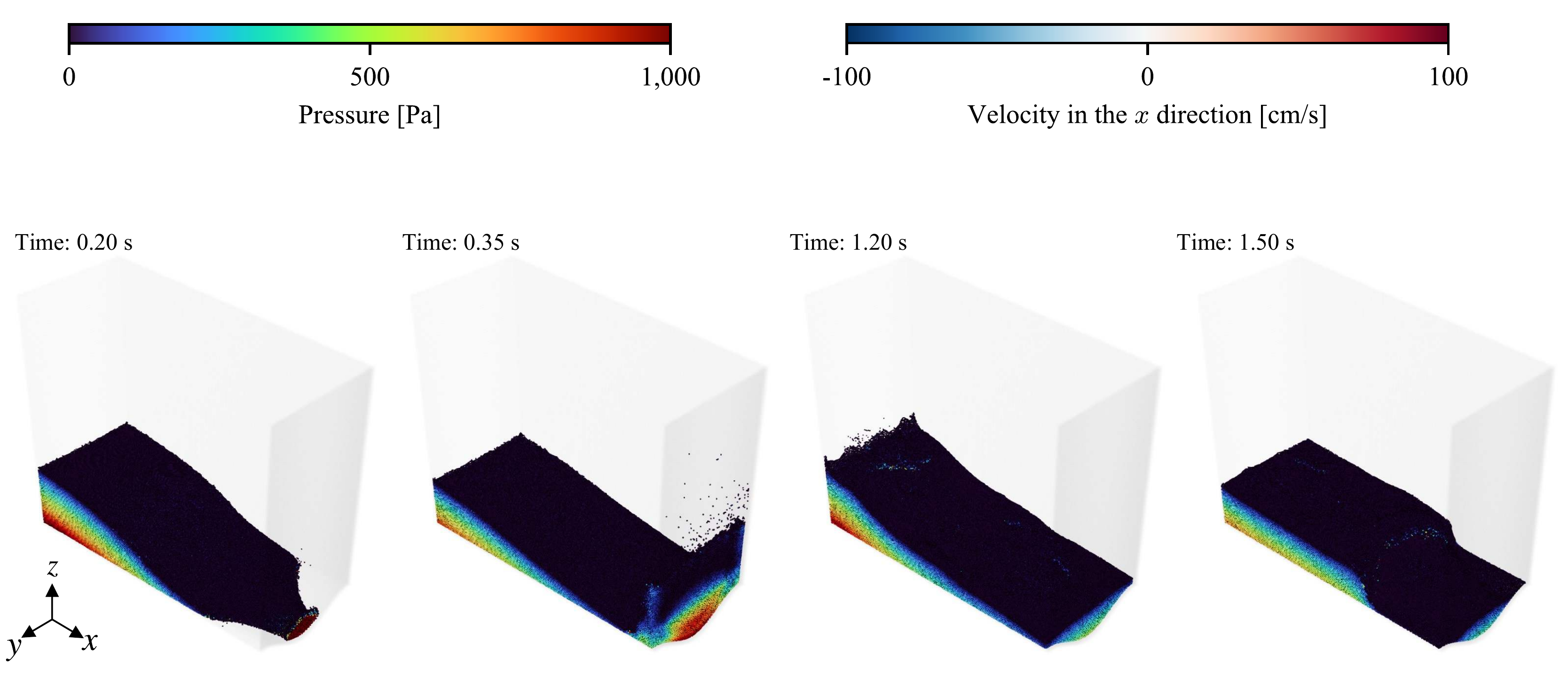}
    \subcaption{Physical space (pressure fields)}\label{fig:s4_s_3D_a}
  \end{minipage}
  \begin{minipage}[b]{\linewidth}
    \centering
    \includegraphics[width=\linewidth]{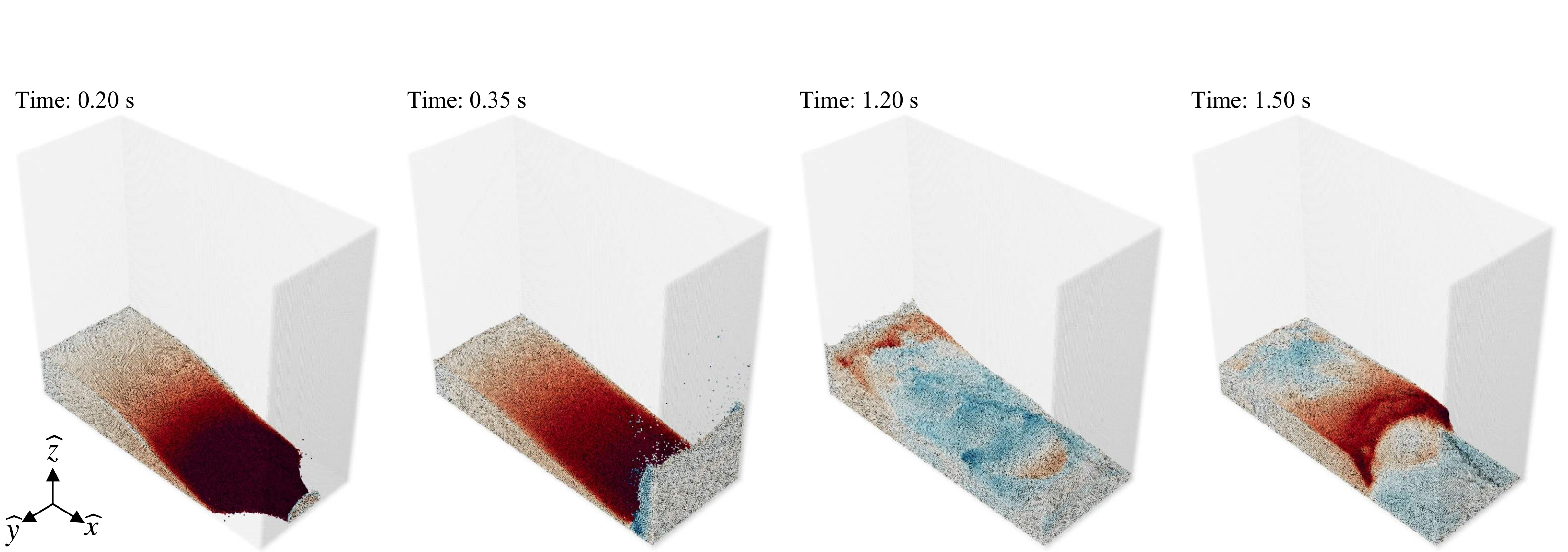}
    \subcaption{Projected space (velocity fields)}\label{fig:s4_s_3D_b}
  \end{minipage}
  \caption{Pressure and $x$ direction velocity fields obtained using the $\sigma$-SPH under Condition~C with $C_V=\D t\times 10^2$ in a 3-D dam break problem}\label{fig:s4_s_3D}
\end{figure}

\refF{fig:s4_s_3D} shows the pressure and $x$-direction velocity fields. At 0.20~s, the water flows into the region with a lower bottom elevation and collides with the wall.
Subsequently, the reflected wave propagates backward, impacting the opposite wall and initiating a run-up at 1.20~s. At this moment, the velocity in the $x$-direction is nearly uniform along the $y$-axis.
At 1.50~s, the surface velocity gradually increases in the shallower region.
These results successfully capture the key three-dimensional characteristics of fluid motion. Additionally, a splash is observed at 0.35~s, confirming that the simulation remains stable even under highly dynamic conditions.
The results also demonstrate robustness in handling inflow into initially dry regions.\par
\begin{figure}[H]
  \centering
  \begin{minipage}[b]{\linewidth}
    \centering
    \includegraphics[width=0.9827\linewidth]{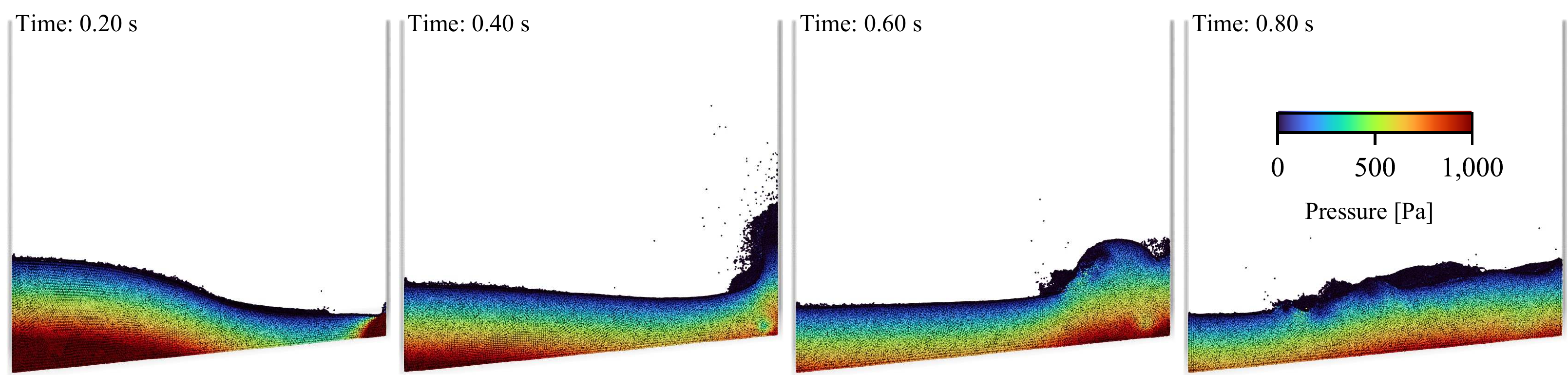}
    \subcaption{Physical space (pressure fields)}\label{fig:s4_s_3D_2D_a}
  \end{minipage}
  \begin{minipage}[b]{\linewidth}
    \centering
    \includegraphics[width=0.9827\linewidth]{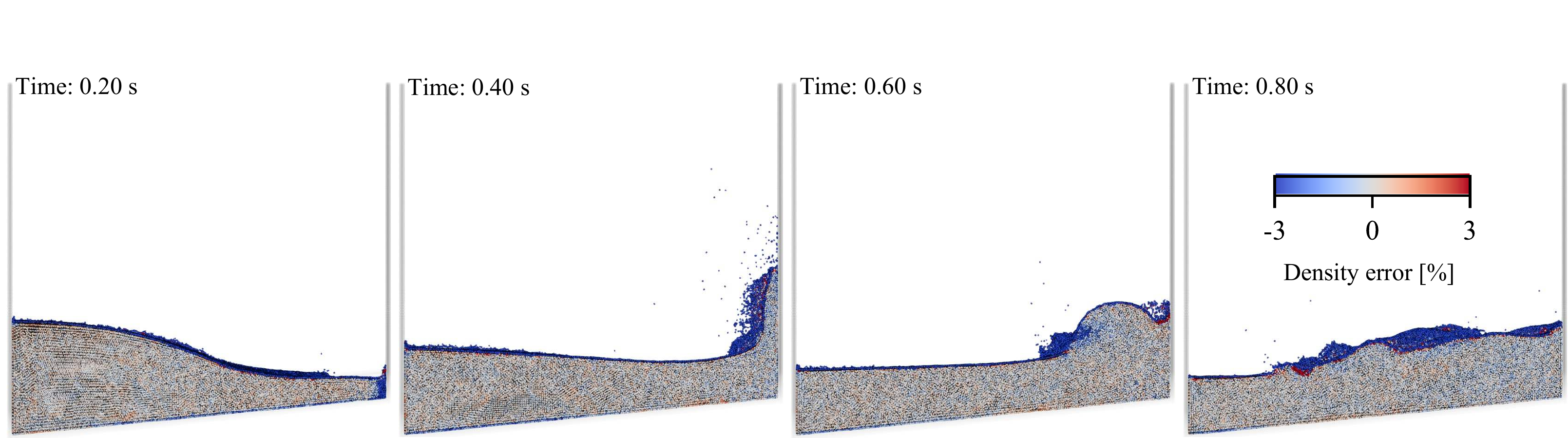}
    \subcaption{Physical space (density error fields)}\label{fig:s4_s_3D_2D_b}
  \end{minipage}
  \begin{minipage}[b]{\linewidth}
    \centering
    \includegraphics[width=0.9827\linewidth]{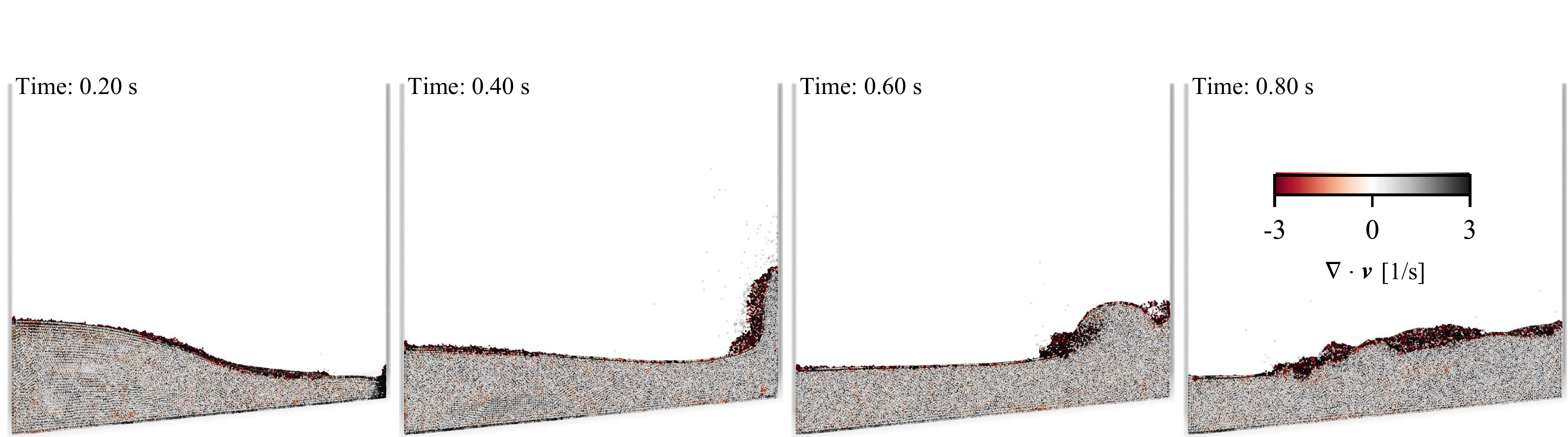}
    \subcaption{Physical space (velocity divergence fields)}\label{fig:s4_s_3D_2D_c}
  \end{minipage}
  \begin{minipage}[b]{\linewidth}
    \centering
    \includegraphics[width=0.9827\linewidth]{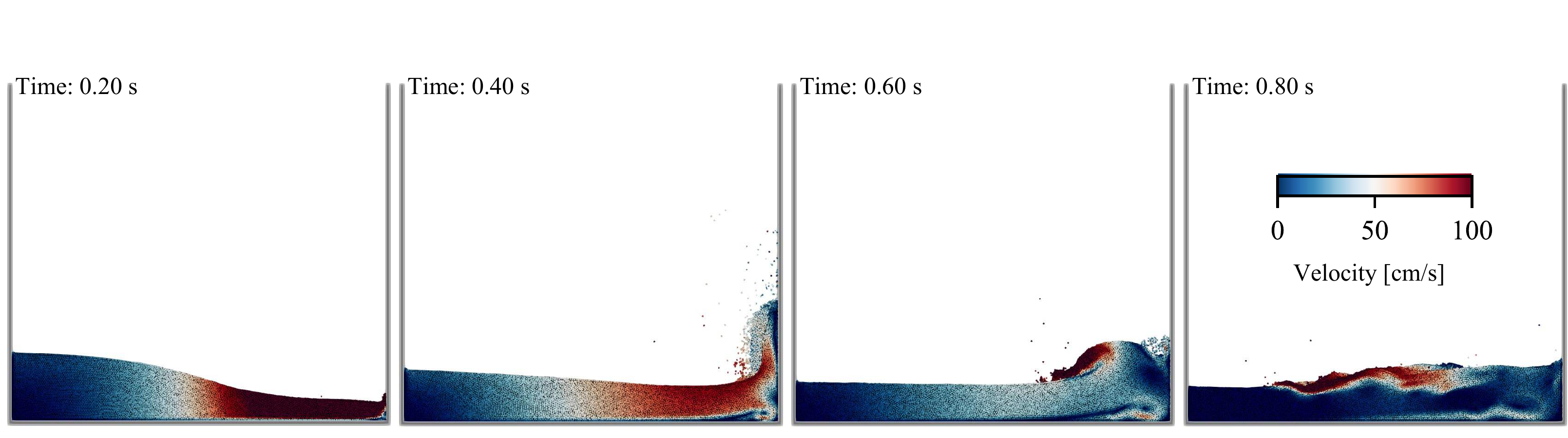}
    \subcaption{Projected space (velocity fields)}\label{fig:s4_s_3D_2D_d}
  \end{minipage}
  \caption{Pressure, density error, velocity divergence, and velocity fields obtained using the $\sigma$-SPH under Condition~C with $C_V=\D t\times 10^2$ in a 3-D dam break problem (cross-section at $y=10$~cm)}\label{fig:s4_s_3D_2D}
\end{figure}

\refF{fig:s4_s_3D_2D} shows the pressure, density error, velocity divergence, and velocity fields on the cross-sectional plane at $y=10$~cm. As shown in \hyperref[fig:s4_s_3D_2D_a]{Fig.~\ref*{fig:s4_s_3D_2D}(a)}, the pressure is smooth and hydrostatically consistent with the water depth.
A relatively low pressure region appears in the lower-right corner, which, as seen in \hyperref[fig:s4_s_3D_2D_d]{Fig.~\ref*{fig:s4_s_3D_2D}(d)}, is attributable to vortex formation.
This conforms to the physical plausibility of the pressure field.
As shown in \hyperref[fig:s4_s_3D_2D_b]{Fig.~\ref*{fig:s4_s_3D_2D}(b)} and \hyperref[fig:s4_s_3D_2D_c]{Fig.~\ref*{fig:s4_s_3D_2D}(c)}, both density error and velocity divergence remain close to zero, indicating good conservation properties.\par
\begin{figure}[H]
    % \centering
    % \includegraphics[width=14.34cm]{Figure_22_a.pdf}
    \includegraphics[width=0.9671\linewidth]{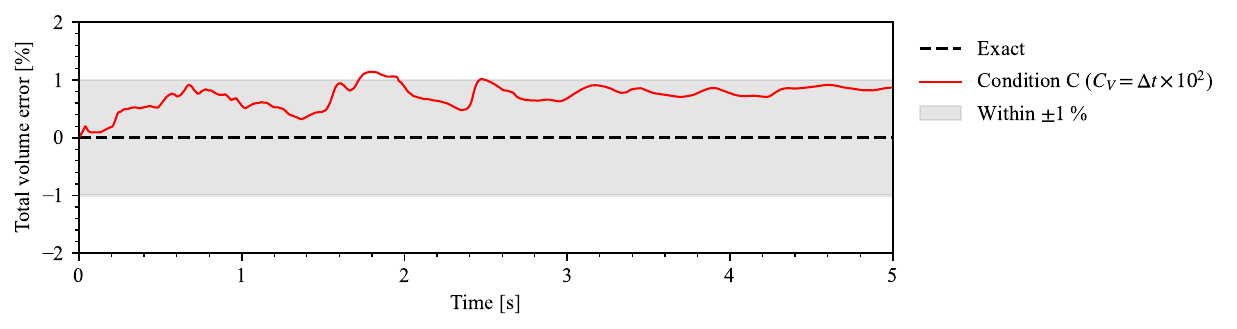}
  \caption{Time histories of total volume errors $e_V$ obtained using the $\sigma$-SPH under Conditions C with $C_V=\D t\times 10^2$ in a 3-D dam break problem}\label{fig:s4_s_vol_error_3D}
\end{figure}

\refF{fig:s4_s_vol_error_3D} shows the time history of the total volume error $e_V$.
Although a slight increase is observed from the initial state, the error remains within $\pm1$~\% throughout most of the simulation period.

%%%%%%%%%%%%%%%%%%%%%%%%%%%%%%%%%%%%%%%%%%%%%%%%%%%%%%%%%%%%%%%%%%%%%
%%  5. Conclusion
%%%%%%%%%%%%%%%%%%%%%%%%%%%%%%%%%%%%%%%%%%%%%%%%%%%%%%%%%%%%%%%%%%%%%
\section{Conclusion}

In this paper, three new particle methods; the bottom boundary-fitted particle method (BF-SPH), the bottom boundary-fitted ellipsoidal particle method (BFE-SPH), and the $\sigma$-SPH method using a $\sigma$-coordinate system are presented by generalizing the vertical coordinate transformation (VCT).\par

The particle method based on the Lagrange description requires the particle configuration to be uniformly distributed in each direction while the particles are moving at physical speed. For this reason, recent improvements have been made to the kernel function update method, using particle rearrangement methods to improve homogeneity and thus ensure accuracy or to update the kernel function in a way that does not cause a loss of accuracy when the particle arrangement is perturbed.
However, these improved methods only address the shortcomings of the particle method, and it is challenging to discretize the space with a bias to improve the computational efficiency associated with the physics problem. Therefore, by using a coordinate transformation, which has been proposed as an improvement method, the idea of the ellipsoidal particle method, which allows the use of kernel functions that are practically ellipsoid in real space, has been generalized to a new particle method that can increase the computational efficiency without reducing the computational accuracy.\par

The BF-SPH method shows that, in practice, the analysis of problems with complex wall geometries can be easily carried out by projecting the wall particles onto a wall with a flat surface instead of placing them uniformly inside the solid boundary. Furthermore, the BFE-SPH method combined with the ellipsoidal particle method allows efficient analysis by reducing the number of particles required for a given accuracy while retaining the advantages of the BF-SPH method by providing only a vertical bias and by increasing or decreasing the resolution. The advantages of the vertical locus transformation were then further exploited by developing it into the $\sigma$-SPH method, which is innovative in the field of particle methods.\par

In the particle method, each particle is given a representative volume, and the material point with this constant representative volume moves with the motion. Using the $\sigma$-coordinates, when particles with the same volume move horizontally, a volume change is associated with the Jacobian, the rate of volume change associated with the coordinate transformation. 
To solve the volume change issue that arises when combining the particle method associated with the Lagrange description with the coordinate transformation technique using the $\sigma$-coordinates, a new volume conservation correction technique has been proposed, which was a significant factor in the success of this innovative $\sigma$-SPH method. The authors also point out that the computational accuracy of the second-order spatial derivative is important for coordinate transformation methods and that high-precision models such as SPH(2), a second-order spatial SPH proposed by the authors in a previous paper, should be adopted.\par

The usefulness of the above new computational techniques can be demonstrated through several computational examples, and in particular, the $\sigma$-SPH method has succeeded in improving the computational efficiency by a factor of five compared to BF-SPH, which has no spatial bias. \EDIT{As this study primarily confirms the feasibility of applying the proposed $\sigma$-SPH method to three-dimensional problems, future work will focus on its application to practical scenarios, such as coastal and marine environments, where its effectiveness is anticipated. These applications will serve to evaluate both the efficiency and accuracy of the method under more complex flow conditions.} \par
\EDITf{Moreover, in addition to enhancing computational efficiency, the proposed BFE-SPH and $\sigma$-SPH methods may provide additional advantages in resolving near-wall phenomena, including wall turbulence, particularly in turbulent flow simulations. This potential arises from their ability to flexibly control particle resolution in the vertical direction, which is advantageous for capturing steep velocity gradients near solid boundaries. Further investigation of this capability is planned for future work.}

\section*{CRediT authorship contribution statement}

\textbf{Shujiro Fujioka}: Writing - original draft, Validation, Software, Methodology, Investigation, Funding acquisition, Formal analysis, Conceptualization; 
\textbf{Kumpei Tsuji}: Writing - review \& editing, Software, Funding acquisition;
\textbf{Naoto Mitsume}: Writing - review \& editing, Methodology, Funding acquisition, Conceptualization;
\textbf{Mitsuteru Asai}: Writing - review \& editing, Supervision, Resources, Project administration, Funding acquisition, Conceptualization.

\section*{Acknowledgement}
This work was supported by JST SPRING~[Japan Grant Number JPMJSP2136]; JSPS~[KAKENHI Grant Number JP-23KK0182, 23K26356, 23K19132, 23H00160, 24K22288, 25KJ1965]; and SECOM Science and Technology Foundation. 

%% The Appendices part is started with the command \appendix;
%% appendix sections are then done as normal sections
% \appendix

%% If you have bib database file and want bibtex to generate the
%% bibitems, please use
%%
\bibliographystyle{elsarticle-num}
% \bibliography{ref_short.bib, references.bib}
\bibliography{ref_long.bib, references.bib}

\end{document}